\documentclass[12pt,a4paper]{article}

\usepackage{latexsym}
\usepackage{amsfonts}
\usepackage{amssymb}
\usepackage{amsmath}
\usepackage{theorem}
\usepackage{enumerate}
\usepackage{cite}
\usepackage{url}
\usepackage[dvips]{color}
\usepackage{graphicx}

\theoremstyle{plain}
\theorembodyfont{\itshape}
\newtheorem{thm}{Theorem}

\newtheorem{prp}{Proposition}

\newtheorem{defi}{Definition}

\theorembodyfont{\upshape}
\newtheorem{exem}{Example}

\newtheorem{rmk}{Remark}

\newcommand{\proof}{\noindent {\bf Proof:} \hspace{0.1in}}
\newcommand{\qed}{\hfill\mbox{\raggedright $\Box$}\medskip}
\newcommand{\smin}{\,\raisebox{0.06em}{${\scriptstyle \in}$}\,}
\newcommand{\nsmin}{\,\raisebox{0.06em}{${\scriptstyle \notin}$}\,}
\newcommand{\ssmin}{\,\raisebox{0.06em}{${\scriptscriptstyle \in}$}\,}
\newcommand{\smsubset}{\,\raisebox{0.06em}{${\scriptstyle \subset}$}\,}

\newcommand{\smcap}{\,\raisebox{0.06em}{${\scriptstyle \cap}$}\,}

\newcommand{\smotimes}{\,{\scriptstyle \otimes}\,}

\newcommand{\smwedge}{{\scriptstyle \wedge\,}}

\newcommand{\kwedge}{\ensuremath{%
 \,\raisebox{-0.5ex}{$\stackrel{{\scriptstyle \wedge}}{,}$}\,}}
\newcommand{\smvee}{{\scriptstyle \vee}}
\newcommand{\smcirc}{{\scriptstyle \,\circ\,}}

\newcommand{\bwedge}{\raisebox{0.2ex}{${\textstyle \bigwedge}$}}
\newcommand{\bvee}{\raisebox{0.2ex}{${\textstyle \bigvee}$}}

\newcommand{\fp}{\hat{\omega}}
\newcommand{\DA}{\ensuremath{D_{\!A\,}}^{}}

\begin{document}

\title{Multisymplectic and Polysymplectic Structures on Fiber Bundles
       \thanks{Work partially supported by CAPES (Coordena\c{c}\~ao de
               Aperfei\c{c}oamento de Pessoal de N\'{\i}vel Superior)
               and CNPq (Conselho Nacional de Desenvolvimento
               Cient\'{\i}fico e Tecno\-l\'o\-gico), Brazil}}
\author{Michael Forger~~and~~Leandro Gomes
        \thanks{Work done in partial fulfillment of the requirements
                for the degree of Doctor in Science}}
\date{Departamento de Matem\'atica Aplicada \\
      Instituto de Matem\'atica e Estat\'{\i}stica \\
      Universidade de S\~ao Paulo \\[2mm]
      BR--05315-970~ S\~ao Paulo, S.P., Brazil}
\maketitle

\thispagestyle{empty}

\begin{abstract}
\noindent
 We introduce the concepts of a multisymplectic structure and a polysymplectic
 structure on a general fiber bundle over a general base manifold, define the
 concept of the symbol of a multisymplectic form, which is a polysymplectic
 form representing its leading order contribution, and prove Darboux theorems
 for the existence of canonical local coordinates.
\end{abstract}

\begin{flushright}
 \parbox{12em}{
  \begin{center}
   Universidade de S\~ao Paulo \\
   RT-MAP-0702 \\
   December 2007
  \end{center}
 }
\end{flushright}

\newpage

\setcounter{page}{1}

\section*{Introduction}

Multisymplectic geometry is increasingly recognized as providing the
appropriate mathe\-matical framework for classical field theory from
the hamiltonian point of view~-- just as symplectic geometry does
for classical mechanics. Unfortunately, the development of this new
area of differential geometry has for a long time been hampered by
the lack of a fully satisfactory definition of the concept of a
multisymplectic structure, which should be mathematically simple
as well as in harmony with the needs of applications to physics;
the same goes for the closely related notion of a polysymplectic
structure.

The main goal of this paper is to provide such a definition and
establish a general relation between the two types of structure.

To set the stage, let us consider a simple analogy. The symplectic forms
encountered in classical mechanics can locally all be written in the form
\begin{equation} \label{eq:SPLFO1}
 \omega~=~dq_{}^i \>\smwedge\, dp_i^{}~,
\end{equation}
where $\; q_{}^1,\ldots,q_{}^N,p_1^{},\ldots,p_N^{} \;$ are a particular
kind of local coordinates on phase space known as canonical coordinates
or Darboux coordinates. Introducing time~$t$ and energy~$E$ as additional
variables (which is essential, e.g., for incorporating non-autonomous
systems into the symplectic framework of hamiltonian mechanics), this
equation is replaced by
\begin{equation} \label{eq:SPLFO2}
 \omega~=~dq_{}^i \>\smwedge\, dp_i^{} \, + \, dE \>\smwedge\, dt~,
\end{equation}
where $\; t,q_{}^1,\ldots,q_{}^N,p_1^{},\ldots,p_N^{},E \;$ can be viewed
as canonical coordinates on an extended phase space. Similarly, the
multisymplectic forms encountered in classical field theory over an
$n$-dimensional space-time manifold $M$ can locally all be written
in the form
\begin{equation} \label{eq:MSPF01}
 \omega~=~dq_{}^i \>\smwedge\, dp\>\!_i^\mu \>\smwedge\, d^{\,n} x_\mu^{} \, - \,
          dp \;\smwedge\, d^{\,n} x~,
\end{equation}
where $\; x_{}^\mu,q_{}^i,p\>\!_i^\mu,p \,$ ($1 \leqslant \mu \leqslant n ,
1 \leqslant i \leqslant N$) can again be viewed as canonical coordinates
on some extended multiphase space. Here, the $x_{}^\mu$ are (local)
coordinates for~$M$, while $p\,$ is still a single energy variable
(except for a sign), $d^{\,n} x$ is the (local) volume form induced
by the $x_{}^\mu$ and $d^{\,n} x_\mu^{}$ is the (local) $(n\!-\!1)$-form
obtained by contracting $d^{\,n} x$ with $\, \partial_\mu^{} \equiv
\partial/\partial x_{}^\mu$:
\[
 d^{\,n} x_\mu^{}~=~\mathrm{i}_{\partial_\mu^{}}^{} d^{\,n} x~.
\]
The idea of introducing ``multimomentum variables'' labelled by an additional
space-time index $\mu$ ($n$ multimomentum variables $p\>\!_i^\mu$ for each
position variable $q_{}^i$) goes back to the work of de Donder~\cite{DD}
and Weyl~\cite{We} in the 1930's (and perhaps even further) and has been
recognized ever since as being an essential and unavoidable ingredient in
any approach to a generally covariant hamiltonian formulation of classical
field theory. Understanding the proper geometric setting for this kind of
structure, however, has baffled both mathematicians and physicists for
decades, as witnessed by the large number of different proposals for
an appropriate global framework that can be found in the literature.

The by now standard example of a globally defined multisymplectic
structure starts out from an arbitrary fiber bundle $E$ over~$M$
called the configuration bundle (since its sections are the basic
fields of the field theory under consideration) and whose typical fiber
is an $N$-dimensional manifold $Q$ representing the configuration space,
as in mechanics. Following Refs~\cite{GIM,CCI}, for example, consider the
vector bundle $\, \bwedge_{\,1}^{\,n} \, T^* E \,$ of $n$-forms over~$E$
that are $(n-1)$-horizontal (i.e., that vanish whenever contracted with
more than $1$ vertical vector field), with projection onto~$E$ denoted
by $\pi_1^n$: its total space carries a naturally defined $n$-form
$\theta$ which we shall refer to as the multicanonical form, given by
\begin{equation} \label{eq:MCANF1}
 \theta_\alpha(v_1, \ldots, v_n)~
 =~\alpha(T_\alpha^{} \pi_1^n \cdot v_1, \ldots, T_\alpha^{} \pi_1^n \cdot v_n)
\end{equation}
for $\, \alpha \smin \bwedge_{\,1}^{\,n} \, T^* E \,$ and $\; v_1,\ldots,v_n
\smin T_\alpha^{} (\bwedge_{\,1}^{\,n} \, T^* E)$, and which gives rise to a
closed $(n+1)$-form $\, \omega = - d\theta \,$: this is the multisymplectic
form for the hamiltonian formalism of classical field theory.
The construction can be easily extended to the more general situation
of the vector bundle $\, \bwedge_{\,r-1}^{\,~k} \, T^* E \,$ of $k$-forms
over~$E$ that are $(k+1-r)$-horizontal (i.e., that vanish whenever con%
tracted with more than $r-1$ vertical vector fields), with projection
onto~$E$ denoted by $\pi_{r-1}^{~k}$, where $\, 1 \leqslant r \leqslant k \,$
and $\, k+1-r \leqslant n \,$: its total space carries a naturally defined
$k$-form $\theta$ which we shall again refer to as the multicanonical form,
given by
\begin{equation} \label{eq:MCANF2}
 \theta_\alpha(v_1, \ldots, v_k)~
 =~\alpha(T_\alpha^{} \pi_{r-1}^{~k} \cdot v_1, \ldots,
          T_\alpha^{} \pi_{r-1}^{~k} \cdot v_k)
\end{equation}
for $\, \alpha \smin \bwedge_{\,r-1}^{\,~k} \, T^* E \,$ and $\; v_1,\ldots,v_k
\smin T_\alpha^{} (\bwedge_{\,r-1}^{\,~k} \, T^* E)$, and which gives rise to
a closed $(k+1)$-form $\, \omega = - d\theta \,$: we propose to call it a
multilagrangian form. This term is motivated by the observation that the
vertical bundle for the projection $\pi_{r-1}^{~k}$ onto~$E$ is a distinguished
lagrangian (maximal isotropic) subbundle of the tangent bundle of
$\, \bwedge_{\,r-1}^{\,~k} \, T^* E$, and the existence of such a
``multilagrangian subbundle'' will turn out to play a central role in
our general definition of multisymplectic and multilagrangian structures.

The standard local coordinate expressions can be obtained by starting out
from local coordinates $(x^\mu,q^i)$ for~$E$ composed of local coordinates
$x^\mu$ for~$M$ and local coordinates $q^i$ for~$Q$ together with a local
trivialization of~$E$ over~$M$: these give rise to canonical local co%
ordinates $(x^\mu,q^i,p\>\!_i^\mu,p\>\!)$ for~$\, \bwedge_{\,1}^{\,n} \,
T^* E \,$ in which
\begin{equation} \label{eq:MCANF3}
 \theta~=~p\>\!_i^\mu \; dq_{}^i \>\smwedge\, d^{\,n} x_\mu^{} \, + \,
          p \; d^{\,n} x~,
\end{equation}
so
\begin{equation} \label{eq:MSPLF1}
 \omega~=~dq_{}^i \>\smwedge\, dp\>\!_i^\mu \>\smwedge\, d^{\,n} x_\mu^{} \, - \,
          dp \;\smwedge\, d^{\,n} x~,
\end{equation}
and more generally, to canonical local coordinates $(x^\mu,q^i,
p\>\!_{i_1 \ldots\, i_s;\,\mu_1 \ldots\, \mu_{k-s}}^{})$ ($0 \leqslant
s \leqslant r-1$) for~$\, \bwedge_{\,r-1}^{\,~k} \, T^* E \,$ in which
\begin{equation} \label{eq:MCANF4}
 \theta~=~\sum_{s=0}^{r-1} \, {\textstyle \frac{1}{s!} \, \frac{1}{(k-s)!}}~
          p\>\!_{i_1 \ldots\, i_s;\,\mu_1 \ldots\, \mu_{k-s}}^{} \,
          dq^{i_1} \,\smwedge \ldots \smwedge\, dq^{i_s} \,\smwedge\,
          dx^{\mu_1} \,\smwedge \ldots \smwedge\, dx^{\mu_{k-s}}~,
\end{equation}
so
\begin{equation} \label{eq:MLAGF1}
 \omega~= \; - \, \sum_{s=0}^{r-1} \,
          {\textstyle \frac{1}{s!} \, \frac{1}{(k-s)!}}~
          dp\>\!_{i_1 \ldots\, i_s;\,\mu_1 \ldots\, \mu_{k-s}}^{} \,\smwedge\,
          dq^{i_1} \,\smwedge \ldots \smwedge\, dq^{i_s} \,\smwedge\,
          dx^{\mu_1} \,\smwedge \ldots \smwedge\, dx^{\mu_{k-s}}~,
\end{equation}
In a more general context, when $\, \bwedge_{\,r-1}^{\,~k} \, T^* E \,$ is
replaced by a manifold~$P$ which is only supposed to be the total space of
a fiber bundle over a base manifold~$M$, Darboux's theorem guarantees the
existence of canonical local coordinates in which $\omega$ is given by the
expression in equation~(\ref{eq:MLAGF1}), under appropriate conditions on
the form $\omega$. The central question is to figure out what precisely
are these conditions.

A naive first guess would be to simply require the form $\omega$ to be
closed and non-degenerate. However, unlike in the symplectic case, these
conditions alone are far too weak to guarantee the validity of a Darboux
theorem, even at the purely algebraic level. For certain purposes, they
may be sufficient to derive results that are of interest (for an example,
see Refs~\cite{FR1,FPR1,FPR2}), but this version of the definition of a
multisymplectic structure~-- even though often adopted in the literature,
mostly for lack of a better one~-- is clearly inadequate. What is needed
is an additional algebraic condition.

An indication of what should be this additional algebraic condition can
be found in Ref.~\cite{Mar}, but the Darboux theorem proved there covers
a special situation which is disjoint from the case of interest for the
applications to physics because the structure of the underlying manifold
as the total space of a fiber bundle over space-time and the corresponding
horizontality conditions are completely ignored. More specifically,
Ref.~\cite{Mar} deals with what we call a multilagrangian form $\omega$
on a manifold, viewed as the total space of a fiber bundle whose base
manifold~$M$ is reduced to a point, so that the pertinent horizontality
condition becomes empty and the corresponding expression for~$\omega$
in canonical local coordinates takes the following form, analogous to
equation~(\ref{eq:MLAGF1}):
\begin{equation} \label{eq:MLAGF2}
 \omega~= \; - \, {\textstyle \frac{1}{k!}}~
          dp\>\!_{i_1 \ldots\, i_k}^{} \,\smwedge\,
          dq^{i_1} \,\smwedge \ldots \smwedge\, dq^{i_k}~.
\end{equation}
It is true that this generalizes the concept of a symplectic form to forms
of higher degree, but in a direction that is far away from the concept of
a multisymplectic form as encountered in classical field theory, and
this discrepancy, which was clearly stated only much later~\cite{CIL}
(see also Ref.~\cite{LDS}, for example), has created a great deal of
confusion in the literature. Taking into account that the term
``multisymplectic'' is already occupied at least since the mid
1970's~\cite{Kij,KS1,KS2}, we believe its use in a quite different
context such as that of Refs~\cite{Mar,CIL} to be inappropriate
and propose to correct this historical misnomer, replacing the term
``multisymplectic'' in that context by the term ``multilagrangian''.

In this paper, we shall follow a different approach, based on a new
and more profound understanding of the link between multisymplectic
and polysymplectic structures.

Polysymplectic structures in the hamiltonian approach to classical field
theory seem to have been introduced in Ref.~\cite{Gu} and have been further
investigated in Ref.~\cite{Aw} (where they were called ``$k$-symplectic
structures''~-- a terminology that we shall not follow in order not to
increase the already existing confusion). Roughly speaking, polysymplectic
forms are vector-valued analogues of symplectic forms. Similarly, the
polylagrangian forms to be introduced in this paper are vector-valued
analogues of the forms studied in Ref.~\cite{Mar}.

The standard example of a globally defined polysymplectic structure is
the one on the bundle $\, T^* E \otimes \hat{T} \,$ of $\hat{T}$-valued
$1$-forms over a manifold~$E$, with projection onto~$E$ denoted by $\pi^1$,
where $\hat{T}$ is a fixed finite-dimensional auxiliary vector space: its
total space carries a naturally defined $\hat{T}$-valued $1$-form $\theta$
which we shall refer to as the polycanonical form, given by
\begin{equation} \label{eq:PCANF1}
 \theta_\alpha(v)~=~\alpha(T_\alpha \pi^1 \cdot v)
\end{equation}
for $\, \alpha \smin T^* E \otimes \hat{T} \,$ and $\; v \smin T_\alpha^{}
(T^* E \otimes \hat{T})$, and which gives rise to a closed $\hat{T}$-valued
$2$-form $\, \omega = - d\theta \,$: this is the type of form called
polysymplectic, for instance, in Ref.~\cite{Gu}. Again, the construction
can be easily extended to the more general situation of the vector bundle
$\, \bwedge^k T^* E \otimes \hat{T} \,$ of $\hat{T}$-valued $k$-forms
over~$E$, with projection onto~$E$ denoted by $\pi^k$: its total space
carries a naturally defined $\hat{T}$-valued $k$-form $\theta$ which we
shall again refer to as the polycanonical form, given by
\begin{equation} \label{eq:PCANF2}
 \theta_\alpha(v_1, \ldots, v_k)~
 =~\alpha(T_\alpha \pi^k \cdot v_1, \ldots, T_\alpha^{} \pi^k \cdot v_k)
\end{equation}
for $\, \alpha \smin \bwedge^k T^* E \otimes \hat{T} \,$ and
$\; v_1,\ldots,v_k \smin T_\alpha (\bwedge^k T^* E \otimes \hat{T})$,
and which gives rise to a closed $\hat{T}$-valued $(k+1)$-form
$\, \omega = - d\theta \,$: we propose to call it a polylagrangian
form.  This term is motivated by the observation that the vertical
bundle for the projection $\pi^k$ onto~$E$ is a distinguished
lagrangian (maximal isotropic) subbundle of the tangent bundle of
$\, \bwedge^k T^* E \otimes \hat{T}$, and the existence of such a
``polylagrangian subbundle'' will turn out to play a central role
in our general definition of polysymplectic and polylagrangian
structures.

In terms of standard local coordinates $(q^i,p\>\!_i)$ for
$T^* E$ and a basis $\, \{ \hat{e}_a \, | \, 1 \leqslant a \leqslant
\hat{n} \}$ of~$\hat{T}$, we have
\begin{equation} \label{eq:PCANF3}
 \hat{\theta}~=~p\>\!_i^a \, dq^i \smotimes \hat{e}_a~,
\end{equation}
so
\begin{equation} \label{eq:PSPLF1}
 \fp~=~dq^i \,\smwedge\, dp\>\!_i^a \smotimes \hat{e}_a~,
\end{equation}
in the polysymplectic case and
\begin{equation} \label{eq:PCANF4}
 \hat{\theta}~=~{\textstyle \frac{1}{k!}}~
                p\>\!_{i_1 \ldots\, i_k}^a \,
                dq^{i_1} \,\smwedge \ldots \smwedge\, dq^{i_k} \,\smotimes\,
                \hat{e}_a~,
\end{equation}
so
\begin{equation} \label{eq:PLAGF1}
 \fp~= \; - \, {\textstyle \frac{1}{k!}}~
       dp\>\!_{i_1 \ldots\, i_k}^a \,\smwedge\,
       dq^{i_1} \,\smwedge \ldots \smwedge\, dq^{i_k} \,\smotimes\,
       \hat{e}_a~,
\end{equation}
in the polylagrangian case. In a more general context, when $\, \bwedge^k
T^* E \otimes \hat{T} \,$ is replaced by a manifold~$P$, Darboux's theorem
guarantees the existence of canonical local coordinates in which $\fp$
is given by the expression in equation~(\ref{eq:PLAGF1}), under appropriate
conditions on the form~$\fp$. Again, the central question is to figure
out what precisely are these conditions.

It should be pointed out that for the polysymplectic case, this problem
has been solved in Ref.~\cite{Aw}, but the fundamental role of what we
call the polylagrangian subbundle is not fully appreciated there. As it
turns out, this object and its basic properties are the key to the
entire subject, allowing to generalize the proof of Darboux's theorem
not only from the polysymplectic to the polylagrangian case but also
from polysymplectic/polylagrangian structures on manifolds to poly%
symplectic/polylagrangian structures on the total spaces of fiber
bundles~-- a concept that conveys a precise mathematical meaning
to the idea of a ``smooth family of polysymplectic/polylagrangian
structures'' (each fiber is a polysymplectic/polylagrangian manifold
in such a way that the entire structure depends smoothly on the points
of the base manifold). This extension is particularly useful in that it
finally allows to formulate in precise mathematical terms what is the
relation between multisymplectic and polysymplectic structures and,
more generally, between multilagrangian and polylagrangian structures:
the latter appear as the ``leading order term'' of the former, through
a simple formal construction that we introduce below and propose to
call the ``symbol'' because it strongly resembles the construction
of the symbol of a differential operator.


\section{Polylagrangian forms on vector spaces}
\label{sec:FPLA}

We begin by briefly recalling a few basic notions involving vector-valued
alternating multi\-linear forms. Given finite-dimensional real vector spaces%
\footnote{In order to simplify the presentation, we assume all vector spaces
involved to be real and finite-dimensional: the extension to vector spaces
over an arbitrary field of characteristic $0$ is straightforward, and
generalization to the infinite-dimensional setting, which requires
imposing appropriate continuity conditions from functional analysis,
will be left to a possible future investigation.}
$V$ and $\hat{T}$, we consider $\hat{T}$-valued $(k+1)$-forms $\,\fp$ on~$V$,
\begin{equation}
 \fp \,\smin\, \bwedge^{k+1} \, V^* \otimes\, \hat{T}~.
\end{equation}
The \textbf{contraction} of such a form~$\,\fp$ is the linear map
$\; \fp^\flat : V \rightarrow \bwedge^k \, V^* \otimes\, \hat{T} \;$
given by
\begin{equation} \label{eq:OMFLT1} 
 \fp^\flat(v)~=~\mathrm{i}_v^{} \fp~,
\end{equation}
and the \textbf{kernel} of $\,\fp$ is defined to be the kernel of
$\,\fp^\flat$: $\ker \, \fp = \ker \, \fp^\flat$. If $\, \ker \, \fp
= \{0\}$ we say that $\,\fp$ is \textbf{non-degenerate}. Given a
linear form $\, \hat{t}^* \smin\, \hat{T}^* \,$ on $\hat{T}$, the
\textbf{projection} of $\,\fp$ along $\hat{t}^*$ is the ordinary
$(k+1)$-form on~$V$ given by%
\footnote{Throughout this paper the symbol $\langle .\,,. \rangle$ will
stand for the natural bilinear pairing between a vector space and its dual.}
\begin{equation} \label{eq:FPROJ1} 
 \omega_{\hat{t}^*}~=~\langle \hat{t}^* , \fp \rangle~.
\end{equation}
Note that $\,\omega_{\hat{t}^*}^{}$ depends linearly on~$\hat{t}^*$, so if we
choose a basis $\, \{ \hat{e}_1,\ldots,\hat{e}_{\hat{n}} \} \,$ of~$\hat{T}$,
with dual basis $\, \{ \hat{e}^1,\ldots,\hat{e}^{\hat{n}} \} \,$ of~%
$\hat{T}^*$, we have
\begin{equation} \label{eq:FPROJ2} 
 \fp~=~\omega^a \otimes \hat{e}_a^{} \quad \text{with}\quad
 \omega^a~=~\omega_{\hat{e}^a}^{} \quad (1 \leqslant a \leqslant \hat{n})~.
\end{equation}
Then it is clear that
\begin{equation} \label{eq:FKERN3} 
 \ker \, \fp~=~\bigcap_{\hat{t}^* \!\ssmin \hat{T}^*}
               \ker \, \omega_{\hat{t}^*}^{}~
             =~\bigcap_{a=1}^{\hat{n}} \; \ker \, \omega^a~.
\end{equation}
Next, suppose that $L$ is a subspace of~$V$ and $\ell$ is an integer
satisfying $\, 1 \leqslant \ell \leqslant k \,$; then extending the
definition given in Ref.~\cite{CIL} from ordinary to vector-valued
forms, we define the \textbf{$\ell$-orthogonal complement} of~$L$
(with respect to~$\,\fp$) to be the subspace $L^{\fp,\ell}$ of~$V$
given by%
\footnote{We discard the trivial case $\, \ell=0$ since extrapolating
the definition to this case would lead to the conclusion that $L^{\fp,0}$
is simply the kernel of~$\,\fp$, independently of the subspace $L$ of~$V$.}
\begin{equation}
 L^{\fp,\ell}~=~\{ \, v \smin V~|~\mathrm{i}_v^{} \mathrm{i}_{v_1}^{} \ldots
                                 \mathrm{i}_{v_{\ell}}^{} \fp = 0~~
                     \mbox{for all $\, v_1,\ldots,v_{\ell} \smin L$} \, \}~.
\end{equation}
Note that these orthogonal complements form an increasing sequence under
inclusion:
\begin{equation}
 L^{\fp,1} \, \subset \ldots \subset \, L^{\fp,k}~.
\end{equation}
We say that $L$ is \textbf{$\ell$-isotropic} (with respect to~$\,\fp$)
if $\, L \subset L^{\fp,\ell} \,$ and is \textbf{maximal $\ell$-isotropic}
or, more briefly, \textbf{$\ell$-lagrangian} (with respect to~$\,\fp$)
if it is maximal in the partially ordered set formed by the $\ell$-%
isotropic subspaces of~$V$; it is a simple exercise to check that,
as usual, this is the case if and only if $\, L = L^{\fp,\ell}$.
If $\ell\!=\!1$ we omit the prefix $1$ and may conclude that a
subspace $L$ of~$V$ will be isotropic if and only if\,%
\footnote{Throughout this paper the symbol $.^\bot$ will denote the
annihilator of a subspace, i.e., given a subspace $L$ of a vector
space~$V$, $L^\bot$ is the subspace of its dual space $V^*$ consisting
of the linear forms on~$V$ that vanish on~$L$.}
\begin{equation} \label{eq:ISOSS2} 
 \fp^\flat \bigl( L \bigr)~
 \subset~\fp^\flat \bigl( V \bigr) \,\cap\,
         \bigl( \bwedge^k L^\bot \bigr) \otimes\, \hat{T}
\end{equation}
and will be maximal isotropic if and only if $\, \ker \, \fp \subset L \,$ and
\begin{equation} \label{eq:ISOSS3} 
 \fp^\flat \bigl( L \bigr)~
 =~\fp^\flat \bigl( V \bigr) \,\cap\,
   \bigl( \bwedge^k L^\bot \bigr) \otimes\, \hat{T}~.
\end{equation}
(For explicit proofs of these elementary statements, the reader may
consult Ref.~\cite{Go}.) At~first~sight, the intersection with the
subspace $\,\fp^\flat \bigl( V \bigr)$ on the rhs of these relations
may seem strange, in particular since in the inclusion stated in
equation~(\ref{eq:ISOSS2}) it is really superfluous, but that is by no
means the case for the equality stated in equation~(\ref{eq:ISOSS3}).
Rather, omitting this intersection leads to a strengthened form
of equation~(\ref{eq:ISOSS3}) which turns out to provide the key to the
theory of polysymplectic and, more generally, of polylagrangian forms:
\begin{defi}~ \label{def:POLI}
 Let\/ $V$ and $\hat{T}$ be finite-dimensional vector spaces, with
 $\, \hat{n} \equiv \dim \hat{T}$, and let $\,\fp$ be a non-vanishing
 $\hat{T}$-valued $(k+1)$-form on~$V$. We say that $\,\fp$ is a
 \textbf{polylagrangian} form of \textbf{rank}~$N$ if\/ $V$ admits
 a subspace $L$ of codimension~$N$ which is \textbf{polylagrangian},
 i.e., such that
 \begin{equation} \label{eq:POLIF} 
  \fp^\flat \bigl( L \bigr)~
  =~\bigl( \bwedge^k L^\bot \bigr) \otimes\, \hat{T}~.
 \end{equation}
 When $\, k=1 \,$ and $\,\fp$ is non-degenerate, we call $\,\fp$ a
 \textbf{polysymplectic} form. If the condition of non-degeneracy
 is dropped, we speak of a \textbf{polypresymplectic} form.
\end{defi}
As a first property of polylagrangian forms, we note that a polylagrangian
subspace, when it exists, contains the kernel of $\,\fp$ and hence really
is a special type of maximal isotropic subspace. But more than that is true.
\begin{prp}~ \label{prp:POLILS1} 
 Let\/ $V$ and $\hat{T}$ be finite-dimensional vector spaces, with
 $\, \hat{n} \equiv \dim \hat{T}$, and let $\,\fp$ be a $\hat{T}$-valued
 polylagrangian $(k+1)$-form on\/~$V$ of rank~$N$, with polylagrangian
 subspace~$L$. Then $\, N \geqslant k$, and $L$ contains the kernel
 of~$\,\fp$ as well as the kernel of each of the projected forms
 $\,\omega_{\hat{t}^*}^{}$ $(\hat{t}^* \!\smin \hat{T}^* \setminus \{0\})$:
 \begin{equation} \label{eq:POLIK1}
  \ker \, \fp \,\subset\, L \qquad , \qquad
  \ker \, \omega_{\hat{t}^*}^{} \,\subset\, L \quad
  \mbox{$\mathrm{for~all}$ $\, \hat{t}^* \!\smin \hat{T}^* \setminus \{0\}$}~.
 \end{equation}
\end{prp}
\proof~
 First we observe that if $\, N<k$, we have $\, \bwedge^k L^\bot = \{0\}$, so
 both sides of the equation~(\ref{eq:POLIF}) vanish, i.e., $L$ is contained
 in $\, \ker \, \fp \,$ and hence $\, \ker \, \fp$ has codimension~$<k$
 in~$V$, implying $\, \fp \equiv 0 \,$, since the $(k+1)$-form on the
 quotient space $\, V/\ker \, \fp \,$ induced by $\fp$ vanishes identically.
 (More generally, this argument shows that a nonvanishing vector-valued
 $(k+1)$-form does not permit isotropic subspaces of codimension $<k$.)
 \linebreak
 Thus supposing that $\, \dim L^\bot = N \geqslant k$, we can for any
 vector $\, v \,\smin\, V \setminus L \,$ find a linearly independent
 set of $1$-forms $\, v_1^*,\ldots,v_k^* \smin L^\bot \,$ such that
 $\, \langle v_1^*,v \rangle = 1 \,$ and $\, \langle v_i^*, v \rangle = 0 \,$
 for $\, i > 1$. \linebreak
 Given $\, \hat{t}^* \!\smin \hat{T}^*$, take $\, \hat{t} \smin \hat{T} \,$
 such that $\, \langle\, \hat{t}^*,\hat{t} \,\rangle = 1$. According to the
 definition of a poly\-lagrangian subspace, there is a vector $\, u \smin L \,$
 such that
 \[
  \mathrm{i}_u^{} \fp~
  =~v_1^* \,\smwedge \ldots \smwedge\, v_k^* \otimes \hat{t}
  \quad \Rightarrow \quad
  \mathrm{i}_v^{} \mathrm{i}_u^{} \omega_{\hat{t}^*}^{}~
  =~v_2^* \,\smwedge \ldots \smwedge\, v_k^*~
  \ne~0
 \]
 and so $\, v \nsmin \ker \, \omega_{\hat{t}^*}^{}$. Hence it follows that
 $\, \ker \, \fp \smsubset \ker \, \omega_{\hat{t}^*}^{} \smsubset L$.
\qed

\noindent
On the other hand, considering the case of main interest, which is that
of $2$-forms, it must be emphasized that, as shown by the counterexamples
presented in Appendix~A, by far not every vector-valued $2$-form is poly%
(pre)symplectic, which means that in contrast to lagrangian subspaces, a
polylagrangian subspace need not exist, and even if it does exist, not every
lagrangian subspace is polylagrangian. However, there is a simple dimension
criterion that allows to decide whether a given isotropic subspace is
polylagrangian:
\begin{prp}~ \label{prp:POLILS2} 
 Let\/ $V$ and $\hat{T}$ be finite-dimensional vector spaces, with $\, \hat{n}
 \equiv \dim \hat{T}$, and let $\,\fp$ be a non-vanishing $\hat{T}$-valued
 $(k+1)$-form on~$V$. Given any subspace $L$ of\/~$V$, with
 $\, N = \dim (V/L)$, such that $\, N \geqslant k$, the
 following statements are equivalent:
 \begin{itemize}
  \item $L$ is a polylagrangian subspace and $\,\fp$ is a polylagrangian
        form of rank $N$.
  \item $L$ contains $\, \ker \, \fp$, is isotropic and has dimension
        \begin{equation} \label{eq:POLID2}
         \dim \, L~=~\dim \, \ker \, \fp \, + \, \hat{n} \, {N \choose k}~.
        \end{equation}
 \end{itemize}
\end{prp}
\proof~
 Taking into account that, for any isotropic subspace $L$ of~$V$ containing
 the kernel of $\,\fp$, the contraction map $\,\fp^\flat$ induces an injective
 linear map of $\, L/\ker \, \fp \,$ into $\, \bigl( \bwedge^k L^\bot \bigr)
 \otimes\, \hat{T}$, the result follows from an elementary dimension count.
\qed

\noindent
The case of non-degenerate ordinary or scalar polylagrangian forms
($\dim \hat{T} = 1$) has been studied in the literature~\cite{Mar,CIL,LDS}
under the label ``multisymplectic forms'' (a terminology that we propose
to abandon since we use this term in a different sense; see Definition~%
\ref{def:MULT} \linebreak below), but the concept of polylagrangian sub%
space appears only implicitly, namely through the dimension criterion
formulated in Proposition~\ref{prp:POLILS2} above, which is employed
as a definition, so that it remains unclear how to extend this purely
numerical recipe to other situations, in particular when $\,\fp$ is
taken to be a truly vector-valued form ($\dim \hat{T} > 1$).
The main statement here is
\begin{prp}~ \label{prp:POLILS3} 
 Let\/ $V$ be a finite-dimensional vector space and let $\,\fp$ be a scalar
 poly\-lagrangian $(k+1)$-form on\/~$V$ of rank~$N$, with polylagrangian
 subspace~$L$. Then any isotropic subspace $\tilde{L}$ containing the
 kernel of $\,\fp$ and such that
 \begin{equation} \label{eq:POLID3}
  \dim \tilde{L}~>~\dim \, \ker \, \fp + {N-1 \choose k} + 1
 \end{equation}
 is contained in~$L$. In particular, if $\, N > k > 1$, $L$ is unique.
\end{prp}
\begin{rmk}~
 Note that the uniqueness statement for $L$ is of course false for symplectic
 forms ($k=1$) and also for volume forms ($N=k$): in both cases, isotropic
 subspaces $\tilde{L}$ satisfying the dimension condition~(\ref{eq:POLID3})
 do not exist, and there is no restriction whatsoever on the relative
 position of lagrangian subspaces (which for a symplectic form on a
 $(2N)$-dimensional space are $N$-dimensional and for a volume form
 are one-dimensional).
\end{rmk}
\proof
 Obviously, passing from $V$ to the quotient space $\, V/\ker \, \fp \,$
 if necessary, and taking into account the previous remark, we may assume
 without loss of generality that $\,\fp$ is non-degenerate and also that
 $\, N > k > 1$. Following Refs~\cite{Mar,LDS}, we begin by showing
 that any isotropic subspace of~$V$ of dimension greater than~$1$ must
 intersect $L$ non-trivially. Indeed, if $v_1$ and $v_2$ are linearly
 independent vectors in~$V$ such that the two-dimensional subspace
 $\mathsf{span}(v_1,v_2)$ generated by $v_1$ and $v_2$ satisfies
 $\; \mathsf{span}(v_1,v_2) \cap L = \{0\}$, we can find a basis
 $\, \{ e_1^{},\ldots, e_N^{} \} \,$ of a subspace of~$V$ complementary
 to~$L$ such that $\, e_1 = v_1$ \linebreak and $\, e_2 = v_2 \,$; then
 denoting the corresponding dual basis of~$L^\bot$ by $\, \{ e^1,\ldots,
 e^N \}$, we use the fact that $L$ is polylagrangian to conclude that
 there exists a vector $u$ in $L$ such that $\, \mathrm{i}_u^{} \fp =
 e^1 \smwedge \ldots\, \smwedge e^k$, so $\; \fp(u,e_1,\ldots,e_k) = 1 \,$
 and hence $\, \mathrm{i}_{v_1} \mathrm{i}_{v_2} \fp \,$ cannot vanish, i.e.,
 $\mathsf{span}(v_1,v_2)$ cannot be isotropic. Using this result, we conclude
 that if $\tilde{L}$ is any isotropic subspace of~$V$, then the codimension
 of $\, \tilde{L} \cap L \,$ in $\tilde{L}$ is at most $1$, so
 \[
  \dim \bigl( \tilde{L} + L \bigr) \, - \, \dim L~
  =~\dim \tilde{L} \, - \, \dim \bigl( \tilde{L} \cap L \bigr)~
  \leqslant~1~.
 \]
 Suppose now that $\tilde{L}$ is an isotropic subspace of~$V$ such that
 \[
  \dim \bigl( \tilde{L} + L \bigr) \, - \, \dim L~
  =~\dim \tilde{L} \, - \, \dim \bigl( \tilde{L} \cap L \bigr)~
  =~1~.
 \]
 Then $\, \tilde{L} + L \,$ has codimension $N-1$ in $V$ and hence
 \[
  \dim \, \bwedge^k \bigl( \tilde{L}+L \bigr)^\bot~=~{N-1 \choose k}~,
 \]
 whereas, by hypothesis,
 \[
  \dim \bigl( \tilde{L} \cap L \bigr)~>~{N-1 \choose k}~,
 \]
 which contradicts the fact that the restriction of the linear isomorphism
 $\; \fp^\flat\vert_L: L \rightarrow \bwedge^k L^\bot$ \linebreak maps the
 subspace $\, \tilde{L} \cap L \,$ of~$L$ injectively into the subspace
 $\, \bwedge^k \bigl( \tilde{L}+L \bigr)^\bot \,$ of~$\bwedge^k L^\bot$.
\qed

\noindent
What is remarkable is that the above uniqueness statement remains valid
for vector-valued polylagrangian forms ($\dim \hat{T} > 1$) and that in
this case it becomes true even when $\, k=1$ \linebreak or $\, N=k \,$:
this is a consequence of the following explicit construction of $L$.
\begin{thm}~ \label{thm:POLILS1} 
 Let $V$ and $\hat{T}$ be finite-dimensional vector spaces, with $\, \hat{n}
 \equiv \dim \hat{T} \geqslant 2$, and let $\,\fp$ be a $\hat{T}$-valued
 polylagrangian $(k+1)$-form of rank~$N$ on~$V$, with polylagrangian
 subspace $L$. Then $L$ is given by
 \begin{equation} \label{eq:POLIK2}
 L~=~\sum_{t^* \!\ssmin\, \hat{T}^* \setminus \{0\}}
      \ker \, \omega_{\hat{t}^*}^{}
 \end{equation}
 and, in particular, is unique. In terms of a basis $\, \{ \hat{e}_1,
 \ldots,\hat{e}_{\hat{n}} \} \,$ of\/~$\hat{T}$ with dual basis
 $\, \{ \hat{e}^1,\ldots,\hat{e}^{\hat{n}} \} \,$ of\/~$\hat{T}^*$,
 \begin{equation} \label{eq:POLIK3}
  L~=~\ker \, \fp \, \oplus K_1^{} \oplus \ldots \oplus K_{\hat{n}}^{}~,
 \end{equation}
 and for $\, 1 \leqslant a \leqslant \hat{n}$,
 \begin{equation} \label{eq:POLIK4}
  L~=~\ker \, \omega^a \oplus K_a^{}~,
 \end{equation}
 where for $\, 1 \leqslant a \leqslant \hat{n}$, $K_a^{}$ is a subspace
 of\/~$V$ chosen so that
 \begin{equation} \label{eq:POLIK5}
  \bigcap_{b=1 \atop b \ne a}^{\hat{n}} \, \ker \, \omega^b~
  =~\ker \, \fp \oplus K_a^{}~.
 \end{equation}
 The dimensions of these various subspaces are given by
 \begin{equation} \label{eq:POLIK6}
  \dim \, \ker \, \omega^a~
  =~\dim \, \ker \, \fp \, + \, (\hat{n}-1) {N \choose k}~~~,~~~
  \dim K_a^{}~=~{N \choose k}~.
 \end{equation}
\end{thm}
\proof~
 Fix a basis $\, \{ \hat{e}_1,\ldots,\hat{e}_{\hat{n}} \} \,$ of~$\hat{T}$
 with dual basis $\, \{ \hat{e}^1,\ldots,\hat{e}^{\hat{n}} \} \,$ of~%
 $\hat{T}^*$ and choose subspaces $K_a^{}$ of~$V$ $(a = 1,\ldots,\hat{n})$
 as indicated above. Then the subspaces $\, \ker \, \fp \,$ and $\, K_1^{},
 \ldots,K_{\hat{n}}^{} \,$ of~$V$ have trivial intersection, so their sum
 is direct and defines a subspace of~$V$ which we shall, for the moment,
 denote by $L'$. According to the previous proposition, $L' \subset L$.
 To show that $\, L' = L$, it is therefore sufficient to prove that
 $\, \fp^\flat \bigl( L \bigr) \smsubset\, \fp^\flat \bigl( L' \bigr)$,
 since both $L$ and $L'$ contain $\, \ker \, \fp$. Using the definition
 of a polylagrangian subspace, we conclude that we must establish the
 inclusion
 \[
  \bigl( \bwedge^k L^\bot \bigr) \otimes\, \hat{T}~
  \subset~\fp^\flat \bigl( L' \bigr)~.
 \]
 But the equality~(\ref{eq:POLIF}) guarantees that for any $\, \alpha \smin
 \bwedge^k L^\bot \,$ and for $\, 1 \leqslant a \leqslant \hat{n}$, there is
 a vector $\, v_a^{} \smin L \,$ such that
 \[
  \mathrm{i}_{v_a}^{} \fp~=~\alpha \otimes \hat{e}_a^{}~.
 \]
 Then
 \[
  \mathrm{i}_{v_a}^{} \fp^b~=~\alpha \, \langle \hat{e}^b ,
  \hat{e}_a^{} \rangle~=~\delta_a^b \, \alpha~,
 \]
 so we see that
 \[
  v_a^{} \in \bigcap_{b=1 \atop b \ne a}^{\hat{n}} \, \ker \, \fp^b~.
 \]
 Decomposing $v_a^{}$ according to equation~(\ref{eq:POLIK5}), we find a
 vector $\, u_a^{} \smin K_a^{} \,$ such that
 \[
  \mathrm{i}_{u_a}^{} \fp~=~\alpha \otimes \hat{e}_a^{}~,
 \]
 so $\; \alpha \otimes \hat{e}_a^{} \,\smin\, \fp^\flat(K_a^{})
 \,\smsubset\, \fp^\flat(L')$. Conversely, if $\, u_a^{} \smin K_a^{}
 \smsubset L$, then $\, \mathrm{i}_{u_a}^{} \fp^b = 0 \,$ for $\, b
 \neq a \,$ and hence $\mathrm{i}_{u_a}^{} \fp$ is of the form
 $\, \alpha \otimes \hat{e}_a^{} \,$ for some $\, \alpha \smin
 \bwedge^k L^\bot$. Thus we conclude that $\,\fp^\flat$ maps $K_a$
 isomorphically onto $\, \bigl( \bwedge^k L^\bot \bigr) \otimes\,
 \mathsf{span}(\hat{e}_a^{})$, where $\mathsf{span}(\hat{e}_a^{})$
 denotes the one-dimensional subspace of~$\hat{T}$ generated by
 $\hat{e}_a^{}$, which proves the second formula in equation~%
 (\ref{eq:POLIK6}). Finally, we observe that combining equations~%
 (\ref{eq:FKERN3}) and~(\ref{eq:POLIK5}) gives
 \[
  \ker \, \omega^a \smcap K_a^{}~=~\{0\}~,
 \]
 so that equation~(\ref{eq:POLIK4}) follows from equation~(\ref{eq:POLIK3}),
 whereas the first formula in equation~(\ref{eq:POLIK6}) is now a direct
 consequence of equations~(\ref{eq:POLID2}) and~(\ref{eq:POLIK4}).
\qed

\noindent
Another fundamental property of polylagrangian forms is that the polylagrangian
subspace has a particular type of direct complement.
\begin{thm}~ \label{thm:POLIDD} 
 Let $V$ and $\hat{T}$ be finite-dimensional vector spaces and let $\,\fp$
 be a $\hat{T}$-valued polylagrangian $(k+1)$-form of rank~$N$ on~$V$, with
 polylagrangian subspace~$L$. Then there exists a $k$-isotropic subspace
 $E$ of\/~$V$ complementary to~$L$, i.e., such that
 \begin{equation} \label{eq:POLIDD}
  V~=~E \oplus L~.
 \end{equation}
\end{thm}
\proof~
 Let $E_0$ be a $k$-isotropic subspace of~$V$ of dimension $N'$ such that
 $\, E_0 \cap L = \{0\}$. (For instance, as long as $\, N' \leqslant k$, $E_0$
 can be any subspace of~$V$ such that $\, E_0 \cap L = \{0\}$.) \linebreak
 If $\, N'=N$, we are done. Otherwise, choose a basis $\, \{ e_1^{},\ldots,
 e_N^{} \} \,$ of a subspace of~$V$ complementary to~$L$ such that the first
 $N'$ vectors constitute a basis of~$E_0$, and denote the corresponding dual
 basis of~$L^\bot$ by $\, \{ e^1,\ldots,e^N \}$. We shall prove that there
 exists a vector $\, u \smin V \setminus (E_0 \oplus L) \,$ such that the
 subspace $E_1$ of~$V$ spanned by $u$ and $E_0$ is $k$-isotropic and
 satisfies $\, E_1 \cap L = \{0\}$; then since $\, \dim E_1 = N'+1$, the
 statement of the theorem follows \linebreak by induction. To this end,
 consider an arbitrary basis $\, \{ \hat{e}_a^{} \,|\; 1 \leqslant a
 \leqslant \hat{n} \} \,$ of~$\hat{T}$ with dual basis $\, \{ \hat{e}^a
 \,|\; 1 \leqslant a \leqslant \hat{n} \} \,$ of~$\hat{T}^*$ and, choosing
 any subspace $L'$ of~$L$ complementary to $\, \ker \, \fp$, use the fact
 that~$L$ is polylagrangian to conclude that there exists a unique basis
 $\, \{\, e_a^{i_1 \ldots\, i_k} \,|\; 1 \leqslant a \leqslant \hat{n} \,,\,
 1 \leqslant i_1 < \ldots <i_k \leqslant N \,\} \,$ of~$L'$ such that
 \[
  \fp^\flat (e_a^{i_1 \ldots\, i_k})~
  =~e^{i_1} \,\smwedge \ldots\, \smwedge\, e^{i_k} \,\smotimes\, \hat{e}_a~.
 \]
 Thus, for $\, 1 \leqslant i_1 < \ldots < i_k \leqslant N \,$ and
 $\, 1 \leqslant j_1 < \ldots <j_k \leqslant N$, we have
 \[
  \omega^b(e_a^{i_1 \ldots\, i_k},e_{j_1}^{},\ldots,e_{j_k}^{})~
  =~\delta_a^b \, \delta_{j_1}^{i_1} \ldots\, \delta_{j_k}^{i_k}~.
 \]
 Therefore, the vector
 \[
  u~=~e_{N'+1} \, - \, {\textstyle \frac{1}{k!}} \;
      \omega^a(e_{N'+1}^{},e_{i_1}^{},\ldots,e_{i_k}^{}) \,
      e_a^{i_1 \ldots\, i_k}
 \]
 does not belong to the subspace $\, E_0 \oplus L \,$ and, for
 $\, 1 \leqslant j_1 < \ldots <j_k \leqslant N$, satisfies
 \begin{eqnarray*}
  \lefteqn{\omega^b(u,e_{j_1}^{},\ldots,e_{j_k}^{})} \\
  &=&\!\! \omega^b(e_{N'+1}^{},e_{j_1}^{},\ldots,e_{j_k}^{}) \, - \,
          {\textstyle \frac{1}{k!}} \;
          \omega^a(e_{N'+1}^{},e_{i_1}^{},\ldots,e_{i_k}^{}) \;
          \omega^b (e_a^{i_1 \ldots\, i_k},e_{j_1}^{},\ldots,e_{j_k}^{})
  \\[1mm]
  &=&\!\! 0~,
 \end{eqnarray*}
 which implies that since the subspace $E_0$ spanned by $e_1^{},\ldots,
 e_{N'}^{}$ is $k$-isotropic, the subspace $E_1$ spanned by $e_1^{},\ldots,
 e_{N'}^{}$ and $u$ is so as well.
\qed

\begin{exem}(\textbf{The canonical form})~~
 Let $E$ and $\hat{T}$ be vector spaces of dimension $N$ and~$\hat{n}$,
 respectively. Set
 \begin{equation} \label{eq:FPLC1} 
  V_0^{}~=~E \,\oplus\, \left( \bigl( \bwedge^k E^* \bigr)
                               \otimes\, \hat{T} \right)~.
 \end{equation}
 The \textbf{canonical polylagrangian form} of rank $N$ is the non-degenerate
 $\hat{T}$-valued $(k+1)$-form $\,\fp_0^{}$ on~$V_0^{}$ defined by
 \begin{equation} \label{eq:FPLC2} 
  \fp_0^{} \bigl( (u_0^{},\alpha_0^{} \smotimes \hat{t}_0^{}) , \ldots ,
                  (u_k^{},\alpha_k^{} \smotimes \hat{t}_k^{}) \bigr)~
  =~\sum_{i=0}^k \, (-1)^i \,
    \alpha_i^{} (u_0^{},\ldots,\widehat{u_i^{}},\ldots,u_k^{}) \;
    \hat{t}_i^{}~.
 \end{equation}
 If $\, k=1 \,$ we shall call $\,\fp_0^{}$ the \textbf{canonical
 polysymplectic form}.

 \noindent
 For scalar forms ($\hat{T} = \mathbb{R}$) this construction can be
 found, e.g., in Refs~\cite{Mar,CIL,LDS}.
\end{exem}
To justify this terminology, note that it is a straightforward exercise
to show that $\fp_0^{}$ is non-degenerate and that, considering $E$ and
\begin{equation} \label{eq:FPLC3} 
 L~=~\bigl( \bwedge^k E^* \bigr) \otimes\, \hat{T}
\end{equation}
as subspaces of~$V_0^{}$, we have the direct decomposition $\,
V_0^{} = E \oplus L \,$ where
\begin{equation} \label{eq:FPLC4} 
 L~\text{is polylagrangian}
 \quad \text{and} \quad
 E~\text{is $k$-isotropic}~.
\end{equation}
In terms of bases, let $\, \{\, \hat{e}_a^{} \,|\; 1 \leqslant a \leqslant
\hat{n} \,\} \,$ be any basis of~$\hat{T}$ with dual basis $\, \{\, \hat{e}^a
\,|\; 1 \leqslant a \leqslant \hat{n} \,\}$ \linebreak of $\hat{T}^*$ and
let $\, \{\, e_i^{} \,|\; 1 \leqslant i \leqslant N \,\} \,$ be any basis
of~$E$ with dual basis $\, \{\, e^i \,|\; 1 \leqslant i \leqslant N \,\} \,$
of~$E^*$. \linebreak For $\, 1 \leqslant a \leqslant \hat{n} \,$ and
$\; 1 \leqslant i_1^{} < \ldots < i_k^{} \leqslant N$, define
\[
 e_a^{i_1 \ldots\, i_k}~
 =~e^{i_1} \smwedge \ldots\, \smwedge\, e^{i_k} \smotimes \hat{e}_a^{}~~~,~~~
 e_{i_1 \ldots\, i_k}^a~
 =~e_{i_1}^{} \smwedge \ldots\, \smwedge\, e_{i_k}^{} \smotimes \hat{e}^a~.
\]
This provides a basis $\, \{\, e_i^{} , e_a^{i_1 \ldots\, i_k} \,|\;
1 \leqslant a \leqslant \hat{n} \,,\, 1 \leqslant i \leqslant N \,,\,
1 \leqslant i_1^{} < \ldots < i_k^{} \leqslant N \,\} \,$ of~$V_0^{}$
with dual basis $\, \{\, e^i , e_{i_1 \ldots\, i_k}^a \,|\;
1 \leqslant a \leqslant \hat{n} \,,\, 1 \leqslant i \leqslant N \,,\,
1 \leqslant i_1^{} < \ldots < i_k^{} \leqslant N \,\} \,$ of~$V_0^*$,
both of which we shall refer to as a \textbf{canonical basis} or
\textbf{Darboux basis}, such that
\begin{equation} \label{eq:FPLCB1} 
 \fp_0^{}~=~{\textstyle \frac{1}{k!}} \;
            \bigl( e_{i_1 \ldots\, i_k}^a \,\smwedge\,
                   e^{i_1} \,\smwedge \ldots\, \smwedge\, e^{i_k} \bigr)
            \,\smotimes\, \hat{e}_a~.
\end{equation}

Now it is easy to derive the algebraic Darboux theorem for general poly%
lagrangian forms: let $\, \{\, \hat{e}_a^{} \,|\; 1 \leqslant a \leqslant
\hat{n} \,\} \,$ be an arbitrary basis of~$\hat{T}$, with dual basis
$\, \{\, \hat{e}^a \,|\; 1 \leqslant a \leqslant \hat{n} \,\}$ \linebreak
of~$\hat{T}^*$, and let $\, \{\, e_i^{} \,|\; 1 \leqslant i \leqslant
N \,\} \,$ be an arbitrary basis of a $k$-isotropic subspace~$E$
complementary to~$L$ in~$V$, with dual basis $\, \{\, e^i \,|\;
1 \leqslant i \leqslant N \,\} \,$ of~$\, L^\bot \cong E^*$. Choosing
an arbitrary subspace $L'$ of~$L$ complementary to~$\, \ker \, \fp \,$
and taking into account the identity~(\ref{eq:POLIF}), we define a basis
$\, \{\, e_a^{i_1 \ldots\, i_k} \,|\; 1 \leqslant a \leqslant \hat{n} \,,\,
1 \leqslant i_1 < \ldots <i_k \leqslant N \,\} \,$ of~$L'$ by
\[
 \fp^\flat (e_a^{i_1 \ldots\, i_k})~
 =~e^{i_1} \,\smwedge \ldots\, \smwedge\, e^{i_k} \,\smotimes\, \hat{e}_a~.
\]
It is easy to see that the union of this basis with that of~$E$
gives a canonical basis of $V$ (or more precisely, of $\, E \oplus
L'$, which is a subspace of~$V$ complementary to~$\, \ker \, \fp$).
Thus we have proved%
\footnote{Clearly, the inductive construction of a $k$-isotropic
subspace $E$ complementary to the poly\-lagrangian subspace $L$, as
presented in the proof of Theorem~\ref{thm:POLIDD}, provides an
explicit iteration procedure for building polylagrangian bases
in a way similar to the well known Gram-Schmidt orthogonalization
process. }
\begin{thm}[Darboux theorem for polylagrangian vector spaces]
\label{thm:DARPA} \mbox{} \\[1mm]
 Every polylagrangian vector space admits a canonical basis.
\end{thm}


\section{Multilagrangian forms on vector spaces} \label{sec:FMLA}

\addtocounter{footnote}{-5}
In this section we deal with ordinary alternating multilinear forms
which are partially horizontal with respect to a given vertical subspace.
To explain what this means, assume that we are given a finite-dimensional
real vector space $W$ together with a fixed subspace~$V$ and a projection
$\pi$ from $W$ to another finite-dimensional real vector space $T$ which
has $V$ as its kernel, so that $\, T \cong W/V$, i.e., we have a short
exact sequence of vector spaces\,\footnotemark\addtocounter{footnote}{4}
\begin{equation} \label{eq:SEQEX1} 
 0~\longrightarrow~V~\longrightarrow~W~\stackrel{\pi}{\longrightarrow}~
                   T~\longrightarrow~0~.
\end{equation}
Motivated by standard jargon of fiber bundle theory, we shall refer to $W$ as
the \textbf{total space}, $V$ as the \textbf{vertical space} and $T$ as the
\textbf{base space}. Then an $r$-form $\, \alpha \,\smin\, \bwedge^r \, W^*$
\linebreak
on~$W$ is said to be \textbf{$(r-s)$-horizontal} (with respect to~$\pi$),
where $\, 0 \leqslant s \leqslant r \,$, if its contraction with more than
$s$ vertical vectors vanishes, i.e., if
\begin{equation} \label{eq:HFORM1} 
 \mathrm{i}_{v_1} \ldots\, \mathrm{i}_{v_{s+1}} \, \alpha~=~0
 \qquad \mbox{for $\, v_1^{},\ldots,v_{s+1}^{} \smin V$}~.
\end{equation}
The vector space of $(r-s)$-horizontal $r$-forms on~$W$ will be denoted
by $\, \bwedge_{\,s}^{\,r} \, W^*$. Note that as $s$ is varied (with $r$
fixed), these spaces form an increasing sequence under inclusion:\,%
\footnote{The first few terms of this sequence may be trivial, since
$\, \bwedge_{\,s}^{\,r} \, W^* = \{0\} \,$ if $\; s < r - \dim T$.}
\begin{equation} \label{eq:INCLU1} 
 \bwedge^r \, T^* \cong \bwedge_{\,0}^{\,r} \, W^* \, \subset \ldots
 \subset \, \bwedge_{\,s}^{\,r} \, W^* \, \subset \ldots
 \subset \, \bwedge_{\,r}^{\,r} \, W^* = \bwedge^r \, W^*~.
\end{equation}
At the two extremes, we have $\, \bwedge_{\,r}^{\,r} \, W^*
= \bwedge^r \, W^* \,$ since the condition of $0$-horizontality
is void, whereas the space $\bwedge_{\,0}^{\,r} \, W^*$ of fully
horizontal $r$-forms on~$W$, which are precisely the horizontal
$r$-forms as defined in \cite[Vol.~2]{GHV} and \cite[Vol.~1]{KN},
is canonically isomorphic to the space $\bwedge^r \, T^*$ of all
$r$-forms on~$T$: $\, \bwedge_{\,0}^{\,r} \, W^* \cong \bwedge^r \, T^*$.
This canonical isomorphism is simply given by pull-back with the
projection $\pi$, i.e., $\alpha_W^{} = \pi^* \alpha_T^{} \,$ or
\begin{equation} \label{eq:CANIS2} 
 \alpha_W^{}(w_1^{},\ldots,w_r^{})~
 =~\alpha_T^{}(\pi(w_1^{}),\ldots,\pi(w_r^{}))
 \qquad \mbox{for $\, w_1^{},\ldots,w_r^{} \smin W$}~.
\end{equation}
Its inverse is given by $\, \alpha_T^{} = \mathfrak{s}^* \alpha_W^{} \,$ or
\begin{equation} \label{eq:CANIS5} 
 \alpha_T^{}(t_1^{},\ldots,t_r^{})~
 =~\alpha_W^{}(\mathfrak{s}(t_1^{}),\ldots,\mathfrak{s}(t_r^{}))
 \qquad \mbox{for $\, t_1^{},\ldots,t_r^{} \smin T$}~,
\end{equation}
where $\mathfrak{s}$ is any splitting of the exact sequence~%
(\ref{eq:SEQEX1}), i.e., any linear mapping from $T$ to~$W$
such that $\, \pi \smcirc \mathfrak{s} = \mathrm{id}_T$.%
\footnote{It is a straightfoward exercise to verify that the expressions
on the rhs of equation~(\ref{eq:CANIS5}) and on the rhs of equation~%
(\ref{eq:SYMBF2}) do not depend on the choice of the splitting
$\mathfrak{s}$ when $\alpha_W^{}$ is horizontal and $\alpha$ is
$(r-s)$-horizontal.}
Extending this construction to partially horizontal forms leads us naturally
to the concept of symbol, which will provide the link between polysymplectic/
polylagrangian and multisymplectic/multilagrangian structures.
\begin{defi}~ \label{def:SYMB}
 Let\/ $W$, $V$ and $T$ be finite-dimensional vector spaces related by
 the short exact sequence~(\ref{eq:SEQEX1}). The \textbf{symbol} of
 an $(r-s)$-horizontal $r$-form $\alpha$ on~$W$, $\alpha \,\smin\,
 \bwedge_{\,s}^{\,r} \, W^*$, is the $\, \bwedge^{r-s} \, T^*$-valued
 $s$-form $\hat{\alpha}$ on~$V$, $\hat{\alpha} \,\smin\, \bwedge^s \, V^*
 \otimes \bwedge^{r-s} \, T^*$, given by
 \begin{equation} \label{eq:SYMBF1} 
  \hat{\alpha}(v_1^{},\ldots,v_s^{})~
  =~\mathrm{i}_{v_1}^{} \ldots\, \mathrm{i}_{v_s}^{} \alpha
  \qquad \mbox{$\mathrm{for}$ $\, v_1^{},\ldots,v_s^{} \smin V$}~.
 \end{equation}
\end{defi}
Explicitly, we may use equation (\ref{eq:CANIS5}) to arrive at the following
formula for $\hat{\alpha}$ in terms of~$\alpha$:
\begin{equation} \label{eq:SYMBF2} 
 \begin{array}{c}
  \hat{\alpha}(v_1^{},\ldots,v_s^{}) \bigl( t_1^{},\ldots,t_{r-s}^{} \bigr)~
  =~\alpha(v_1^{},\ldots,v_s^{},
           \mathfrak{s}(t_1^{}),\ldots,\mathfrak{s}(t_{r-s}^{})) \\[2mm]
  \mbox{for $\, v_1^{},\ldots,v_s^{} \smin V$,
             $t_1^{},\ldots,t_{r-s}^{} \smin T$}~,
 \end{array}
\end{equation}
where, once again, $\mathfrak{s}$ is any splitting of the exact
sequence~(\ref{eq:SEQEX1}).\,\addtocounter{footnote}{-1}\footnotemark\ 
It follows that passage to the symbol can be regarded as a projection, from
the space $\bwedge_{\,s}^{\,r} \, W^*$ of $(r-s)$-horizontal $r$-forms on~$W$
to the space $\, \bwedge^s \, V^* \otimes \bwedge^{r-s} \, T^* \,$ of $s$-forms
on~$V$ with values in the space $\bwedge^{r-s} \, T^*$ of $(r-s)$-forms on~$T$,
whose kernel is the subspace $\bwedge_{\,s-1}^{\,~r} W^*$ of $(r-s+1)$-%
horizontal $r$-forms on~$W$.

If we introduce a basis $\, \{ e_1^V,\ldots,e_m^V,e_1^T,\ldots,e_n^T \} \,$
of~$W$ such that the first $m$ vectors span~$V$ while the last $n$ vectors
span a subspace complementary to~$V$ and hence isomorphic to~$T$, then in
terms of the dual basis $\, \{ e_V^1,\ldots, e_V^m,e_T^1,\ldots,e_T^n \} \,$
of~$W^*$, an arbitrary form $\, \alpha \,\smin\, \bwedge_{\,s}^{\,r} \, W^* \,$
is represented as\,%
\footnote{The expansion in equation~(\ref{eq:EXPAN1}) explains why
forms in $\bwedge_{\,s}^{\,r} \, W^*$ are called $(r-s)$-horizontal:
they are represented as linear combinations of exterior products of
$1$-forms such that, in each term of the sum, at least $r-s$ of them
are horizontal.}
\begin{equation} \label{eq:EXPAN1} 
 \alpha~=~\sum_{t=0}^s \, {\textstyle \frac{1}{t!} \, \frac{1}{(r-t)!}} \;
          \omega_{i_1\ldots\,i_t;\,\mu_1\ldots\,\mu_{r-t}}^{} \,
          e_V^{i_1} \,\smwedge \ldots\, \smwedge\, e_V^{i_t} \;\smwedge\,
          e_T^{\mu_1} \smwedge \ldots\, \smwedge\, e_T^{\mu_{r-t}}~,
\end{equation}
while its symbol $\, \hat{\alpha} \,\smin\, \bwedge^s \, V^* \otimes
\bwedge^{r-s} \, T^* \,$ is represented as
\begin{equation} \label{eq:EXPAN3} 
 \hat{\alpha}~=~{\textstyle \frac{1}{s!} \, \frac{1}{(r-s)!}} \;
          \alpha_{i_1\ldots\,i_s;\,\mu_1\ldots\,\mu_{r-s}}^{} \,
          e_V^{i_1} \,\smwedge \ldots\, \smwedge\, e_V^{i_s} \;\smotimes\,
          e_T^{\mu_1} \smwedge \ldots\, \smwedge\, e_T^{\mu_{r-s}}~,
\end{equation}
which also shows that
\begin{equation} \label{eq:HEPDIM} 
 \dim \, \bwedge_{\,s}^{\,r} \, W^*~
 =~\sum_{t=0}^s {\dim V \choose t} {\dim T \choose r-t}~,
\end{equation}
where it is to be understood that $\; {k \choose l} = 0 \;$ if $\, l > k \,$.

As in the polysymplectic and, more generally, the polylagrangian case, our
definition of a multisymplectic and, more generally, of a multilagrangian
form will be based on the existence of a special type of maximal isotropic
subspace, the only restriction being that we consider only isotropic sub%
spaces of the vertical space~$V$; correspondingly, the concept of maximality
should in this context be understood to mean maximality in the partially
ordered set formed by the isotropic subspaces of~$V$ (not~$W$).
More precisely, suppose we are given a fixed $(k+1-r)$-horizontal
$(k+1)$-form $\,\omega$ on~$W$,
\begin{equation}
 \omega \,\smin\, \bwedge_{\,~r}^{k+1} \, W^*~,
\end{equation}
where $\, 1 \leqslant r \leqslant k+1$,\,%
\footnote{The extreme case of fully horizontal forms ($r=0$) will be excluded
right from the start since it can be reduced to the other extreme case where
the horizontality condition is void ($r=k+1$), substituting the total space
$W$ by the quotient space $T$. Additional restrictions that serve to exclude
other trivial or uninteresting cases will be imposed as we go along.}
and denote by $\,\omega^\flat$ the contraction of $\,\omega$, as defined in
the previous section, which is a linear map $\; \omega^\flat: W \rightarrow
\bwedge^k \, W^*$, as well as its restriction to the vertical subspace $V$,
which is a linear map $\; \omega^\flat: V \rightarrow \bwedge_{\,r-1}^{\,~k}
W^*$. Then defining, for any subspace $L$ of~$V$,
\begin{equation} \label{eq:HORAN1} 
 \bwedge_{\,r-1}^{\,~k} \, L^\bot~
 =~\bwedge^k L^\bot \cap\, \bwedge_{\,r-1}^{\,~k} W^*~,
\end{equation}
we conclude as in the previous section that a subspace $L$ of~$V$ will be
isotropic (with respect to $\,\omega$) if and only if
\begin{equation} \label{eq:ISOSS6} 
 \omega^\flat \bigl( L \bigr)~
 \subset~\omega^\flat \bigl( V \bigr) \,\cap\,
         \bwedge_{\,r-1}^{\,~k}\, L^\bot~,
\end{equation}
and will be maximal isotropic (with respect to $\,\omega$) if and only if
$\; V \cap\, \ker \, \omega \subset L \;$ and
\begin{equation} \label{eq:ISOSS7} 
 \omega^\flat \bigl( L \bigr)~
 =~\omega^\flat \bigl( V \bigr) \,\cap\,
   \bwedge_{\,r-1}^{\,~k}\, L^\bot~.
\end{equation}
As in the previous section, the intersection with the subspace $\omega^\flat
\bigl( V \bigr)$ on the rhs of these relations is superfluous in equation~%
(\ref{eq:ISOSS6}) but not in equation~(\ref{eq:ISOSS7}), and omitting it
here leads to a strengthened form of equation~(\ref{eq:ISOSS7}) which
turns out to provide the key to the theory of multisymplectic and,
more generally, of multilagrangian forms:
\begin{defi}~ \label{def:MULT}
 Let\/ $W$, $V$ and $T$ be finite-dimensional vector spaces related
 by the short exact sequence~(\ref{eq:SEQEX1}), with $\, n = \dim T$,
 and let $\,\omega$ be a non-vanishing $(k+1-r)$-horizontal $(k+1)$-%
 form on\/~$W$, where $\, 1 \leqslant r \leqslant k+1$. We say that
 $\,\omega$ is a \textbf{multilagrangian} form of \textbf{rank}~$N$ and
 \textbf{horizontality degree} $k+1-r$ if\/ $V$ admits a subspace $L$
 of codimension~$N$ which is \textbf{multilagrangian}, i.e., such that
 \begin{equation} \label{eq:MULTF} 
  \omega^\flat \bigl( L \bigr)~=~\bwedge_{\,r-1}^{\,~k}\, L^\bot~.
 \end{equation}
 When $\, k=n$, $r=2 \,$ and $\,\omega$ is non-degenerate, we call
 $\,\omega$ a \textbf{multisymplectic} form. If~the condition of non-%
 degeneracy is dropped, we speak of a \textbf{multipresymplectic} form.
\end{defi}
The first two propositions on multilagrangian forms are entirely analogous
to the corresponding ones for polylagrangian forms. To begin with, we note
that a multilagrangian subspace, when it exists, contains the kernel of
$\,\omega$ (implying, in particular, that $\, \ker \, \omega \subset V$)
and hence really is a special type of maximal isotropic subspace.
\begin{prp}~ \label{prp:MULTLS1} 
 Let\/ $W$, $V$ and $T$ be finite-dimensional vector spaces related by
 the short exact sequence~(\ref{eq:SEQEX1}), with $\, n = \dim T$, and
 let $\,\omega$ be a multilagrangian $(k+1)$-form on~$W$ of rank~$N$ and
 horizontality degree $k+1-r$, where $\, 1 \leqslant r \leqslant k+1$,
 and with multilagrangian subspace $L$. Then $\, N+n \geqslant k \,$
 and $\, k+1-r \leqslant n$, and $L$ contains the kernel of~$\omega$:
 \begin{equation} \label{eq:MULTK1}
  \ker \, \omega~\subset~L~.
 \end{equation}
\end{prp}
\proof~
 First we observe that if $\, N+n<k$, we have $\, \bwedge^k L^\bot = \{0\} \,$
 while if $\, k+1-r>n$, we have $\, \bwedge_{\,r-1}^{\,~k} W^* = \{0\}$, so
 in either case, both sides of the equation~(\ref{eq:MULTF}) vanish, i.e.,
 $L$ is contained in $\, \ker \, \omega \,$ and hence $\, \ker \, \omega \,$
 has codimension $<k$ in~$W$, implying $\, \omega \equiv 0 \,$, since the
 $(k+1)$-form on the quotient space $\, W/\ker \, \omega \,$ induced by
 $\omega$ vanishes identically. Thus supposing that $\, \dim L^\bot
 = N+n \geqslant k \,$ and using that $\, \dim V^\bot = n \,$ and
 $\, V^\bot \smsubset L^\bot$, we conclude that we can, for any
 vector $\, w \,\smin\, W \setminus L$, find a linearly independent
 set of $1$-forms $\, w_1^*,\ldots,w_k^* \smin L^\bot \,$ such that
 $\, w_r^*,\ldots,w_k^* \smin V^\bot$, $\langle w_1^*,w \rangle = 1 \,$
 and $\, \langle w_i^*,w \rangle = 0$ for $\, i>1$. \linebreak According
 to the definition of a multilagrangian subspace, there is a vector
 $\, u \smin L \,$ such that
 \[
  \mathrm{i}_u^{} \omega~
  =~w_1^* \,\smwedge \ldots \smwedge\, w_k^*
  \quad \Rightarrow \quad
  \mathrm{i}_w^{} \mathrm{i}_u^{} \omega~
  =~w_2^* \,\smwedge \ldots \smwedge\, w_k^*~
  \ne~0
 \]
 and so $\, w \nsmin \ker \, \omega$. Hence it follows that
 $\, \ker \, \omega \smsubset L$.
\qed

\noindent
The second gives a simple dimension criterion that allows to decide whether
a given isotropic subspace of~$V$ is multilagrangian:
\begin{prp}~ \label{prp:MULTLS2} 
 Let\/ $W$, $V$ and $T$ be finite-dimensional vector spaces related by
 the short exact sequence~(\ref{eq:SEQEX1}), with $\, n = \dim T$, and
 let $\,\omega$ be a non-vanishing $(k+1-r)$-horizontal $(k+1)$-form on~$W$,
 where $1 \leqslant r \leqslant k+1 \,$ and $\, k+1-r \leqslant n$. Given
 any subspace $L$ of\/~$V$, with $\, N = \dim (V/L)$, such that $\, N+n
 \geqslant k$, the following statements are equivalent:
 \begin{itemize}
  \item $L$ is a multilagrangian subspace and $\,\omega$ is a multilagrangian
        form of rank $N$.
  \item $L$ contains $\, \ker \, \omega$, is isotropic and has dimension
        \begin{equation} \label{eq:MULTD1}
         \dim \, L~=~\dim \, \ker \, \omega \, + \,
                     \sum_{s=0}^{r-1} {N \choose s} {n \choose k-s}~,
        \end{equation}
        where it is to be understood that $\, {N \choose s} = 0 \;$ if
        $\, s > N \,$.
 \end{itemize}
\end{prp}
\proof~
 Taking into account that, for any isotropic subspace $L$ of~$V$
 containing the kernel of $\omega$, the contraction map $\,\omega^\flat$
 induces an injective linear map of $\, L/\ker \, \omega \,$ into
 $\, \bwedge_{\,r-1}^{\,~k}\, L^\bot$, we obtain
 \[
  \dim L \, - \, \dim \, \ker \, \omega~
  =~\dim \; \omega^\flat \bigl( L \bigr)~
  =~\dim \, \bwedge_{\,r-1}^{\,~k}\, L^\bot~.
 \]
 To calculate this dimension, we introduce a basis $\, \{ e_1^L,\ldots,e_l^L,
 e_1^{L'},\ldots,e_N^{L'},e_1^T,\ldots,e_n^T \}$ \linebreak of~$W$ such that
 the first $l$ vectors form a basis of~$L$, the following $N$ vectors form
 a basis of a subspace $L'$ complementary to~$L$ in~$V$ and the last $n$
 vectors form a basis of a subspace $H$ complementary to~$V$ in~$W$,
 which is isomorphic to~$T$. \linebreak Then in terms of the dual basis
 $\, \{ e_L^1,\ldots,e_L^l,e_{L'}^1,\ldots,e_{L'}^N,e_T^1,\ldots,e_T^n \} \,$
 of~$W^*$, we conclude that $\; \{ \, e_{L'}^{i_1}\,\smwedge \ldots
 \smwedge\, e_{L'}^{i_s} \,\smwedge\, e_T^{\mu_1} \smwedge \ldots
 \smwedge\, e_T^{\mu_{k-s}}~|~0 \leqslant s \leqslant r-1 \,,\,
 1 \leqslant i_1 < \ldots < i_s \leqslant N \,,\linebreak
 1 \leqslant \mu_1 < \ldots <  \mu_{k-s} \leqslant n \, \} \;$
 is a basis of $\, \bwedge_{\,r-1}^{~k} \, L^\bot$.
\qed

\noindent
The relation between multilagrangian and polylagrangian forms is established
through the symbol:
\begin{thm}~ \label{thm:SIMBA} 
 Let\/ $W$, $V$ and $T$ be finite-dimensional vector spaces related by
 the short exact sequence~(\ref{eq:SEQEX1}), with $\, n = \dim T$, let
 $\,\omega$ be a non-vanishing $(k+1-r)$-horizontal $(k+1)$-form on~$W$,
 where $1 \leqslant r \leqslant k+1 \,$ and $\, k+1-r \leqslant n$, and
 let $\,\fp$ be its symbol, which is a $\bwedge^{k+1-r} \, T^*$-valued
 $r$-form on~$V$. Suppose that $\,\omega$ is  multilagrangian, with
 multi\-lagrangian subspace $L$. Then $\,\fp$ is polylagrangian, with
 polylagrangian subspace $L$, and
 \begin{equation} \label{eq:SIMBA1}
  \ker \, \omega~\subset~\ker \, \fp~.
 \end{equation}
 If $\,\omega$ is multipresymplectic, then $\,\fp$ is polypresymplectic, with
 \begin{equation} \label{eq:SIMBA2}
  \dim \ker \, \fp \, - \, \dim \ker \, \omega~\leqslant~1~.
 \end{equation}
\end{thm}
\proof~
 Fixing an arbitrary horizontal subspace $H$ of $W$ and using the direct
 decompositions $\, W = V \oplus H \,$ and $\, W^* = H^\bot \oplus V^\bot$,
 with $\, H^\bot \cong V^* \,$ and $\, V^\bot \cong H^*$, we note that in
 order to show that $L$ is polylagrangian with respect to $\,\fp$, we must
 establish the equality
 \[
  \fp^\flat \bigl( L \bigr)~
  \cong~\bwedge^{r-1} (L^\bot \cap H^\bot) \,\otimes\, \bwedge^{k+1-r} \, T^*~.
 \]
 To do so, we use the canonical isomorphism $\; \bwedge^{k+1-r} \, T^*
 \cong \bwedge_{~~~0}^{k+1-r} \, W^* = \bwedge^{k+1-r} \, V^\bot \;$ and
 the~inclusion $\, V^\bot \smsubset L^\bot$, together with the fact that
 the space $\, \bwedge^{r-1} (L^\bot \cap H^\bot) \,\otimes\, \bwedge^{k+1-r} \,
 V^\bot$ is generated by elements which can be written in the form
 \[
  \hat{\alpha}~
  =~\bigl( w_1^* \,\smwedge \ldots\, \smwedge\, w_{r-1}^* \bigr) \otimes
    \bigl( w_r^* \,\smwedge \ldots\, \smwedge\, w_k^* \bigr)
 \qquad \mbox{where} \qquad
  \alpha~=~w_1^* \,\smwedge \ldots\, \smwedge\, w_k^*
 \]
 with $\, w_1^*,\ldots,w_{r-1}^* \,\smin\, L^\bot \cap H^\bot \,$ and
 $\, w_r^*,\ldots,w_k^* \,\smin\, V^\bot \smsubset L^\bot$. Since $L$
 is multilagrangian with respect to~$\omega$, there is a vector
 $\, u \smin L \,$ such that $\, \alpha = \mathrm{i}_u^{} \omega \,$
 and so $\, \hat{\alpha} = \mathrm{i}_u^{} \fp$, showing that $L$ is
 polylagrangian with respect to~$\fp$. The inclusion~(\ref{eq:SIMBA1})
 is obvious from the definition~(\ref{eq:SYMBF1}) or~(\ref{eq:SYMBF2}),
 and the estimate~(\ref{eq:SIMBA2}) follows by observing that the linear
 map $\; \omega^\flat: V \rightarrow \bwedge_{\,r-1}^{\,~k} W^*$ \linebreak
 induces an injective linear map from $\, V / \ker \, \omega \,$ to
 $\bwedge_{\,r-1}^{\,~k} W^*$ that takes $\, \ker \, \fp / \ker \, \omega \,$
 into $\bwedge_{\,r-2}^{\,~k} W^*$, and this space is one-dimensional when
 $\,\omega$ is multipresymplectic.
\qed

\noindent
In particular, we note that when $\,\omega$ is multisymplectic, there are two
possibilities: either $\,\fp$ is already non-degenerate or else $\,\fp$ has a
one-dimensional kernel. Both cases do arise in practice, so all that can be
said in general about the relation between multisymplectic and polysymplectic
forms (in addition to the last statement of Theorem~\ref{thm:SIMBA}) is that
the symbol $\,\fp$ of a multisymplectic form $\,\omega$ on~$W$ induces a
polysymplectic form on the quotient space $\, \hat{W} = W/\ker \, \fp$,
which by abuse of notation will again be denoted by~$\,\fp$.

Before going on to explore general consequences of the relation between
multi\-lagrangian and polylagrangian forms that we have just established,
let us pause for a moment to comment on the situations encountered when
the parameter $r$ is assumed to take on one of its extreme values.
On the one hand, when $\, r = k+1$, the horizontality condition
becomes void, the choice of $V$ and $T$ becomes irrelevant, the
space $\bwedge^{k+1-r} \, T^*$ is one-dimensional and the form
$\,\omega$ coincides with its symbol $\,\fp$. This is, once again,
the particular case that has been studied before in the literature~%
\cite{Mar,CIL,LDS} under the label ``multi\-symplectic forms'' (a
terminology that we propose to abandon since we use this term in a
different sense; see Definition~\ref{def:MULT} above). Thus it becomes
clear that these ``multisymplectic \linebreak forms'' appear as the
intersection between polylagrangian and multilagrangian forms: they
are polylagrangian without being vector-valued as well as multi%
lagrangian without being horizontal. On the other hand, note that
we have already excluded the case $r\!=\!0$, where the $(k+1)$-form
$\,\omega$ is fully horizontal (i.e., $(k+1)$-horizontal), since this
situation can be reduced to the previous one if we substitute $W$
by~$T$, but even the next case $r\!=\!1$, where the $(k+1)$-form
$\,\omega$ is $k$-horizontal, is essentially trivial, since every
$k$-horizontal multilagrangian $(k+1)$-form has rank $0$ and multi%
lagrangian subspace $V$. (In fact, the condition that $\,\omega$ should
be $k$-horizontal is equivalent to the condition that $V$ should be
isotropic, and in this case, $V$ does satisfy the remaining criteria
of Proposition~\ref{prp:MULTLS2}.)

Another fundamental property of multilagrangian forms that can be derived
from the corresponding property for polylagrangian forms is that the multi%
lagrangian subspace has a particular type of direct complement.
\begin{thm}~ \label{thm:MULTDD} 
 Let\/ $W$, $V$ and $T$ be finite-dimensional vector spaces related by
 the short exact sequence~(\ref{eq:SEQEX1}), with $\, n = \dim T$, and
 let $\,\omega$ be a multilagrangian $(k+1)$-form on~$W$ of rank~$N$ and
 horizontality degree $k+1-r$, where $\, 1 \leqslant r \leqslant k+1 \,$
 and $\, k+1-r \leqslant n$, and with multilagrangian subspace $L$.
 Then there exists a $k$-isotropic subspace $F$ of\/~$W$ such that
 the intersection $\, E = V \smcap F \,$ is an $(r-1)$-isotropic
 subspace of\/~$V$ and
 \begin{equation} \label{eq:MULTDD}
   W~=~F \oplus L~~~,~~~V~=~E \oplus L~.
 \end{equation}
\end{thm}
\proof~
 First we construct an $(r-1)$-isotropic subspace $E$ of~$V$ of dimension~$N$
 which is complementary to~$L$ in~$V$. If $\, r=1 \,$ there is nothing to prove
 since in this case the vertical subspace $V$ is isotropic and so we have
 $\, L = V$, $E = \{0\}$ and $N=0 \,$. If $\, r>1 \,$ we apply Theorem~%
\ref{thm:POLIDD} to the symbol $\fp$ of~$\omega$ to conclude that there
 is a subspace $E$ of~$V$ of dimension~$N$ which is complementary to~$L$
 in~$V$ and is $(r-1)$-isotropic with respect to~$\fp$. Now taking into
 account that the whole vertical subspace $V$ is $r$-isotropic with respect
 to~$\omega$, it follows that $E$ is $(r-1)$-isotropic with respect to~%
 $\omega$ as well.\ --\, Now let $F_0$ be a subspace of~$W$ of dimension
 $N+n'$ which is $k$-isotropic with respect to~$\omega$ and such that
 $\, F_0 \cap V = E$. (For instance, if $\, n'=0$, $F_0 = E$).
 If $\, n'=n$, we are done. Otherwise, choose a basis $\, \{ e_1^E,
 \ldots,e_N^E,e_1^{},\ldots,e_n^{} \} \,$ of a subspace of~$W$
 complementary to~$L$  such that the first $N$ vectors constitute a
 basis of~$E$ and the first $N+n'$ vectors constitute a basis of~$F_0$,
 and denote the corresponding dual basis of~$L^\bot$ by $\, \{ e_E^1,
 \ldots,e_E^N,e^1,\ldots,e^n \}$. We shall prove that there exists
 a vector $\, u \smin W \setminus (F_0 \oplus L) \,$ such that the
 subspace $F_1$ spanned by $u$ and $F_0$ is $k$-isotropic and
 satisfies $\, F_1 \cap V = E$; then since $\, \dim F_1 = N+n'+1$,
 the statement of the theorem follows by induction. To this end,
 choose any subspace $L'$ of~$L$ complementary to $\, \ker \, \omega$
 and use the fact that $L$ is multilagrangian to conclude that there
 exists a unique basis
 \[
  \{ \, e^{i_1 \ldots\, i_s;\,\mu_1 \ldots\, \mu_{k-s}} \, |~
     0 \leqslant s \leqslant r-1 \,,\,
     \parbox{44mm}{$~~1 \leqslant i_1^{} < \ldots < i_s^{} \leqslant N~~ \\
                      1 \leqslant \mu_1 < \ldots < \mu_{k-s}^{}
                      \leqslant n$}
     \, \}
 \]
 of $L'$ such that
 \[
  \omega_V^\flat (e^{i_1 \ldots\, i_s;\,\mu_1 \ldots\, \mu_{k-s}})~
  =~e_E^{i_1} \>\smwedge \ldots\, \smwedge\, e_E^{i_s} \>\smwedge\,
    e^{\mu_1} \,\smwedge \ldots\, \smwedge\, e^{\mu_{k+1-s}}~.
 \]
 Thus, for $\, 0 \leqslant s,t \leqslant r-1$, $1 \leqslant i_1 < \ldots
 < i_s \leqslant N$, $1 \leqslant j_1 < \ldots <j_t \leqslant N$, \linebreak
 $1 \leqslant \mu_1 < \ldots < \mu_{k-s} \leqslant n$, $1 \leqslant \nu_1
 < \ldots <\nu_{k-t} \leqslant n$, we have
 \[
   \omega(e^{i_1 \ldots\, i_s;\,\mu_1 \ldots\, \mu_{k-s}},
          e_{j_1}^E,\ldots,e_{j_t}^E,e_{\nu_1}^{},\ldots,e_{\nu_{k-t}}^{})~
   =~\left\{
      \begin{array}{ccc}
       0 & \mbox{if} & s \neq t \\[2mm]
       \delta_{j_1}^{i_1} \ldots\, \delta_{j_s}^{i_s} \,
       \delta_{\nu_1^{}}^{\mu_1^{}} \ldots\, \delta_{\nu_{k-s}}^{\mu_{k-s}} &
           \mbox{if} & s=t
      \end{array}
     \right\}~.
 \]
 Therefore, the vector
 \[
  u~=~e_{n'+1}~-~\sum_{s=0}^{r-1} \,
      {\textstyle \frac{1}{s!} \, \frac{1}{(k-s)!}} \;
      \omega(e_{n'+1}^{},e_{i_1}^E,\ldots,e_{i_s}^E,
                         e_{\mu_1}^{},\ldots,e_{\mu_{k-s}}^{}) \;
      e^{i_1 \ldots\, i_s;\,\mu_1 \ldots\, \mu_{k-s}}
 \]
 does not belong to the subspace $\, F_0 \oplus L \,$ and, for
 $\, 0 \leqslant t \leqslant r-1$, $1 \leqslant j_1 < \ldots <j_t
 \leqslant N \,$ and $\, 1 \leqslant \nu_1 < \ldots < \nu_{k-t}
 \leqslant n$, satisfies
 \begin{eqnarray*}
  \lefteqn{\omega(u,e_{j_1}^E,\ldots,e_{j_t}^E,
                    e_{\nu_1}^{},\ldots,e_{\nu_{k-t}}^{})} \\[2mm]
  &=&\!\! \omega(e_{n'+1}^{},e_{j_1}^E,\ldots,e_{j_t}^E,
                             e_{\nu_1}^{},\ldots,e_{\nu_{k-t}}^{}) \\
  & &     \mbox{} - \, \sum_{s=0}^{r-1} \,
          {\textstyle \frac{1}{s!} \, \frac{1}{(k-s)!}} \;
          \omega(e_{n'+1}^{},e_{i_1}^E,\ldots,e_{i_s}^E,
                             e_{\mu_1}^{},\ldots,e_{\mu_{k-s}}^{}) \\[-2mm]
  & & \hspace{3cm} \times \;
          \omega(e^{i_1 \ldots\, i_s;\,\mu_1 \ldots\, \mu_{k-s}},
                 e_{j_1}^E,\ldots,e_{j_t}^E,
                 e_{\nu_1}^{},\ldots,e_{\nu_{k-t}}^{})
  \\[2mm]
  &=&\!\! 0~.
 \end{eqnarray*}
 which implies that since the subspace $F_0$ spanned by $\, e_1^E,\ldots,e_N^E,
 e_1^{},\ldots,e_{n'}^{} \,$ is $k$-isotropic, the subspace $F_1$ spanned by
 $\, e_1^E, \ldots,e_N^E,e_1^{},\ldots,e_{n'}^{} \,$ and $u$ is so as well.
\qed

\begin{exem}(\textbf{The canonical form})~~
 Let $F$ be a vector space of dimension $N+n$ and $E$ be a fixed $N$-%
 dimensional subspace of~$F$. Denoting the $n$-dimensional quotient
 space $F/E$ by~$T$ and the canonical projection of~$F$ onto~$T$ by
 $\rho$, we obtain the following exact sequence of vector spaces:
 \begin{equation} \label{eq:SEQEX2} 
  0~\longrightarrow~E~\longrightarrow~F~
    \stackrel{\rho}{\longrightarrow}~T~\longrightarrow~0~.
 \end{equation}
 Set
 \begin{equation} \label{eq:FMLC1} 
  W_0^{}~=~F \,\oplus\, \bwedge_{\,r-1}^{\;~k} F^*~~~,~~~
  V_0^{}~=~E \,\oplus\, \bwedge_{\,r-1}^{\;~k} F^*~~~,~~~
  \pi_0^{}~=~\rho \smcirc \mathrm{pr}_1~,
 \end{equation}
 where $\, \mathrm{pr}_1:  W_0^{} \rightarrow F \,$ is the canonical
 projection, which leads us to the following exact sequence of vector
 spaces:
\begin{equation} \label{eq:SEQEX3} 
 0~\longrightarrow~V_0^{}~\longrightarrow~W_0^{}~
   \stackrel{\pi_0^{}}{\longrightarrow}~T~\longrightarrow~0~.
\end{equation}
 The \textbf{canonical multilagrangian form} of rank $N$ and horizontality
 degree $k+1-r$ is the $(k+1-r)$-horizontal $(k+1)$-form $\,\omega_0^{}$
 on~$W_0^{}$ defined by
 \begin{equation} \label{eq:FMLC2} 
   \omega_0^{} \bigl( (u_0^{},\omega_0^{}), \ldots ,
                      (u_k^{},\omega_k^{}) \bigr)~
   =~\sum_{i=0}^k \, (-1)^i \,
     \omega_i^{}(u_0^{},\ldots,\widehat{u_i^{}},\ldots,u_k^{})~.
 \end{equation}
 If $\, k=n \,$ and $\, r=2 \,$ we shall call $\,\omega_0^{}$ the
 \textbf{canonical multisymplectic form}.
\end{exem}
It is a straightforward exercise to show that $\,\omega_0^{}$ is $(k+1-r)$-%
horizontal and is degenerate when $\, r=1$, with
\begin{equation} \label{eq:FMLCK} 
 \ker \, \omega_0^{}~=~E \qquad \mbox{if $\, r=1$}~,
\end{equation}
but is non-degenerate when $\, r>1$. In what follows, we shall assume that
$\, N>0$, since when $\, E = \{0\}$, we are back to the polylagrangian case
on~$\; T \oplus \bwedge^k T^*$, with $\, \hat{T} = \mathbb{R}$, which has
already been studied in Refs~\cite{Mar,CIL,LDS}. For the same reason,
we shall also assume that $\, r>1$, since for $\, r=1 \,$ we have
$\, \bwedge_{\,0}^{\,k} F^* \cong \bwedge^k T^*$, so that after passing
to the quotient by the kernel of~$\omega_0^{}$, we are once again back
to the polylagrangian case on~$\; T \oplus \bwedge^k T^*$, with
$\, \hat{T} = \mathbb{R}$. Then considering $F$ and
\begin{equation} \label{eq:FMLC3} 
 L~=~\bwedge_{\,r-1}^{\;~k} F^*
\end{equation}
as subspaces of~$W_0^{}$, we have the direct decompositions
$\, W_0^{} = F \oplus L \,$ and $\, V_0^{} = E \oplus L \,$ where
\begin{equation} \label{eq:FMLC4} 
 L~\text{is multilagrangian}
 \quad \text{,} \quad
 F~\text{is $k$-isotropic}
 \quad \text{,} \quad
 E~\text{is $(r-1)$-isotropic}~.
\end{equation}
In terms of bases, let $\, \{\, e_i^{} , e_\mu^{} \,|\; 1 \leqslant i
\leqslant N \,,\, 1 \leqslant \mu \leqslant n \,\} \,$ be a basis of~$F$
with dual basis $\, \{\, e^i , e^\mu \,|\; 1 \leqslant i \leqslant N \,,\,
1 \leqslant \mu \leqslant n \,\} \,$ of~$F^*$ such that $\, \{\, e_i^{} \,|\;
1 \leqslant i \leqslant N \,\} \,$ is a basis of~$E$ and $\, \{\, e_\mu^{}
\,|\; 1 \leqslant \mu \leqslant n \,\} \,$ is a basis of a subspace $H$
of~$F$ complementary to~$E$, isomorphic to~$T$. For $\, 0 \leqslant s
\leqslant r$, $1 \leqslant i_1^{} < \ldots < i_s^{} \leqslant N \,$ and
$\, 1 \leqslant \mu_1^{} < \ldots < \mu_{k-s}^{} \leqslant N$, define
\[
 \begin{array}{c}
 e_{i_1 \ldots\, i_s;\,\mu_1 \ldots\, \mu_{k-s}}^{}~
 =~e_{i_1}^{} \,\smwedge \ldots\, \smwedge\, e_{i_s}^{} \,\smwedge\,
   e_{\mu_1}^{} \,\smwedge \ldots\, \smwedge\, e_{\mu_{k-s}}^{} \\[2mm]
 e^{i_1 \ldots\, i_s;\,\mu_1 \ldots\, \mu_{k-s}}~
 =~e^{i_1} \,\smwedge \ldots\, \smwedge\, e^{i_s} \,\smwedge\,
   e^{\mu_1} \,\smwedge \ldots\, \smwedge\, e^{\mu_{k-s}}
\end{array}~.
\]
This provides a basis
\[
 \{ \, e_i^{} , e_\mu^{} , e^{i_1 \ldots\, i_s;\,\mu_1 \ldots\, \mu_{k-s}} \,
    |~0 \leqslant s \leqslant r-1 \,,\,
    \parbox{67mm}{$~~1 \leqslant i \leqslant N \,,\,
                     1 \leqslant i_1^{} < \ldots < i_s^{} \leqslant N~~ \\
                     1 \leqslant \mu \leqslant n \,,\,
                     1 \leqslant \mu_1 < \ldots < \mu_{k-s}^{}
                     \leqslant n$}
    \, \}
\]
of~$W_0^{}$ with dual basis
\[
 \{ \, e^i , e^\mu , e_{i_1 \ldots\, i_s;\,\mu_1 \ldots\, \mu_{k-s}}^{} \,
    |~0 \leqslant s \leqslant r-1 \,,\,
    \parbox{67mm}{$~~1 \leqslant i \leqslant N \,,\,
                     1 \leqslant i_1^{} < \ldots < i_s^{} \leqslant N~~ \\
                     1 \leqslant \mu \leqslant n \,,\,
                     1 \leqslant \mu_1 < \ldots < \mu_{k-s}^{}
                     \leqslant n$}
    \, \}
\]
of~$W_0^*$, both of which we shall refer to as a \textbf{canonical basis}
or \textbf{Darboux basis}, such that
\begin{equation} \label{eq:FMLCB1} 
 \omega_0^{}~
 =~\sum_{s=0}^{r-1} \, {\textstyle \frac{1}{s!} \frac{1}{(k-s)!}} \,
   e_{i_1 \ldots\, i_s;\,\mu_1 \ldots\, \mu_{k-s}}^{} \,\smwedge\,
   e^{i_1} \,\smwedge \ldots\, \smwedge\, e^{i_s} \,\smwedge\,
   e^{\mu_1} \,\smwedge \ldots\, \smwedge\, e^{\mu_{k-s}}~.
\end{equation}
and for the symbol
\begin{eqnarray} \label{eq:SMLCB}
 \fp_0^{} \;
 =~{\textstyle \frac{1}{(r-1)!} \, \frac{1}{(k+1-r)!}} \;
   \bigl( e_{i_1 \ldots\, i_{r-1};\,\mu_1 \ldots\, \mu_{k+1-r}}^{} \,\smwedge\,
          e^{i_1} \,\smwedge \ldots\, \smwedge\, e^{i_{r-1}} \bigr) \smotimes
   \bigl( e^{\mu_1} \,\smwedge \ldots\, \smwedge\, e^{\mu_{k+1-r}} \bigr) \, .~
\end{eqnarray}

\pagebreak

Now it is easy to derive the algebraic Darboux theorem for general
multilagrangian forms: let $\, \{\, e_i^{} , e_\mu^{} \,|\; 1 \leqslant
i \leqslant N \,,\, 1 \leqslant \mu \leqslant n \,\} \,$ be a basis of
a $k$-isotropic subspace~$F$ complementary to~$L$ in~$W$, with dual
basis $\, \{\, e^i , e^\mu \,|\; 1 \leqslant i \leqslant N \,,\,
1 \leqslant \mu \leqslant n \,\} \,$ of~$\, L^\bot \cong F^*$, such that
$\, \{\, e_i^{} \,|\; 1 \leqslant i \leqslant N \,\} \,$ is a basis of the
$(r-1)$-isotropic subspace $\, E = V \smcap F \,$ which is complementary
to~$L$ in~$V$. Choosing any subspace $L'$ of~$L$ complementary to~%
$\, \ker \, \omega$ and taking into account the identity~%
(\ref{eq:MULTF}), we define a basis
\[
 \{ \, e^{i_1 \ldots\, i_s;\,\mu_1 \ldots\, \mu_{k-s}} \,
    |~0 \leqslant s \leqslant r-1 \,,\,
    \parbox{44mm}{$~~1 \leqslant i_1^{} < \ldots < i_s^{} \leqslant N~~ \\
                     1 \leqslant \mu_1 < \ldots < \mu_{k-s}^{}
                     \leqslant n$}
 \, \}
\]
of~$L'$ by
\[
 \omega_V^\flat (e^{i_1 \ldots\, i_s;\,\mu_1 \ldots\, \mu_{k-s}})~
 =~e^{i_1} \,\smwedge \ldots \smwedge\, e^{i_s} \,\smwedge\,
   e^{\mu_1} \,\smwedge \ldots \smwedge\, e^{\mu_{k+1-s}}~.
\]
It is easy to see that the union of this basis with that of~$F$ gives a
canonical basis of $W$ (or more precisely, of $\, F \oplus L'$, which is a
subspace of~$W$ complementary to~$\, \ker \, \omega$). Thus we have proved%
\footnote{Once again, the inductive construction of a $k$-isotropic
subspace $F$ complementary to the multi\-lagrangian subspace $L$,
as presented in the proof of Theorem~\ref{thm:MULTDD}, provides an
explicit iteration procedure for building multilagrangian bases in
a way similar to the well known Gram-Schmidt orthogonalization
process.}
\begin{thm}[Darboux theorem for multilagrangian vector spaces]
\label{thm:DARMA} \mbox{} \\[1mm]
 Every multilagrangian vector space admits a canonical basis.
\end{thm}


\section{Cartan calculus and the symbol}
\label{sec:CCFV}

In order to extend the structures studied in the previous two sections and,
in particular, the concept of symbol that interrelates them, from a purely
algebraic setting to the realm of differential geometry, we shall need a 
variant of Cartan's calculus, which in its standard formulation deals with
differential forms on manifolds, to handle vertical differential forms on
total spaces of fiber bundles.

Let $P$ be a fiber bundle over a base manifold~$M$, with projection
$\, \pi : P \rightarrow M$. Then the vector bundle $\bwedge_s^r \, T^* P$
over~$P$ whose fiber at any point $p$ in $P$ is the space $\bwedge_s^r \,
T_p^* P$ of $(r-s)$-horizontal $r$-forms on the tangent space $T_p P$ to~$P$
at~$p$ is called the \textbf{bundle of $(r-s)$-horizontal $r$-forms} on~$P$,
and its sections are called \textbf{$(r-s)$-horizontal differential $r$-forms}
or simply \textbf{$(r-s)$-horizontal $r$-forms} on~$P$; the space of such forms
will be denoted by $\Omega_{\,s}^{\,r}(P)$. Similarly, assuming in addition
that $\hat{T}$ is a vector bundle over the same base manifold~$M$, with
projection $\, \hat{\tau} : \hat{T} \rightarrow M$, and denoting the
pull-back of $\hat{T}$ to $P$ by $\pi^* \hat{T}$ and the vertical bundle
of~$P$ by $VP$ (both are vector bundles over~$P$), the vector bundle
$\, \bwedge^r V^* P \otimes \pi^* \hat{T} \,$ over~$P$ is called the
\textbf{bundle of vertical $r$-forms} on~$P$, and its sections are called
\textbf{vertical differential $r$-forms} or simply \textbf{vertical $r$-forms}
on~$P$, with \textbf{values} or \textbf{coefficients} in $\pi^* \hat{T}$ or,
by abuse of language, in $\hat{T}$: the space of such forms will be denoted
by $\Omega_V^{\,r}(P;\pi^* \hat{T})$. Finally, the sections of the vertical
bundle $VP$ itself are called \textbf{vertical vector fields} or simply
\textbf{vertical fields} on~$P$: the space of such fields will be denoted
by $\mathfrak{X}_V^{}(P)$. Obviously, $\Omega_{\,s}^{\,r}(P)$, $\Omega_V^{\,r}%
(P;\pi^* \hat{T})$ and $\mathfrak{X}_V^{}(P)$ are (locally finite) modules
over the algebra $\mathfrak{F}(P)$ of functions on~$P$.

It should be noted that speaking of vertical forms constitutes another
abuse of language because these ``forms'' are really equivalence classes
of differential $r$-forms on~$P$: $\Omega_V^{\,r}(P;\pi^* \hat{T})$ is not
a subspace of the space $\Omega^r(P;\pi^* \hat{T})$ of all differential
$r$-forms on~$P$ but rather its quotient space
\[
 \Omega_V^{\,r}(P;\pi^* \hat{T})~
 =~\Omega^r(P;\pi^* \hat{T}) / \Omega_{r-1}^{\,~r}(P;\pi^* \hat{T})
\]
by the subspace $\Omega_{r-1}^{\,~r}(P;\pi^* \hat{T})$ of all
$1$-horizontal differential $r$-forms on~$P$.

An interesting aspect of this construction is that it is possible to
develop a variant of the usual Cartan calculus for differential forms
on the manifold $P$ in which vector fields on~$P$ are replaced by
vertical fields $X$ on~$P$ and differential forms on~$P$ (taking values
in some fixed vector space) are replaced by vertical differential forms
$\alpha$ on~$P$ (taking values in some fixed vector bundle over the base
manifold), in such a way that all operations of this calculus such as
exterior multiplication, contraction, Lie derivative and exterior
derivative continue to be well defined and to satisfy the standard
rules. (See~\cite[Vol.~1, Problem~8, p.~313]{GHV} for the special
case where the vector bundle of coefficients is the trivial line
bundle $\, M \times \mathbb{R}$.) Here, we shall only need the
\textbf{vertical exterior derivative}
\begin{equation} \label{eq:VEXTD1} 
 \begin{array}{cccc}
  d_V^{} : & \Omega_V^{\,r}(P;\pi^* \hat{T})
           & \longrightarrow & \Omega_{~V}^{\,r+1}(P;\pi^* \hat{T}) \\[1mm]
           & \alpha & \longmapsto & d_V^{} \alpha
 \end{array}~,
\end{equation}
which is defined by exactly the same formula as in the standard case, namely
\begin{equation} \label{eq:VEXTD2} 
 \begin{array}{rcl}
 d_V^{} \alpha \, (X_0^{},\ldots,X_r^{}) \!\! 
 &=&\!\! {\displaystyle
          \sum_{i=0}^r \, (-1)^i \; X_i^{} \cdot
          \left( \alpha(X_0^{},\ldots,\hat{X}_i^{},\ldots,X_r^{}) \right)}
 \\[4mm]
 & & +\, {\displaystyle
          \sum_{0 \leqslant i < j \leqslant r} (-1)^{i+j} \,
          \alpha([X_i^{},X_j^{}],X_0^{},\ldots,
                 \hat{X}_i^{},\ldots,\hat{X}_j^{},\ldots,X_r^{})}
 \end{array}~,
\end{equation}
where $\; X_0^{},X_1^{},\ldots,X_r^{} \smin\, \mathfrak{X}_V^{}(P) \,$:
this makes sense since $V\!P$ is an involutive distribution on~$P$,
provided we correctly define the \textbf{vertical directional
derivative}
\begin{equation} \label{eq:VDIRD1} 
 \begin{array}{ccc}
  \mathfrak{X}_V^{}(P) \times \Gamma(\pi^* \hat{T})
  & \longrightarrow & \Gamma(\pi^* \hat{T}) \\[1mm]
  (X,\varphi) & \longmapsto & X \!\cdot \varphi
 \end{array}~.
\end{equation}
as an $\mathbb{R}$-bilinear operator which is $\mathfrak{F}(P)$-linear
in the first entry and satisfies a Leibniz rule in the second entry,
\begin{equation} \label{eq:VDIRD2} 
 X \!\cdot (f \varphi)~
 =~(X \cdot f) \, \varphi + f (X \!\cdot \varphi)~.
\end{equation}
Explicitly, for $\, X \smin\, \mathfrak{X}_V^{}(P) \,$ and $\, \varphi
\,\smin\, \Gamma(\pi^* \hat{T})$, $X \!\cdot \varphi \,\smin\,
\Gamma(\pi^* \hat{T}) \,$ is defined as the standard directional
derivative of vector valued functions along the fibers, that is,
for any point $m$ in~$M$, $\, (X \!\cdot \varphi) \big|_{P_m} \smin\,
C^\infty(P_m,\hat{T}_m) \,$ is given in terms of $\, X \big|_{P_m}
\smin\, \mathfrak{X}(P_m)$ and $\, \varphi \big|_{P_m} \smin\,
C^\infty(P_m,\hat{T}_m) \,$ by
\begin{equation} \label{eq:VDIRD3} 
 (X \!\cdot \varphi) \big|_{P_m}~=~X \big|_{P_m} \!\cdot \varphi \big|_{P_m}~.
\end{equation}
Since the Lie bracket is natural under restriction to submanifolds,
we have
\[
 X \!\cdot \bigl( Y \!\cdot \varphi \bigr) \, - \,
 Y \!\cdot \bigl( X \!\cdot \varphi \bigr)~=~[X,Y] \!\cdot \varphi
 \qquad \mbox{for $\, X,Y \smin\, \mathfrak{X}_V^{}(P)$,
                   $\varphi \smin\, \Gamma(\pi^* \hat{T})$}~,
\]
which implies that $\, d_V^2 = 0$. On the other hand, sections
$\varphi$ of $\pi^* \hat{T}$ obtained from sections $\hat{t}$
of~$\hat{T}$ by composing with $\pi$ are constant along the fibers
and hence their vertical directional derivative vanishes:
\begin{equation} \label{eq:VDIRD4} 
 X \cdot (\hat{t} \smcirc \pi)~=~0
 \qquad \mbox{for $\, X \smin\, \mathfrak{X}_V^{}(P)$,
                   $\hat{t} \smin\, \Gamma(\hat{T})$}~.
\end{equation}
In the same way, substituting $\hat{T}$ by $\hat{T}^*$, we get
\begin{equation} \label{eq:VDIRD5} 
 X \cdot (\hat{t}^* \smcirc \pi)~=~0
 \qquad \mbox{for $\, X \smin\, \mathfrak{X}_V^{}(P)$,
                   $\hat{t}^* \smin\, \Gamma(\hat{T}^*)$}~,
\end{equation}
which means that
\begin{equation} \label{eq:VDIRD6} 
 X \cdot \langle \hat{t}^* \smcirc \pi \,,\, \varphi \rangle~
 =~\langle \hat{t}^* \smcirc \pi \,,\, (X \cdot \varphi) \rangle
 \qquad \mbox{for $\, X \smin\, \mathfrak{X}_V^{}(P)$,
                   $\varphi \smin\, \Gamma(\pi^* \hat{T})$}~.
\end{equation}
More generally, given a $\pi^* \hat{T}$-valued vertical $r$-form~%
$\hat{\alpha}$ on~$P$ and a section $\hat{t}^*$ of the dual vector
bundle $\hat{T}^*$ of~$\hat{T}$, we define the \textbf{projection}
of~$\hat{\alpha}$ along~$\hat{t}^*$ to be the ordinary vertical
$r$-form~$\hat{\alpha}_{\hat{t}^*}^{}$ on~$P$ given by\,%
\footnote{Here, it is important that $\hat{t}^*$ be a section of
$\hat{T}^*$ and not of $\pi^* \hat{T}^*$.}
\begin{equation}
 \hat{\alpha}_{\hat{t}^*}^{}(p)~
 =~\langle \hat{t}^*(\pi(p)) , \hat{\alpha}(p) \rangle
 \qquad \mbox{for $\, p \,\smin P$}~,
\end{equation}
and obtain
\begin{equation} \label{eq:VEXTD3} 
 d_V^{} \hat{\alpha}_{\hat{t}^*}^{}~
 =~\bigl( d_V^{} \hat{\alpha} \bigr)_{\hat{t}^*}~.
\end{equation}
Hence, if $\hat{\alpha}$ is closed, $\hat{\alpha}_{\hat{t}^*}^{}$ will
be closed as well, and a standard argument shows that, in this case,
the \textbf{kernel} of~$\hat{\alpha}$ and of~$\hat{\alpha}_{\hat{t}^*}^{}$,
given by
\begin{equation}
 \ker_p^{} \hat{\alpha}~=~\ker \, \hat{\alpha}_p^{}~
 =~\{ \, u_p^{} \smin V_p^{} P~|~
      \mathrm{i}_{u_p}^{} \hat{\alpha}_p^{} = 0 \, \}
 \qquad \mbox{for $\, p \,\smin P$}~,
\end{equation}
and
\begin{equation}
 \ker_p^{} \hat{\alpha}_{\hat{t}^*}^{}~
 =~\ker \, (\hat{\alpha}_{\hat{t}^*}^{})_p^{}~
 =~\{ \, u_p^{} \smin V_p^{} P~|~
      \mathrm{i}_{u_p}^{} (\hat{\alpha}_{\hat{t}^*}^{})_p^{} = 0 \, \}
 \qquad \mbox{for $\, p \,\smin P$}~,
\end{equation}
respectively, define involutive distributions on~$P$, provided they have
constant dimension.

With these tools at our disposal, we can now formulate the construction
of the symbol in the differential geometric setting. To deal directly with
the case of interest, suppose that $\,\omega$ is a $(k+1-r)$-horizontal
$(k+1)$-form on~$P$, where $\, 1 \leqslant r \leqslant k+1 \,$ and
$\, k+1-r \leqslant n$, that is,
\begin{equation}
 \omega~\smin~\Omega_{\,\>~r}^{\,k+1} \bigl( P \bigr)~
 =~\Gamma \bigl( \bwedge_{\>~r}^{k+1} \, T^* P \bigr)~.
\end{equation}
Explicitly, the horizontality condition means that contraction of $\,\omega$
with more than $r$ vertical fields on~$P$ gives zero. Then the \textbf{symbol}
$\,\fp$ of $\omega$, whose value at every point $p$ of~$P$ is defined to be
the symbol $\,\fp_p$ of $\,\omega_p$, is a vertical $r$-form on~$P$ taking
values in the vector bundle $\, \pi^* \bigl( \bwedge^{k+1-r} \, T^* M \bigr)$,
that is,
\begin{equation}
 \fp~\smin~\Omega_{\,V}^{\,r}
           \bigl( P , \pi^* \bigl( \bwedge^{k+1-r} \, T^* M \bigr) \bigr)~
 =~\Gamma \bigl( \bwedge^r \, V^* P \otimes
                 \pi^* \bigl( \bwedge^{k+1-r} \, T^* M \bigr) \bigr)~.
\end{equation}
Using the canonical isomorphism $\, \pi^* \bigl( \bwedge^{k+1-r} \, T^* M \bigr)
\cong \bwedge_{~~~0}^{k+1-r} \, T^* P \,$ of vector bundles over~$P$ as an
identification, we have
\begin{equation}
 \fp(X_1,\ldots,X_r)~=~\mathrm{i}_{X_1}^{} \ldots\, \mathrm{i}_{X_r}^{} \omega
 \qquad \mbox{for $\, X_1^{},\ldots,X_r^{} \smin\, \mathfrak{X}_V^{}(P)$}~.
\end{equation}
More explicitly, using the horizontal lift
\[
 \begin{array}{ccc}
  \mathfrak{X}(M) & \longrightarrow & \mathfrak{X}(P) \\[1mm]
        \xi       &   \longmapsto   &      \xi^H
 \end{array}
\]
of vector fields induced by some fixed connection on $P$, we have
\begin{equation}
 \begin{array}{c}
  \fp(X_1^{},\ldots,X_r^{}) \cdot
  \bigl( \xi_1^{} \smcirc \pi ,\ldots, \xi_{k+1-r}^{} \smcirc \pi \bigr)~
  =~\omega(X_1^{},\ldots,X_r^{},\xi_1^H,\ldots,\xi_{k+1-r}^H) \\[2mm]
  \mbox{for $\, X_1^{},\ldots,X_r^{} \smin\, \mathfrak{X}_V^{}(P)$,
            $\xi_1^{},\ldots,\xi_{k+1-r}^{} \smin\, \mathfrak{X}(M)$}~,
 \end{array}
\end{equation}
where it should be noted that although composition of vector fields
on~$M$ with the projection $\pi$ provides only a subspace of the
vector space of all sections of the pull-back $\pi^* TM$ of $TM$
by~$\pi$, this formula is sufficient to fix the value of $\,\fp$ at
each point $p$ of~$P$.
\begin{thm}~ \label{thm:SIMBG1} 
 With the same notations as above, suppose that the form $\,\omega$ satisfies
 \[
  d\omega~\smin~\Omega_{\,\>~r}^{\,k+2} \bigl( P \bigr)~
  =~\Gamma \bigl( \bwedge_{\>~r}^{k+2} \, T^* P \bigr)~,
 \]
 i.e., $d\omega$ is $(k+2-r)$-horizontal. Then the form $\,\fp$ is vertically
 closed:
 \[
  d_V^{} \fp~=~0~.
 \]
\end{thm}
In particular, for $\,\fp$ to be vertically closed, it is sufficient (but not
necessary) that $\,\omega$ be closed.

\noindent
\proof~
 Let $\, X_0^{},\ldots,X_r \,$ be vertical fields on~$P$ and $\, \xi_1^{},
 \ldots,\xi_{k+1-r}^{} \,$ be vector fields on~$M$, and denote the horizontal
 lifts of the latter with respect to some fixed connection on~$P$ by
 $\, \xi_1^H,\ldots,\xi_{k+1-r}^H$, respectively. Then we have the
 relation
\begin{equation*}
 d\omega \bigl( X_0^{},\ldots,X_r^{},\xi_1^H,\ldots,\xi_{k+1-r}^H\bigr)~
 =~\Bigl( d_V^{} \fp \bigl( X_0^{},\ldots,X_r^{} \bigr) \Bigr)
                     \bigl( \xi_1^{} \smcirc \pi ,\ldots,
                            \xi_{k+1-r}^{} \smcirc \pi \bigr)~,
\end{equation*}
 which can be derived from Cartan's formula for $d\omega$ and its analogue
 for $d_V^{} \fp$, equation~(\ref{eq:VEXTD2}), using that $\,\omega$ is a
 $(k+1-r)$-horizontal $(k+1)$-form and the Lie brackets $[X_i^{},\xi_j^H]$
 are vertical fields, by applying equation~(\ref{eq:VDIRD6}) to expressions
 of the form $\; \varphi = \fp(X_0^{},\ldots,\hat{X}_i^{},\ldots,X_r^{})$
 \linebreak
 and $\; \hat{t}^* = \xi_1^{} \,\smwedge \ldots \smwedge\, \xi_{k+1-r}^{}$.
\qed


\section{Polylagrangian and multilagrangian fiber bundles}
\label{sec:FPML}

Now we are ready to transfer the poly- and multilagrangian structures
introduced in the first two sections from the algebraic to the differential
geometric context. All we need to do is add the appropriate integrability
condition, which is the expected one: the differential forms in question
should be closed. For the sake of brevity, we begin directly with the
notion of a polylagrangian fiber bundle, which formalizes the idea of
a ``family of polylagrangian manifolds smoothly parametrized by the
points of a base manifold~$M$'' and includes that of a polylagrangian
manifold as a special case.
\begin{defi}~ \label{def:FibPoliLag}
 A \textbf{polylagrangian fiber bundle} is a fiber bundle $P$ over
 an $n$-dimensional manifold~$M$ equipped with a vertical $(k+1)$-form
 $\,\fp$ of constant rank on the total space~$P$ taking values in a fixed
 $\hat{n}$-dimensional vector bundle $\hat{T}$ over the same manifold~$M$,
 called the \textbf{polylagrangian form along the fibers} of~$P$, or simply
 the \textbf{polylagrangian form}, and said to be of \textbf{rank}~$N$,
 such that $\,\fp$ is vertically closed,
 \begin{equation}
  d_V^{} \fp~=~0~,
 \end{equation}
 and such that at every point $p$ of~$P$, $\fp_p$ is a polylagrangian form
 of rank $N$ on the vertical space~$V_p P$. If the polylagrangian subspaces
 at the different points of~$P$ fit together into a distribution $L$ on~$P$
 (which is contained in the vertical bundle $VP$ of~$P$), we call it the
 \textbf{polylagrangian distribution} of $\,\fp$. \\
 When $\, k=1 \,$ and $\,\fp$ is non-degenerate, we say that $P$ is a
 \textbf{polysymplectic fiber bundle} and $\,\fp$ is a \textbf{poly%
 symplectic form along the fibers} of~$P$, or simply a \textbf{poly%
 symplectic form}. If the condition of non-degeneracy is dropped, we
 call $P$ a  \textbf{polypresymplectic fiber bundle} and $\,\fp$ a
 \textbf{polypresymplectic form along the fibers} of~$P$, or simply
 a \textbf{polypresymplectic form}. \\
 If $M$ reduces to a point, we speak of a \textbf{polylagrangian manifold}
 or \textbf{poly\-(pre)\-symplectic manifold}, respectively.
\end{defi}
The notion of a multilagrangian fiber bundle is defined similarly.
\begin{defi}~ \label{def:FibMultLag}
 A \textbf{multilagrangian fiber bundle} is a fiber bundle $P$ over
 an $n$-dimensional manifold~$M$ equipped with a $(k+1-r)$-horizontal
 $(k+1)$-form $\,\omega$ of constant rank on the total space~$P$, where
 $\, 1 \leqslant r \leqslant k+1 \,$ and $\, k+1-r \leqslant n$, called
 the \textbf{multi\-lagrangian form} and said to be of \textbf{rank} $N$
 and \textbf{horizontality degree} $k+1-r$, such that $\,\omega$ is closed,
 \begin{equation}
  d \omega~=~0~,
 \end{equation}
 and such that at every point $p$ of~$P$, $\omega_p$ is a multilagrangian form
 of rank $N$ on the tangent space~$T_p P$. If the multilagrangian subspaces
 at the different points of~$P$ fit together into a distribution $L$ on~$P$
 (which is contained in the vertical bundle $VP$ of~$P$), we call it the
 \textbf{multilagrangian distribution} of~$\,\omega$.

 \pagebreak

 \noindent
 When $\, k=n$, $r=2 \,$ and $\,\omega$ is non-degenerate, we say
 that $P$ is a \textbf{multisymplectic fiber bundle} and $\,\omega$
 is a \textbf{multisymplectic form}. If the condition of non-degeneracy
 is dropped, we call $P$ a  \textbf{multipresymplectic fiber bundle} and
 $\,\omega$ a \textbf{multipresymplectic form}. \\
 If $M$ reduces to a point, we speak of a \textbf{multilagrangian manifold}.
\end{defi}
Note that within the logic underlying this terminology, the concept of a
multisymplectic manifold by itself is meaningless, since multisymplectic
structures ``live'' on total spaces of fiber bundles over some non-trivial
base space (representing space-time). The ``multi\-symplectic manifolds'' of
Refs~\cite{Mar,CIL,LDS} correspond to what we call multilagrangian manifolds.

\noindent
Combining Theorem~\ref{thm:SIMBA} and Theorem~\ref{thm:SIMBG1}, we obtain
\begin{thm}~ \label{thm:SIMBG2}
 Let $P$ be a fiber bundle over an $n$-dimensional manifold~$M$, with
 projection $\, \pi : P \rightarrow M$, let $\,\omega$ be a $(k+1-r)$-%
 horizontal $(k+1)$-form of constant rank on~$P$, where $\, 1 \leqslant
 r \leqslant k+1 \,$ and $\, k+1-r \leqslant n$, and let $\,\fp$ be its
 symbol, which is a vertical $r$-form on~$P$ taking values in the bundle
 of $(k+1-r)$-forms on~$M$. Suppose that $\,\omega$ is multilagrangian,
 with multilagrangian distribution~$L$. Then $\,\fp$ is polylagrangian,
 with polylagrangian distribution~$L$, and
 \begin{equation} \label{eq:SIMBG1}
  \ker \, \omega~\subset~\ker \, \fp~.
 \end{equation}
 If $\,\omega$ is multipresymplectic, then $\,\fp$ is polypresymplectic, with
 \begin{equation} \label{eq:SIMBG2}
  \dim \ker \, \fp \, - \, \dim \ker \, \omega~\leqslant~1~.
 \end{equation}
\end{thm}
In particular, we note again that when $\,\omega$ is multisymplectic, there are
two possibilities: either $\,\fp$ is already non-degenerate or else $\,\fp$ has
a one-dimensional kernel. Both cases \linebreak do arise in practice, so all
that can be said in general about the relation between multi\-symplectic and
polysymplectic structures on fiber bundles (in addition to the last statement
of Theorem~\ref{thm:SIMBG2}) is that the symbol $\,\fp$ of a multisymplectic
form $\,\omega$ on~$P$ induces a polysymplectic form on the quotient bundle
$\, \hat{P} = P/\ker \, \fp \,$ (provided that this quotient is really a
fiber bundle over~$M$), which by abuse of notation will again be denoted
by~$\,\fp$.

\begin{rmk}~
 Note that when $\, \hat{n} \geqslant 2$, in the case of polylagrangian
 structures, or when $\, 1 \leqslant r \leqslant k$, in the case of multi%
 lagrangian structures, the explicit construction of~$L$ in terms of the
 kernels of the projections of $\,\fp$ (see Theorem~\ref{thm:POLILS1})
 implies that the poly\-lagrangian or multilagrangian subspaces at the
 different points of~$P$ do fit together to form a uniquely determined,
 smooth distribution $L$ on~$P$. (The proof uses the fact that sums and
 intersections of smooth vector subbundles of a vector bundle are again
 smooth vector subbundles if they have constant rank and, more generally,
 that kernels and images of smooth vector bundle homomorphisms of constant
 rank are smooth vector subbundles; \linebreak

 \pagebreak

 \noindent
 see, e.g., \cite[Exercise 1.6F(c), p.~51]{AM}.) But in the case of scalar
 polylagrangian structures ($\hat{n}=1$) and of multilagrangian structures
 for which the horizontality condition is void ($r=k+1$), such a distribution
 may fail to exist, even if we assume the base manifold to be trivial
 (i.e., $M$ reduces to a point, $n=0$) and the form $\,\fp$ or $\,\omega$
 to be non-degenerate. Indeed, for symplectic forms ($k=1$) or volume forms
 ($N=k$), it is easy to find examples of manifolds whose tangent bundle
 does not admit any (smooth) lagrangian subbundle: \linebreak the simplest
 of them all is just the $2$-sphere $S^2$, according to the ``no hair
 theorem''. \linebreak And even in the remaining cases ($N > k > 1$),
 where according to Proposition~\ref{prp:POLILS3}, the polylagrangian or
 multilagrangian subspaces at the different points of~$P$ are unique, an
 explicit construction for them does not seem to be available, so we do
 not know whether it is true that smoothness of $\,\fp$ or $\,\omega$ by
 itself implies smoothness of~$L$. Thus when we \linebreak simply refer
 to ``the polylagrangian distribution~$L$'' or ``the multilagrangian
 distribution~$L$'' \linebreak without further specification, as will
 often be done in what follows, existence and smoothness of $L$ is in
 these cases tacitly assumed.
\end{rmk}


\section{Integrability}
\label{sec:VPL}

A further remarkable property of polylagrangian/multilagrangian structures
is that, generically, the polylagrangian/multilagrangian distribution is
automatically integrable and thus gives rise to a \textbf{polylagrangian}/%
\textbf{multilagrangian} \textbf{foliation}.
\begin{thm}[Integrability Theorem]~ \label{thm:PMLInt}
 Let $P$ be a polylagrangian fiber bundle over an $n$-dimensional manifold~$M$
 with polylagrangian $(k+1)$-form $\,\fp$ of rank~$N$ taking values in a fixed
 $\hat{n}$-dimensional vector bundle $\hat{T}$ over the same manifold~$M$.
 Then if
 \begin{equation}
  \hat{n} \geqslant 3~,
 \end{equation}
 the polylagrangian distribution $L$ is integrable. Similarly, let $P$ be
 a multilagrangian fiber bundle over an $n$-dimensional manifold~$M$ with
 multilagrangian $(k+1)$-form $\,\omega$ of rank~$N$ and horizontality degree
 $k+1-r$, where $\, 1 \leqslant r \leqslant k+1 \,$ and $\, k+1-r \leqslant n$.
 Then if
 \begin{equation}
  {n \choose k+1-r} \geqslant 3~,
 \end{equation}
 the multilagrangian distribution $L$ is integrable.
\end{thm}
\proof~
 Using Theorem~\ref{thm:SIMBG2}, the second statement is easily reduced
 to the first. To prove this, suppose that $X$ and $Y$ are vector fields
 on~$P$ which are sections of $\, L \,\smsubset\, VP$. Using the decomposition
 \[
  L~=~K_0^{} \oplus K_1^{} \oplus \ldots \oplus K_{\hat{n}}^{}~,
 \]
 with $\, K_0^{} = \ker \, \fp$, as in equation~(\ref{eq:POLIK3}), we can
 decompose $X$ and $Y$ according to
 \[
  X~=~\sum_{a=0}^{\hat{n}} X_a^{}~~~,~~~Y~=~\sum_{b=0}^{\hat{n}} Y_b^{}~,
 \]
 where $X_a^{}$ and $Y_b^{}$ are sections of $K_a^{}$ and $K_b^{}$,
 respectively. Using that $\, \hat{n} \geqslant 3$, we can for each value
 of $a$ and $b$ find a value $c \neq 0$ such that $\, c \neq a \,$ and
 $\, c \neq b \,$; then $\, K_a^{} \,\smsubset \ker \, \omega^c \,$ and
 $\, K_b^{} \,\smsubset \ker \, \omega^c$. Since $\omega^c$ is vertically
 closed and has constant rank, $\ker \, \omega^c \,\smsubset\, VP \,$ is
 involutive. Therefore the vector field $[X_a^{},Y_b^{}]$ is a section
 of  $\, \ker \, \omega^c \,\smsubset\, L$.
\qed

\noindent
It must be emphasized that when the above inequalities are not satisfied,
the statement of Theorem~\ref{thm:PMLInt} is false, i.e., $L$ may fail to
be involutive. For the multilagrangian case, it is instructive to spell out
explicitly under what circumstances this may happen:
\begin{enumerate}
 \item $r=k+1$, $n = 0,1,2,3,\ldots$ arbitrary: this is the extreme case
       discussed before in which the horizontality condition is void.
 \item $r=k+1-n$, $n = 1,2,3,\ldots$ arbitrary: this includes the
       symplectic case, obtained by choosing $\, k=n+1 \,$ and
       $\, r=2$ (i.e., $\omega \,\smin\, \Omega_{\;~2}^{\,n+2}
       \bigl( P \bigr)$) and supposing in addition that
       $\, \ker \, \fp = \{0\} \,$ and that $M$ is orientable,
       so that $\; \fp \,\smin\, \Omega_{\,V}^{\,2} \bigl( P , \pi^*
       \bigl( \bwedge^n T^*\!M \bigr) \bigr)$ \linebreak represents
       a ``family of symplectic forms smoothly parametrized by the
       points of a base manifold~$M$ of dimension~$n$''. Here, it is
       not difficult to construct examples of lagrangian distributions
       which are not involutive.
 \item $r=k$, $n=2$ (the only possibility to have $\, {n \choose k+1-r}
       \leqslant 2 \,$ with $\, 1 \leqslant r \leqslant k \,$ and
       $\, k+1-r < n$): this includes the multisymplectic case over
       a two-dimensional base manifold~$M$, obtained by choosing
       $\, r=k=n=2$. An explicit example of this situation is given
       in Example~\ref{exe:MSB-LNINV} below.
\end{enumerate}
A simple example of a polysymplectic manifold with a two-dimensional
non-involutive polylagrangian distribution of rank $1$ is the following.
\begin{exem}~ \label{exe:PSM-LNINV}
 Let $P$ be $SU(2)$ and let $\, \alpha \smin \Omega^1(P,\mathfrak{su}(2)) \,$
 be the left invariant Maurer-Cartan form on~$P$, where $\mathfrak{su}(2)$ is
 the Lie algebra of~$SU(2)$. Consider the left invariant $\mathfrak{su}(2)^*$-%
 valued $2$-form $\beta$ on~$P$ obtained by taking the exterior product of
 component forms whose values are multiplied using the commutator in~%
 $\mathfrak{su}(2)$ and finally passing to the dual using the invariant
 scalar product $(.\,,.)$ on~$\mathfrak{su}(2)$, i.e.,
 \[
  \langle \beta \,, X \rangle~
  =~{\textstyle \frac{1}{2}} \, ( [ \, \alpha \kwedge \alpha \, ] \,, X )
  \qquad \mbox{for $\, X \smin \mathfrak{su}(2)$}~.
 \]
 Using the isomorphism $\, \mathfrak{su}(2) \cong \mathbb{R}^3 \,$ induced
 by employing the standard orthonormal basis $\, \{ \, \sigma_a^{}/2i \; | \;
 a \smin \{1,2,3\} \} \,$ of~$\mathfrak{su}(2)$,%
 \footnote{Note that under the isomorphism $\, \mathfrak{su}(2) \cong
 \mathbb{R}^3 \,$ given by mapping the basis $\, \{ \, \sigma_a^{}/2i \; | \;
 a \smin \{1,2,3\} \}$ \linebreak of $\mathfrak{su}(2)$ to the standard basis
 of~$\mathbb{R}^3$, where
 \[
  \sigma_1^{}~=~\left( \begin{array}{cc}
                        0 & 1 \\
                        1 & 0
                       \end{array} \right)~,~~
  \sigma_2^{}~=~\left( \begin{array}{cc}
                        0 & -i \\
                        i & 0
                       \end{array} \right)~,~~
  \sigma_3^{}~=~\left( \begin{array}{cc}
                        1 & 0 \\
                        0 & -1
                       \end{array} \right)
 \]
 are the usual Pauli matrices, the commutator in~$\mathfrak{su}(2)$
 corresponds to the vector product in~$\mathbb{R}^3$ whereas the
 invariant scalar product in~$\mathfrak{su}(2)$ given by $\, (X,Y)
 = - 2 \, \mathrm{tr} (XY) \,$ corresponds to the standard scalar
 product in $\mathbb{R}^3$. With respect to this orthonormal basis,
 the (totally covariant) structure constants for~$\mathfrak{su}(2)$
 are given by the components of the $\epsilon$-tensor, which are
 those of the standard volume form on~$\mathbb{R}^3$.}
 and working in components, we see that $\alpha$ is represented by a
 triplet of left invariant $1$-forms $\,\alpha^1$, $\alpha^2$, $\alpha^3$
 and $\beta$ by a triplet of left invariant $2$-forms $\,\beta_1^{}$,
 $\beta_2^{}$, $\beta_3^{}$ on~$P$ such that
 \[
  \beta_a^{}~=~{\textstyle \frac{1}{2}} \, \epsilon_{abc}^{} \,
               \alpha^b \,\smwedge\, \alpha^c~.
 \]
 Moreover, taking $\, \hat{T} = (\mathbb{R}^2)^*$, let $e$ be a given
 $2$-frame in $\mathbb{R}^3$, i.e., a given injective linear map $e$
 from~$\mathbb{R}^2$ to~$\mathbb{R}^3$. Its transpose will be a sur%
 jective linear map (projection) $\mathrm{pr}$ from~$(\mathbb{R}^3)^*$
 to~$(\mathbb{R}^2)^*$,%
 \footnote{It may seem overly pedantic not to identify the spaces
 $\mathbb{R}^2$ and~$\mathbb{R}^3$ with their respective duals
 $(\mathbb{R}^2)^*$ and~$(\mathbb{R}^3)^*$, but we have refrained
 from doing so right from the start since this turns out to facilitate
 the understanding of the generalizations to be discussed in the next
 two examples. Moreover, maintaining this distinction is not completely
 irrelevant since $e$ is \emph{not} assumed to be isometric: the space
 $\mathbb{R}^2$ is not even supposed to carry a scalar product.}
 and we define
 \begin{equation} \label{eq:PSFEX1}
  \fp~=~\mathrm{pr} \>\smcirc\, \beta~.
 \end{equation}
 In components, $e$ is represented by a triplet of vectors $e^1$, $e^2$, $e^3$
 in~$(\mathbb{R}^2)^*$ (not a basis, of course), and we have
 \begin{equation} \label{eq:PSFEX2}
  \fp~=~{\textstyle \frac{1}{2}} \, \epsilon_{abc}^{} \;
        \alpha^a \,\smwedge\, \alpha^b \,\smotimes\, e^c~.
 \end{equation}
 It follows immediately from the Maurer-Cartan structure equations that each
 of the $2$-forms $\,\beta_1^{}$, $\beta_2^{}$, $\beta_3^{}$ is closed, and
 hence so is $\,\fp$. Moreover, it is easily verified that the image of~%
 $\mathbb{R}^2$ under $e$, regarded as a subspace of $\, \mathfrak{su}(2)
 \cong \mathbb{R}^3$, generates a two-dimensional left invariant
 distribution $L$ on~$P$ which is isotropic under $\,\fp$ (since
 the expression $\, \epsilon_{abc}^{} \, e^a(u) \, e^b(v) \, e^c(w) \,$
 vanishes identically, for any three vectors $u,v,w$ in $\mathbb{R}^2$),
 so $\,\fp$ is polysymplectic and~$L$ is the corresponding polylagrangian
 distribution. Of course, $L$ is not involutive.
\end{exem}
Similar but somewhat more complicated constructions lead to examples of
poly\-symplectic and multisymplectic fiber bundles over an arbitrary two-%
dimensional base manifold~$M$ whose polylagrangian or multilagrangian
distribution, again of rank $1$, is not involutive. \linebreak
All of them are principal bundles over~$M$, with structure group
$SU(2)$ in the poly\-symplectic case and $U(2)$ in the multisymplectic
case, carrying an invariant polysymplectic or multisymplectic form built
from data that can be interpreted in terms of concepts from Yang-Mills-%
Higgs field theory, namely an $SU(2)$-connection form $A$ and a Higgs field
$\phi$ in the adjoint representation of~$SU(2)$ in the polysymplectic case
and a $U(1)$-connection form $A$ and a Higgs field $\phi$ in the truncated
adjoint representation of~$U(2)$ in the multisymplectic case.%
\footnote{By the truncated adjoint representation of a Lie group~$G$,
we mean the restriction of the adjoint representation of~$G$ on its Lie
algebra~$\mathfrak{g}$ to its derived algebra $[\mathfrak{g},\mathfrak{g}]$,
which is an $\mathrm{Ad}(G)$-invariant subspace (in particular, an ideal)
of~$\mathfrak{g}$. Note that the derived algebra of $\mathfrak{u}(n)$ is
$\mathfrak{su}(n)$.}
\addtocounter{footnote}{-3}
The details follow.
\begin{exem}~ \label{exe:PSB-LNINV}
 Let $P$ be the total space of a principal bundle over a two-dimensional
 manifold~$M$ with structure group~$SU(2)$ and bundle projection $\; \pi:
 P \rightarrow M \,$ and let $\, A \smin \Omega^1(P,\mathfrak{su}(2)) \,$
 be a given connection form on~$P$, where $\mathfrak{su}(2)$ is the Lie
 algebra of~$SU(2)$. Consider the equivariant $\mathfrak{su}(2)^*$-valued
 $2$-form $B$ on~$P$ obtained by taking the exterior product of component
 forms whose values are multiplied using the commutator in~$\mathfrak{su}(2)$
 and finally passing to the dual using the invariant scalar product $(.\,,.)$
 on~$\mathfrak{su}(2)$, i.e.,
 \[
  \langle B \,, X \rangle~
  =~{\textstyle \frac{1}{2}} \, ( [ \, A \kwedge A \, ] \,, X )
  \qquad \mbox{for $\, X \smin \mathfrak{su}(2)$}~.
 \]
 Using the isomorphism $\, \mathfrak{su}(2) \cong \mathbb{R}^3 \,$ as before,%
 \footnotemark \addtocounter{footnote}{2}
 and working in components, we see that $A$ is represented by a triplet of
 $1$-forms $A^1$, $A^2$, $A^3$ and $B$ by a triplet of $2$-forms $B_1^{}$,
 $B_2^{}$, $B_3^{}$ \linebreak on $P$ such that
 \[
  B_a^{}~=~{\textstyle \frac{1}{2}} \, \epsilon_{abc}^{} \,
           A^b \,\smwedge\, A^c~.
 \]
 Moreover, taking $\, \hat{T} = T^* M$, let $e$ be a given $1$-form on~$M$
 taking values in the adjoint bundle $\, P \times_{SU(2)} \mathbb{R}^3 \,$
 which is an ``immersion'' in the sense that, when interpreted as a vector
 bundle homomorphism from~$TM$ to~$\, P \times_{SU(2)} \mathbb{R}^3$, it
 is fiberwise injective. When pulled~back to~$P$, $e$ corresponds to an
 equivariant horizontal $\mathbb{R}^3$-valued $1$-form $\pi^* e$ on~$P$
 which, once again, is an ``immersion'' in the sense that, when inter%
 preted as an equivariant vector bundle homomorphism from~$\pi^*(TM)$
 to~$\, P \times \mathbb{R}^3$, it is fiberwise injective. Its transpose
 will be a fiberwise surjective equivariant vector bundle homomorphism
 (equivariant projection) $\mathrm{pr}$ from~$\, P \times (\mathbb{R}^3)^*$
 to~$\pi^*(T^* M)$, and we define
 \begin{equation} \label{eq:PSFEX3}
  \fp~=~\mathrm{pr} \>\smcirc\, B \big|_{V\!P}~,
 \end{equation}
 where $.\,\big|_{V\!P}$ denotes restriction to the vertical bundle $V\!P$
 of~$P$. In components, $\pi^* e$ \linebreak is represented by a triplet
 of sections $\pi^* e^1$, $\pi^* e^2$, $\pi^* e^3$ of~$\pi^*(T^* M)$,
 and we have
 \begin{equation} \label{eq:PSFEX4}
  \fp~=~{\textstyle \frac{1}{2}} \, \epsilon_{abc}^{} \;
        A^a \big|_{V\!P} \>\smwedge\, A^b \big|_{V\!P} \,\smotimes\,
        \pi^* e^c~.
 \end{equation}
 It follows immediately from the structure equation
 \[
  d_V^{} A^a \big|_{V\!P}~=~(dA^a) \big|_{V\!P}~
  = \; - \, {\textstyle \frac{1}{2}} \, \epsilon_{\phantom{\,a}bc}^{\,a} \;
    A^b \big|_{V\!P} \>\smwedge\, A^c \big|_{V\!P}~,
 \]
 together with the formula $\; \epsilon_{abc}^{} \,
 \epsilon_{\phantom{\,a}de}^{\,a} = \delta_{bd}^{} \, \delta_{ce}^{} \, - \,
 \delta_{be}^{} \, \delta_{cd}^{} \,$, that $\,\fp$ is vertically closed:
 \begin{eqnarray*}
  d_V^{} \fp \!\!
  &=&\!\! {\textstyle \frac{1}{2}} \, \epsilon_{abc}^{} \,\,
          d_V^{} A^a \big|_{V\!P} \>\smwedge\,
          A^b \big|_{V\!P} \,\smotimes\, \pi^* e^c \, - \,
          {\textstyle \frac{1}{2}} \, \epsilon_{abc}^{} \;
          A^a \big|_{V\!P} \>\smwedge\,
          d_V^{} A^b \big|_{V\!P} \,\smotimes\, \pi^* e^c \\[2mm]
  &=&\!\! \mbox{} - \, {\textstyle \frac{1}{2}} \,
          \epsilon_{abc}^{} \, \epsilon_{\phantom{\,a}de}^{\,a} \;
          A^d  \big|_{V\!P} \>\smwedge\, A^e \big|_{V\!P} \>\smwedge\,
          A^b \big|_{V\!P} \,\smotimes\, \pi^* e^c \\
  & &\!\! \mbox{} + \, {\textstyle \frac{1}{2}} \,
          \epsilon_{abc}^{} \, \epsilon_{\phantom{\,b}de}^{\,b} \;
          A^a  \big|_{V\!P} \>\smwedge\, A^d \big|_{V\!P} \>\smwedge\,
          A^e \big|_{V\!P} \,\smotimes\, \pi^* e^c \\[2mm]
  &=&\!\! 0~.
 \end{eqnarray*}
 Moreover, it is easily verified that the image of~$\pi^*(TM)$ under
 $\pi^* e$, when transferred from $P \times \mathbb{R}^3$ to $V\!P$
 by means of the canonical isomorphism that exists between the two
 since $P$ is a principal bundle, constitutes a two-dimensional
 invariant distribution~$L$ on~$P$ which is isotropic under~$\,\fp$,
 so $\,\fp$ is polysymplectic and~$L$ is the corresponding polylagrangian
 distribution. Of course, $L$ is not involutive.
\end{exem}
\begin{exem}~ \label{exe:MSB-LNINV}
 Let $P$ be the total space of a principal bundle over a two-dimensional
 manifold~$M$ with structure group~$U(2)$ and bundle projection $\; \pi:
 P \rightarrow M \,$ and let $\, A \smin \Omega^1(P,\mathfrak{u}(2)) \,$
 be a given connection form on~$P$, where $\mathfrak{u}(2)$ is the Lie
 algebra of~$U(2)$. Consider the equivariant $\mathfrak{su}(2)^*$-valued
 $2$-form $B$ on~$P$ obtained by taking the exterior product of component
 forms whose values are multiplied using the commutator in~$\mathfrak{u}(2)$
 (which maps into the derived algebra $\mathfrak{su}(2)$) and finally
 passing to the dual using the invariant scalar product $(.\,,.)$ on~%
 $\mathfrak{su}(2)$, i.e.,
 \[
  \langle B \,, X \rangle~
  =~{\textstyle \frac{1}{2}} \, ( [ \, A \kwedge A \, ] \,, X )
  \qquad \mbox{for $\, X \smin \mathfrak{su}(2)$}~.
 \]
 Using the isomorphism $\, \mathfrak{u}(2) \cong \mathbb{R}^4 \,$ induced
 by employing the standard orthonormal basis $\, \{ \, \sigma_a^{}/2i \; | \;
 a \smin \{0,1,2,3\} \} \,$ of~$\mathfrak{u}(2)$,%
 \footnote{Note that under the isomorphism $\, \mathfrak{u}(2) \cong
 \mathbb{R}^4 \,$ given by mapping the basis $\, \{ \, \sigma_a^{}/2i \; | \;
 a \smin \{0,1,2,3\} \}$ \linebreak of~$\mathfrak{u}(2)$ to the standard basis
 of~$\mathbb{R}^4$, where $\sigma_0^{}$ is the unit matrix and $\sigma_1^{}$,
 $\sigma_2^{}$, $\sigma_3^{}$ are the Pauli matrices as before, the commutator
 in~$\mathfrak{u}(2)$ corresponds to that induced by the (associative)
 quaternion product in~$\mathbb{R}^4$ whereas the invariant scalar product
 in~$\mathfrak{u}(2)$ given by $\, (X,Y) = - 2 \, \mathrm{tr} (XY) \,$
 corresponds to the standard scalar product in~$\mathbb{R}^4$. With respect
 to this orthonormal basis, the (totally covariant) structure constants for
 $\mathfrak{u}(2)$ are again given by the components of the $\epsilon$-tensor,
 which are those of the standard volume form on~$\mathbb{R}^3$: all structure
 constants for which one of the indices takes the value~$0$ vanish.}
 and working in components, we see that $A$ is represented by a quartet of
 $1$-forms $A^0$, $A^1$, $A^2$, $A^3$ and $B$ by a triplet of $2$-forms
 $B_1^{}$, $B_2^{}$, $B_3^{}$ on~$P$ such that
 \[
  B_a^{}~=~{\textstyle \frac{1}{2}} \, \epsilon_{abc}^{} \,
           A^b \,\smwedge\, A^c~.
 \]
 Moreover, let $\phi$ be a given section of the truncated adjoint bundle
 $\, P \times_{U(2)} \mathbb{R}^3 \,$ which is an ``immersion'' in the sense
 that its covariant derivative $\DA \phi$ with respect to the given connection,
 which is a $1$-form on~$M$ taking values in $\, P \times_{U(2)} \mathbb{R}^3$,
 when interpreted as a vector bundle homomorphism from~$TM$ to~$\, P \times%
 _{U(2)} \mathbb{R}^3$, is fiberwise injective. When pulled back to~$P$,
 $\phi$ corresponds to an equivariant $\mathbb{R}^3$-valued function
 $\, \pi^* \phi \,$ on~$P$ and $\DA \phi$ corresponds to an equivariant
 horizontal $\mathbb{R}^3$-valued $1$-form $\, \pi^* \DA \phi = \DA \,
 \pi^* \phi \,$ on~$P$ which, once again, is an ``immersion'' in the sense
 that, when interpreted as an equivariant vector bundle homomorphism from~%
 $\pi^*(TM)$ to~$\, P \times \mathbb{R}^3$, it is fiberwise injective.
 In components, $\pi^* \phi \,$ is represented by a triplet of functions
 $\, \pi^* \phi^1$, $\pi^* \phi^{\,2}$, $\pi^* \phi^{\,3} \,$ and
 $\, \DA \, \pi^* \phi \,$ by a triplet of horizontal $1$-forms
 $\, \DA \, \pi^* \phi^1$, $\DA \, \pi^* \phi^{\,2}$,
 $\DA \, \pi^* \phi^{\,3} \,$ on~$P$, where
 \[
  \DA \, \pi^* \phi^{\,a}~
  =~d \, \pi^* \phi^{\,a} \, + \, \epsilon_{\phantom{\,a}bc}^{\,a} \,
    A^b \, \pi^* \phi^{\,c}~.
 \]
 By standard results~\cite{GHV,KN}, the covariant exterior derivative
 of~$\, \DA \, \pi^* \phi \,$, given by
 \[
  d_A^{} \, \DA \, \pi^* \phi^{\,a}~
  =~d \, \DA \, \pi^* \phi^{\,a} \, + \, \epsilon_{\phantom{\,a}bc}^{\,a} \;
    A^b \,\smwedge\,  \DA \, \pi^* \phi^{\,c}~,
 \]
 can also be written in the form
 \[
  d_A^{} \, \DA \, \pi^* \phi^{\,a}~
  =~\epsilon_{\phantom{\,a}bc}^{\,a} \, F^b \, \pi^* \phi^{\,c}~,
 \]
 where the $F^a$ are the components of the $\mathfrak{su}(2)$-part of the
 curvature form of the connection form $A$, defined by
 \[
  F^0~=~dA^0~~~,~~~
  F^a~=~dA^a \, + \, {\textstyle \frac{1}{2}} \,
        \epsilon_{\phantom{\,a}bc}^{\,a} \; A^b \>\smwedge\, A^c~,
 \]
 and there is a product rule relating the ordinary exterior derivative to
 the covariant exterior derivative which in the case of importance here reads
 \begin{eqnarray*}
 \lefteqn{d \, \bigl( \epsilon_{abc}^{} \; \pi^* \phi^{\,a} \;
                      \DA \, \pi^* \phi^{\,b} \,\smwedge\,
                      \DA \, \pi^* \phi^{\,c} \bigr)} \hspace*{2cm} \\[2mm]
  &=&\!\! \epsilon_{abc}^{} \; D_A^{} \, \pi^* \phi^{\,a} \,\smwedge\,
          \DA \, \pi^* \phi^{\,b} \,\smwedge\, \DA \, \pi^* \phi^{\,c} \\[1mm]
  & &\!\! \mbox{} + \, \epsilon_{abc}^{} \; \pi^* \phi^{\,a} \;
          d_A^{} \, \DA \, \pi^* \phi^{\,b} \,\smwedge\,
          \DA \, \pi^* \phi^{\,c} \\[1mm]
  & &\!\! \mbox{} - \, \epsilon_{abc}^{} \; \pi^* \phi^{\,a} \;
          \DA \, \pi^* \phi^{\,b} \,\smwedge\,
          d_A^{} \, \DA \, \pi^* \phi^{\,c}~.
 \end{eqnarray*}
 In fact, this expression vanishes identically since it represents an invariant
 horizontal $3$-form on~$P$ which corresponds to a $3$-form on~$M$, but there
 exists no non-zero $3$-form on a two-dimensional manifold. Now define
 \begin{equation} \label{eq:MSFEX1}
  \omega~=~{\textstyle \frac{1}{2}} \, \epsilon_{abc}^{} \;
           A^a \,\smwedge\, A^b \,\smwedge\, \DA \, \pi^* \phi^{\,c} \, - \,
           {\textstyle \frac{1}{2}} \, \epsilon_{abc}^{} \;
           \pi^* \phi^{\,a} \, A^0 \,\smwedge\,
           \DA \, \pi^* \phi^{\,b} \,\smwedge\, \DA \, \pi^* \phi^{\,c}~.
 \end{equation}
 Then using the above relations, together with the formula
 $\; \epsilon_{abc}^{} \, \epsilon_{\phantom{\,a}de}^{\,a} =
 \delta_{bd}^{} \, \delta_{ce}^{} \, - \,
 \delta_{be}^{} \, \delta_{cd}^{} \,$,
 we can calculate the exterior derivative of~$\,\omega$:
 \begin{eqnarray*}
  d\omega \!\!
  &=&\!\! {\textstyle \frac{1}{2}} \, \epsilon_{abc}^{} \;
          dA^a \,\smwedge\, A^b \,\smwedge\, \DA \, \pi^* \phi^{\,c} \, - \,
          {\textstyle \frac{1}{2}} \, \epsilon_{abc}^{} \;
          A^a \,\smwedge\, dA^b \,\smwedge\, \DA \, \pi^* \phi^{\,c} \\[1mm]
  & &\!\! \mbox{} + \, {\textstyle \frac{1}{2}} \, \epsilon_{abc}^{} \;
          A^a \,\smwedge\, A^b \,\smwedge\, d_A^{} \, \DA \, \pi^* \phi^{\,c}
          \, - \, {\textstyle \frac{1}{2}} \, \epsilon_{abc}^{} \,
          \epsilon_{\phantom{\,c}de}^{\,c} \; A^a \,\smwedge\,
          A^b \,\smwedge\, A^d \,\smwedge\, \DA \, \pi^* \phi^{\,e} \\[1mm]
  & &\!\! \mbox{} - \, {\textstyle \frac{1}{2}} \, \epsilon_{abc}^{} \;   
          \pi^* \phi^{\,a} \; dA^0 \,\smwedge\,
          \DA \, \pi^* \phi^{\,b} \,\smwedge\, \DA \, \pi^* \phi^{\,c} \\[1mm]
  & &\!\! \mbox{} + \, {\textstyle \frac{1}{2}} \, A^0 \,\smwedge\,
          d \, \bigl( \epsilon_{abc}^{} \; \pi^* \phi^{\,a} \;
                      \DA \, \pi^* \phi^{\,b} \,\smwedge\,
                      \DA \, \pi^* \phi^{\,c} \bigr) \\[2mm]
  &=&\!\! {\textstyle \frac{1}{2}} \, \epsilon_{abc}^{} \;
          F^a \,\smwedge\, A^b \,\smwedge\, \DA \, \pi^* \phi^{\,c} \, - \,
          {\textstyle \frac{1}{2}} \, \epsilon_{abc}^{} \;
          A^a \,\smwedge\, F^b \,\smwedge\, \DA \, \pi^* \phi^{\,c} \\[1mm]
  & &\!\! \mbox{} - \, {\textstyle \frac{1}{4}} \, \epsilon_{abc}^{} \,
          \epsilon_{\phantom{\,a}de}^{\,a} \; A^d \,\smwedge\,
          A^e \,\smwedge\, A^b \,\smwedge\, \DA \, \pi^* \phi^{\,c} \, + \,
          {\textstyle \frac{1}{4}} \, \epsilon_{abc}^{} \,
          \epsilon_{\phantom{\,b}de}^{\,b} \; A^a \,\smwedge\,
          A^d \,\smwedge\, A^e \,\smwedge\, \DA \, \pi^* \phi^{\,c} \\[1mm]
  & &\!\! \mbox{} + \, {\textstyle \frac{1}{2}} \, \epsilon_{abc}^{} \,
          \epsilon_{\phantom{\,c}de}^{\,c} \; 
          \pi^* \phi^{\,e} \, A^a \,\smwedge\, A^b \,\smwedge\, F^d \\[1mm]
  & &\!\! \mbox{} - \, {\textstyle \frac{1}{2}} \,\epsilon_{abc}^{} \,
          \epsilon_{\phantom{\,c}de}^{\,c} \; A^a \,\smwedge\,
          A^b \,\smwedge\, A^d \,\smwedge\, \DA \, \pi^* \phi^{\,e} \\[1mm]
  & &\!\! \mbox{} - \, {\textstyle \frac{1}{2}} \, \epsilon_{abc}^{} \;   
          \pi^* \phi^{\,a} \; F^0 \,\smwedge\,
          \DA \, \pi^* \phi^{\,b} \,\smwedge\, \DA \, \pi^* \phi^{\,c}~.
 \end{eqnarray*}
 The terms in the first and fifth line of the second equation vanish by the
 same argument of horizontality as before, while the terms in the second
 and fourth line vanish by symmetry. However, the term in the third line
 survives, i.e., we have
 \[
  d\omega~=~\delta_{ac}^{} \, \delta_{bd}^{} \; \pi^* \phi^{\,a} \;
            A^b \,\smwedge\, A^c \,\smwedge\, F^d~.
 \]
 This means that $\,\omega$ is closed if (and also only if) the
 $\mathfrak{su}(2)$-part of the curvature form $F$ vanishes,
 whereas there is no restriction on its $\mathfrak{u}(1)$-part.
 Moreover, it is easily verified that the image of~$\pi^*(TM)$ under
 $\DA \, \pi^* \phi$ together with the (one-dimensional) orthogonal
 complement of $P \times \mathbb{R}^3$ in $P \times \mathbb{R}^4$,
 when transferred from $P \times \mathbb{R}^4$ to $V\!P$ by means
 of the canonical isomorphism that exists between the two since $P$
 is a principal bundle, constitutes a three-dimensional invariant
 distribution~$L$ on~$P$ which is isotropic under $\,\omega$,
 \linebreak so provided $\,\omega$ is closed, $\omega$ is
 multisymplectic and~$L$ is the corresponding multilagrangian
 distribution. Of course, $L$ is not involutive.
\end{exem}

\begin{rmk}~
 We emphasize again that the input data for the construction of the poly\-%
 symplectic and multisymplectic forms in the last two examples allow for
 a natural interpretation in terms of Yang-Mills-Higgs field theory.
 Indeed, in the polysymplectic case, $A$ is an arbitrary $SU(2)$-%
 connection and so we only have to make one assumption that is more
 restrictive than stated in Example~\ref{exe:PSB-LNINV}: namely that
 the $1$-form $e$ is (covariantly) holonomous, that is, of the form
 $\, e = \DA \phi \,$ where the Higgs field $\phi$ is a section of
 the adjoint bundle $\, P \times_{SU(2)} \mathbb{R}^3$. Moreover, since
 only the vertical part of $A$ appears in equation~(\ref{eq:PSFEX4}),
 it is clear that $\,\fp$ is in fact completely independent of the
 choice of connection! In the multisymplectic case, we can argue that
 the relevant input data are (a)~a principal $U(1)$-bundle $P_0^{}$
 over~$M$ (or equivalently, a complex line bundle $L_0^{}$ over~$M$
 with fixed hermitean fiber metric) and a given $U(1)$-connection
 form $A_0^{}$ on~$P_0^{}$ (or equivalently, a given linear connection
 in~$L_0^{}$ preserving this fiber metric) and (b)~a given embedding
 $\phi$ of~$M$ into~$\mathbb{R}^3$. Indeed, the principal $U(2)$-bundle
 $P$ over~$M$ (corresponding to a rank $2$ complex vector bundle $E_0^{}$
 over~$M$ with fixed hermitean fiber metric) and the $U(2)$-connection
 form $A$ on~$P$ (corresponding to a linear connection in~$E_0^{}$ pre%
 serving this fiber metric) that appear in Example~\ref{exe:MSB-LNINV}
 are then obtained by the process of extension of structure group,
 using the embedding homomorphism of~$U(1)$ into~$U(2)$ that takes
 a phase to that phase times the unit $(2 \times 2)$-matrix; moreover,
 this construction guarantees that the $\mathfrak{su}(2)$-part of the
 curvature form $F$ will vanish. Conversely, if this condition is
 satisfied, then according to the Ambrose-Singer theorem, $P$ can
 be reduced to the holonomy bundle of~$A$, whose structure group is
 the holonomy group of~$A$ which has connected one-component $U(1)$:
 this means that except for a possible discrete part, we are back to
 the previous situation. It should also be emphasized that from a
 purely topological point of view, reduction from $U(2)$ to $U(1)$
 is always possible, e.g., over compact Riemann surfaces~\cite%
 [Theorem~10 \& Corollary, p.~63]{Gun}, so the restriction made
 here really only concerns the connection form~$A$, not the bundle~$P$.
 Moreover, in this situation, the truncated adjoint bundle $\, P \times_{U(2)}
 \mathbb{R}^3 \,$ will be trivial and so the Higgs field $\phi$ with its
 injective (covariant = ordinary) derivative will provide an immersion
 of~$M$ into~$\mathbb{R}^3$: it is only slightly more restrictive to
 assume that this immersion is in fact an embedding. Finally, an
 important observation that applies to both settings is that, as
 already mentioned, the polysymplectic form $\,\fp$ of Example~%
 \ref{exe:PSB-LNINV} and the multisymplectic form $\,\omega$ of
 Example~\ref{exe:MSB-LNINV} are both invariant under the action
 of the respective structure groups $SU(2)$ and $U(2)$ on the total
 space~$P$: expressed in physical language, this means that they are
 \emph{gauge invariant}.
\end{rmk}

We expect it to be possible to apply similar procedures, with $SU(2)$
replaced by some compact three-dimensional Lie group which is solvable
and $U(2)$ replaced by a compact four-dimensional Lie group which is a
suitable one-dimensional extension thereof, to construct polysymplectic
and multisymplectic fiber bundles $P$ where $P$ is still a compact
manifold but $L$ is now involutive. However, in our view the real
challenge is to come up with examples of polysymplectic and multi%
symplectic fiber bundles that appear naturally as multiphase spaces of
physically realistic covariant hamiltonian field theories but, just as
in the above cases, cannot be obtained by taking the affine dual of the
first order jet bundle of some configuration bundle. This question is
presently under investigation.


\section{The Darboux Theorem}
\label{sec:DT}

Now we are able to prove the Darboux theorem for poly- and multilagrangian
forms. Here, the specific algebraic structure of poly- and multilagrangian
subspaces identified in the first two sections turns out to be crucial,
in the sense that this central theorem can, in all cases, be proved by
appropriately adapting the procedure used to prove the classical Darboux
theorem for symplectic forms. (See, for instance, Ref.~\cite{AM}).
\begin{thm}[Polylagrangian Darboux theorem]~
 \label{thm:DARPG}
 Let $P$ be a polylagrangian \linebreak fiber bundle over an $n$-dimensional
 manifold~$M$ with polylagrangian $(k+1)$-form $\,\fp$ of rank~$N$ taking
 values in a fixed $\hat{n}$-dimensional vector bundle $\hat{T}$ over the
 same manifold~$M$ and with polylagrangian distribution $L$, which is
 assumed to be involutive (recall this is automatic if $\, \hat{n}
 \geqslant 3$), and let $\, \{\, \hat{e}_a \,|\; 1 \leqslant a
 \leqslant \hat{n} \,\} \,$ be a basis of local sections of\/~$\hat{T}$.
 Then around any point of~$P$ (within the domain of the given basis of
 local sections), there exists a system of local coordinates $(x^\mu,q^i,
 p\>\!_{i_1 \ldots\, i_k}^a,r^\kappa)$ $(1 \leqslant \mu \leqslant n$,
 $1 \leqslant a \leqslant \hat{n}$, $1 \leqslant i \leqslant N$,
 $1 \leqslant i_1^{} < \ldots < i_k^{} \leqslant N$, $1 \leqslant
 \kappa \leqslant \dim \, \ker \, \fp)$, called \textbf{Darboux
 coordinates} or \textbf{canonical coordinates}, such that
 \begin{equation} \label{eq:FPLBDG}
  \fp~=~{\textstyle \frac{1}{k!}}~
        dp\>\!_{i_1 \ldots\, i_k}^a \,\smwedge\,
        dq^{i_1} \,\smwedge \ldots\, \smwedge\, dq^{i_k} \,\smotimes\,
        \hat{e}_a~,
 \end{equation}
 and such that (locally) $L$ is spanned by the vector fields
 $\, \partial/\partial p\>\!_{i_1 \ldots\, i_k}^a \,$ and
 $\, \partial/\partial r^\kappa \,$ while \linebreak $\ker \, \fp \,$
 is spanned by the vector fields $\, \partial/\partial r^\kappa$.
 In the polysymplectic case this expansion takes the form
 \begin{equation} \label{eq:FPSBDG}
  \fp~=~\bigl( dq^i \,\smwedge\, dp\>\!_i^a \bigr) \,\smotimes\, \hat{e}_a~.
 \end{equation}
\end{thm}
Similarly we have
\begin{thm}[Multilagrangian Darboux theorem]~ \label{thm:DARMG}
 Let $P$ be a multilagrangian fiber bundle over an $n$-dimensional manifold~$M$
 with multilagrangian $(k+1)$-form $\,\omega$ of rank~$N$ and horizontality
 degree $k+1-r$, where $\, 1 \leqslant r \leqslant k \,$ and $\, k+1-r
 \leqslant n$, and with multilagrangian distribution $L$, which is assumed
 to be involutive (recall this is automatic if $\, {n \choose k+1-r}
 \geqslant 3$). Then around any point of~$P$, there exists a system of local
 coordinates $(x^\mu,q^i,p\>\!_{i_1 \ldots\, i_s;\,\mu_1 \ldots\,\mu_{k-s}}^{},
 r^\kappa)$ $(0 \leqslant s \leqslant r-1$,  $1 \leqslant \mu \leqslant n$,
 $1 \leqslant i \leqslant N$, $1 \leqslant i_1^{} < \ldots < i_s^{}
 \leqslant N$, $1 \leqslant \mu_1^{} < \ldots < \mu_{k-s}^{} \leqslant N$,
 $1 \leqslant \kappa \leqslant \dim \,  \ker \, \omega)$, called
 \textbf{Darboux coordinates} or \textbf{canonical coordinates},
 such that
 \begin{equation} \label{eq:FMLBDG}
  \omega~=~\sum_{s=0}^{r-1} \, {\textstyle \frac{1}{s!} \, \frac{1}{(k-s)!}}~
           dp\>\!_{i_1 \ldots\, i_s;\,\mu_1 \ldots\, \mu_{k-s}}^{} \,\smwedge\,
           dq^{i_1} \,\smwedge \ldots\, \smwedge\, dq^{i_s} \,\smwedge\,
           dx^{\mu_1} \,\smwedge \ldots\, \smwedge\, dx^{\mu_{k-s}}~,
 \end{equation}
 and such that (locally) $L$ is spanned by the vector fields
 $\, \partial/\partial p\>\!_{i_1 \ldots\, i_s;\,\mu_1 \ldots\, \mu_{k-s}}^{} \,$
 and $\, \partial/\partial r^\kappa$ \linebreak while $\ker \, \omega \,$ is
 spanned by the vector fields $\, \partial/\partial r^\kappa$. In these
 coordinates, its symbol is given by
 \begin{eqnarray} \label{eq:FMLBDS}
  \fp \;
  = \; {\textstyle \frac{1}{(r-1)!} \, \frac{1}{(k+1-r)!}} \;
       \bigl( dp\>\!_{i_1 \ldots\, i_{r-1};\,\mu_1 \ldots\, \mu_{k+1-r}}^{}
              \,\smwedge\, dq^{i_1} \,\smwedge \!\ldots \smwedge\,
                           dq^{i_{r-1}} \bigr) \smotimes
       \bigl( dx^{\mu_1} \,\smwedge \ldots \smwedge\,
       dx^{\mu_{k+1-r}} \bigr) .~
 \end{eqnarray}
 In the multisymplectic case these expansions take the form
 \begin{equation} \label{eq:FMSBDG}
  \omega~=~dq^i \,\smwedge\, dp\>\!_i^\mu \,\smwedge\, d^{\,n} x_\mu^{} \, - \,
   dp \;\smwedge\, d^{\,n} x~,
 \end{equation}
 and
 \begin{equation} \label{eq:FMSBDS}
  \fp~=~\bigl( dq^i \,\smwedge\, dp\>\!_i^\mu \bigr) \,\smotimes\,
        d^{\,n} x_\mu^{}~,
 \end{equation}
 where
 \begin{equation}
   d^{\,n} x_\mu^{}~=~\mathrm{i}_{\partial_\mu^{}} \, d^{\,n} x~
   =~\frac{1}{(n-1)!} \; \epsilon_{\mu \mu_1 \ldots\, \mu_{n-1}}^{} \,
     dx^{\mu_1} \,\smwedge\, \ldots\, \smwedge\, dx^{\mu_{n-1}}~.
 \end{equation}
\end{thm}
\proof~
 For the sake of simplicity, we concentrate on the multilagrangian case:
 the proof for the other case is entirely analogous, requiring only small
 and rather obvious modifications.

 \noindent
 Due to the local character of this theorem and since the kernel of $\omega$,
 the multilagrangian subbundle $L\,$ and the vertical bundle $V\!P$ are all
 involutive, with $\ker \omega \smsubset L \smsubset V\!P \,$, we can
 without loss of generality work in a local chart of the manifold~$P$
 around the chosen reference point in which the corresponding foliations
 are ``straightened out'', so we may assume that $P \cong \mathbb{R}^n
 \oplus \mathbb{R}^N \oplus L_0^{} \oplus K_0^{} \,$ with $\, V\!P \cong
 \mathbb{R}^N \oplus L_0^{} \oplus K_0^{}$, $L \cong L_0^{} \oplus K_0^{} \,$
 and $\, \ker \omega \cong K_0^{} \,$ with fixed subspaces $L_0^{}$ and
 $K_0^{}$ and such that the aforementioned reference point corresponds to
 the origin. We also take $\omega_0^{}$ to be the constant multilagrangian
 form, with multilagrangian distribution $L$, obtained by spreading
 $\omega(0)$, the value of the multilagrangian form $\omega$ at the
 origin, all over the space $P$; then the existence of canonical
 coordinates for $\omega_0$, in the form given by equation~%
 (\ref{eq:FMLBDG}), is already guaranteed by the algebraic
 Darboux theorem of the previous chapter (Theorem~\ref{thm:DARMA}).
 \\[2mm]
 Now consider the family of $(k+1)$-forms given by $\, \omega_t^{} =
 \omega_0^{} + t(\omega-\omega_0^{})$, for every $\, t \smin \mathbb{R}$.
 Obvious\-ly, $\omega_t^{}(0) = \omega_0^{} \,$ for every $\, t \smin
 \mathbb{R}$, which is non-degenerate on $\, K_0' = \mathbb{R}^n \oplus
 \mathbb{R}^N \oplus L_0^{} \,$ (a complement of $K_0^{}$ in~$P$).
 Since non-degeneracy is an open condition, and using a compactness
 argument with respect to the parameter $t$, we conclude that there
 is an open neighborhood of $0$ such that, for all $t$ satisfying
 $\, 0 \leqslant t \leqslant 1 \,$ and all points $p$ in this neighborhood,
 $\omega_t^{}(p)$ is non-degenerate on $\, K_0' = \mathbb{R}^n \oplus
 \mathbb{R}^N \oplus L_0^{} \,$, that is, its kernel equals~$K_0^{}$.
 Moreover, for all $t$ satisfying $\, 0 \leqslant t \leqslant 1 \,$
 and all points $p$ in this neighborhood, the subspace $L_0^{}$, being
 isotropic for $\omega_0^{}$ as well as for $\omega(p)$, is also isotropic
 for $\omega_t^{}(p)$ and, since it contains the kernel of~$\omega_t^{}(p)$
 and has the right dimension as given by equation~(\ref{eq:MULTD1}), is
 even multilagrangian for $\omega_t^{}(p)$, according to Proposition~%
 \ref{prp:MULTLS2}. On the other hand, we have $\, d\omega_0^{} = 0$
 (trivially) and $\, d\omega = 0$ (by hypothesis), so we can apply
 an appropriate version of the Poincar\'e lemma (see Appendix~B) to
 prove, in some open neighborhood of the point~$0$ in~$P$ (contained
 in the previous one), existence of a $k$-form $\alpha$ satisfying
 $\, d\alpha = \omega_0^{} - \omega \,$ and $\, \alpha^\flat(L) = 0$.
 Now take $X_t^{}$ to be the unique time dependent vector field on~$P$
 taking values in $L_0^{}$\,%
 \footnote{It is at this point that we make essential use of the
 hypothesis that $L_0^{}$ is multilagrangian and not just isotropic
 or even maximal isotropic (with respect to $\omega_t^{}(p)$, in
 this case).}
 defined by
 \[
  \mathrm{i}_{X_t}^{} \omega_t^{}~=~\alpha~.
 \]
 Let $\, F_t^{} \equiv F_{(0,t)}$ be its flux beginning at $0$, once again
 defined, for $\, 0 \leqslant t \leqslant 1$, in some open neighborhood of
 the point $0$ in~$P$ (contained in the previous one). Then it follows that
 \[
 \begin{array}{rcl}
  {\displaystyle \frac{d}{ds} \, \bigg|_{s=t} F_s^* \omega_s^{}}
  &=& {\displaystyle F_t^* \left( \frac{d}{ds} \bigg|_{s=t} \omega_s^{} \right) +
                     \frac{d}{ds} \bigg|_{s=t} F_s^* \omega_t^{}} \\[5mm]
  &=& F_t^* \bigl( \omega - \omega_0^{} + L_{X_t}^{} \omega_t^{} \bigr) \\[3mm]
  &=& F_t^* \bigl( \omega - \omega_0^{} +
                   d(\mathrm{i}_{X_t}^{} \omega_t^{}) \bigr) \\[3mm]
  &=& F_t^* \bigl( \omega - \omega_0^{} + d\alpha \bigr) \\[3mm]
  &=& 0
 \end{array}
 \]
 Therefore, $F_1^{}$ is the desired coordinate transformation, since
 $\, F_1^* \omega = F_1^* \omega_1^{} = F_0^*\omega_0^{} = \omega_0^{}$.
 \hspace*{\fill}
\qed


\section{Conclusions and Outlook}

In this paper we have presented a concise definition of a new class
of geometric structures which we propose to call polylagrangian or
multilagrangian structures and which include as a special case the
familiar polysymplectic or multisymplectic structures encountered in
the hamiltonian formulation of classical field theory. All of them
are defined in terms of differential forms satisfying an algebraic
condition that amounts to postulating the existence of an isotropic
subbundle which is ``sufficiently large''~-- a condition which, when
combined with the standard integrability conditions that the pertinent
differential form should be closed and the aforementioned ``sufficiently
large'' isotropic subbundle should be involutive (only in two space-time
dimensions is this not automatic), allows to derive a Darboux theorem
assuring the existence of ``canonical'' local coordinates around each
point. Another characteristic feature of all these structures is that
they are naturally defined on the total spaces of fiber bundles whose
base space is interpreted as the space-time manifold of field theory.
Moreover, there is a standard class of examples defined by bundles of
partially horizontal forms over the total space of another fiber bundle
(the so-called configuration bundle), which includes the multiphase
spaces of interest in physics and is the analogue of cotangent bundles
of manifolds as a standard class of examples of symplectic manifolds.
To our knowledge, this is the first example of a natural geometric
structure with important physical applications that is defined by a
differential form (or even a tensor field) of degree strictly larger
than $2$ and strictly smaller than the dimension of the underlying
manifold.

Starting from this basis, there is a number of rather obvious questions
that arise, most of which are closely interrelated. One of them concerns
the structure of the under\-lying Lie group: interpreting these structures
as $G$-structures, what is the nature of~$G$? Other questions refer to
the definition of Poisson brackets (see the discussion in Refs~%
\cite{FR1,FPR1,FPR2,FR2}), the definition of actions of Lie groups
and, more generally, of Lie groupoids on poly\-symplectic/polylagrangian
or multisymplectic/multilagrangian fiber bundles, the construction of a
corresponding momentum map (which would provide a general framework for
the construction of Noether currents and the energy-momentum tensor
within a direct and manifestly covariant hamiltonian approach), the
formulation of a Marsden-Weinstein reduction procedure and, last but by
no means least, the explicit construction of other classes of examples,
in particular, analogues of the coadjoint orbit construction of
symplectic geometry. All these problems are completely open and
certainly will be the subject of much research in the future.


\begin{appendix}

\section{Some Counterexamples}

In this appendix, we wish to provide additional evidence for the conceptual
simplicity and usefulness of our definition of polylagrangian and multi%
lagrangian structures by investigating, in a purely algebraic setting
and for the simplest case of vector-valued $2$-forms, various structural
properties of polysymplectic forms that can be introduced directly for
general vector-valued $2$-forms. By constructing explicit counterexamples,
we will show, however, that none of them is sufficiently strong to replace
the requirement of existence of a polylagrangian subspace: this condition
must therefore be imposed separately and is then sufficient to imply all the
others, so that~-- in contrast to what is done, e.g., in Ref.~\cite{Gr}~--
we have refrained from including any of them into our definition of a
poly\-symplectic structure.

Suppose, as in Section 1, that $V$ and $\hat{T}$ are finite-dimensional
real vector spaces, with $\, \hat{n} \equiv \dim \hat{T}$, and assume
that $\; \fp \smin\, \bwedge^{\!2\,} V^* \otimes\, \hat{T} \,$ is an
arbitrary $\hat{T}$-valued $2$-form on~$V$. Given any linear form
$\, \hat{t}^* \smin \hat{T}^* \,$ on $\hat{T}$, we consider the
projection $\; \omega_{\hat{t}^*}^{} = \langle \hat{t}^* , \fp
\rangle \,$ of $\,\fp$ along $\hat{t}^*$ which is an ordinary
$2$-form on~$V$, and we define its \textbf{rank} to be equal to
half the dimension of its \textbf{support}, which in turn can be
defined as the annihilator of its kernel~\cite{Mat}:\,%
\footnote{Thus our definition of rank differs from that
Ref.~\cite{Mat} by a factor of 2.}
\begin{equation}
 \mathrm{rk} \bigl( \omega_{\hat{t}^*}^{} \bigr)~
 =~{\textstyle \frac{1}{2}} \; \dim \, \mathrm{supp} \; \omega_{\hat{t}^*}^{}
 =~{\textstyle \frac{1}{2}} \,
   \bigl( \dim V \, - \, \dim \, \ker \, \omega_{\hat{t}^*}^{} \bigr)~.
\end{equation}
Now note that the linear mapping
\begin{equation} \label{eq:OMMON1} 
 \begin{array}{ccc}
  \hat{T}^* & \longrightarrow & \bwedge^{\!2\,}  V^* \\[1mm]
  \hat{t}^* &   \longmapsto   & \omega_{\hat{t}^*}
 \end{array}
\end{equation}
induces, for every integer $\, k \geqslant 1$, a canonically defined
linear mapping
\begin{equation} \label{eq:OMPOL1} 
 \begin{array}{ccc}
  \bvee^{\,k\,} \hat{T}^* & \longrightarrow & \bwedge^{\!2k\,} V^* \\[1mm]
             P            &   \longmapsto   &        P(\fp)
 \end{array}
\end{equation}
where we have identified the space $\, \bvee^{\,k\,} \hat{T}^* \,$ of covariant
symmetric tensors of degree $k$ over~$\hat{T}$ with the space of homogeneous
polynomials $P$ of degree $k$ on~$\hat{T}$. Explicitly, in terms of a basis
$\, \{\, \hat{e}_a^{} \,|\; 1 \leqslant a \leqslant \hat{n} \,\} \,$ of~%
$\hat{T}$, with dual basis $\, \{\, \hat{e}^a \,|\; 1 \leqslant a \leqslant
\hat{n} \,\} \,$ of~$\hat{T}^*$, we write $\, \omega^a = \omega_{\hat{e}^a}^{}
= \langle \hat{e}^a , \fp \rangle$ $(1 \leqslant a \leqslant \hat{n}) \,$
and obtain
\begin{equation} \label{eq:OMPOL2} 
 P~=~P_{a_1 \ldots\, a_k} \;
     \hat{e}^{a_1} \,\smvee \ldots\, \smvee\, \hat{e}^{a_k}
 \quad \Longrightarrow \quad
 P(\fp)~=~P_{a_1 \ldots\, a_k} \;
          \omega^{a_1} \,\smwedge \ldots\, \smwedge\, \omega^{a_k}~.
\end{equation}
This allows us to introduce the following terminology:
\begin{defi}~ \label{def:PSPLF2}
 Let $V$ and $\hat{T}$ be finite-dimensional vector spaces and let $\,\fp$ be
 a $\hat{T}$-valued $2$-form on~$V$. We say that $\,\fp$ has \textbf{constant
 rank}~$N$ if\/ $\, \mathrm{rk}(\omega_{\hat{t}^*}) = N \,$ for every
 $\, \hat{t}^* \smin\, \hat{T}^* \setminus \{0 \}$ \linebreak and that
 $\,\fp$ has \textbf{uniform rank}~$N$ if the linear mapping~(\ref{eq:OMPOL1})
 is injective for $\, k=N \,$ and identically zero for $\, k = N+1$.
\end{defi}
Using multi-indices $\, \alpha = (\alpha_1^{},\ldots,\alpha_{\hat{n}}^{})
\smin\, \mathbb{N}^{\hat{n}}$, we set
\[
 \hat{e}^\alpha~=~(\hat{e}^1)^{\alpha_1} \,\smvee \ldots\, \smvee\,
                  (\hat{e}^{\hat{n}})^{\alpha_{\hat{n}}} 
 \qquad \mbox{where} \qquad
 (\hat{e}^a)^{\alpha_a}~=~\hat{e}^a \,\smvee \ldots\, \smvee\, \hat{e}^a~~~
 \mbox{($\alpha_a$ times)} 
\]
and
\[
 \omega^\alpha~=~(\omega^1)^{\alpha_1} \,\smwedge \ldots\, \smwedge\,
                 (\omega^{\hat{n}})^{\alpha_{\hat{n}}}
 \qquad \mbox{where} \qquad
 (\omega^a)^{\alpha_a}~=~\omega^a \,\smwedge \ldots\, \smwedge\, \omega^a~~~
 \mbox{($\alpha_a$ times)} 
\]
to rewrite equation~(\ref{eq:OMPOL2}) in the form
\begin{equation} \label{eq:OMPOL3} 
 P~=~\sum_{|\alpha|=k} \, P_\alpha \; \hat{e}^\alpha
 \quad \Longrightarrow \quad
 P(\fp)~=~\sum_{|\alpha|=k} \, P_\alpha \; \omega^\alpha~.
\end{equation}
Since $\, \{ \hat{e}^\alpha \, | \; |\alpha| = k \} \,$ is a basis of
$\bvee^{\,k\,} \hat{T}^*$, requiring $\fp$ to have uniform rank~$N$
amounts to imposing the following conditions:
\begin{equation} \label{eq:UNIFR1} 
 \begin{array}{c}
  \mbox{$\{ \omega^\alpha \, | \; |\alpha| = N \} \,$
        is linearly independent} \\[2mm]
  \mbox{$\omega^\alpha = 0 \,$ for $\, |\alpha| = N+1$}
 \end{array}~.
\end{equation}
It is in this form that the requirement of uniform rank appears in the
definition of a poly\-symplectic form adopted in Ref.~\cite{Gr}.

To gain a better understanding for the conditions of constant rank and
of uniform rank introduced above, we note first of all that they both
generalize the standard notion of rank for ordinary forms. Indeed, when
$\, \hat{n} = 1$, that is, given an ordinary $2$-form $\omega$ of rank~$N$
on~$V$, we can choose a canonical basis $\, \{ e_1^{},\ldots,e_N^{},f^1,
\ldots,f^N \} \,$ of a subspace of~$V$ complementary to the kernel of~%
$\omega$, with dual basis $\, \{ e^1,\ldots,e^N,f_1^{},\ldots,f_N^{} \} \,$
of the subspace~$\, \mathrm{supp} \; \omega \,$ of~$V^*$, to conclude that
$\; \omega = e^i \,\smwedge f_i^{} \;$ and therefore
\[
 \omega^N~=~\pm \; e^1 \smwedge \ldots \smwedge\, e^N \smwedge
 f_1 \,\smwedge \ldots \smwedge f_N~\neq~0~~~,~~~
 \omega^{N+1}~=~0~.
\]
In other words, the rank of~$\omega$ can be characterized as that positive
integer~$N$ for which $\; \omega^N \neq 0 \,$ but $\; \omega^{N+1} = 0$.
From this observation, it follows that,  in the general case considered
before, the requirement of uniform rank implies that of constant rank
because it guarantees that for every $\, \hat{t}^* \smin\, \hat{T}^*
\setminus \{0\}$, we have $\, \omega_{\hat{t}^*}^N \neq 0 \,$ and
$\, \omega_{\hat{t}^*}^{N+1} = 0 \,$, that is,
$\, \mathrm{rk}(\omega_{\hat{t}^*}) = N$. However,
the converse does not hold, as shown by the following
\begin{exem}~
 ($\, \hat{n} = 2 \,$, $\, N = 2 \,$, $\, \dim V = 4 \,$,
  $\, \ker \, \fp = \{0\} \,$)\,: \\[2mm]
 Let $\, V = \mathbb{R}^4 \,$, $\, \hat{T} = \mathbb{R}^2 \,$ and consider the
 $\mathbb{R}^2$-valued $2$-form $\fp$ built from the following two ordinary
 $2$-forms:
 \[
  \omega^1~=~dx \>\smwedge\, dy \, + \, du \>\smwedge\, dv~~~,~~~
  \omega^2~=~dx \>\smwedge\, du \, - \, dy \>\smwedge\, dv~.
 \]
 Then for $\; \hat{t}^* = \hat{t}_a^* \, \hat{e}^a \smin\, (\mathbb{R}^2)^*$,
 we have
 \[
  \omega_{\hat{t}^*}~=~\hat{t}_a^* \, \omega^a~
  =~dx \>\smwedge \bigl( \hat{t}_1^* \, dy \, + \, \hat{t}_2^* \, du \bigr)
    \, + \,
    dv \>\smwedge \bigl( \hat{t}_2^* \, dy \, - \, \hat{t}_1^* \, du \bigr)~.
 \]
 Thus we obtain, for every $\, \hat{t}^* \neq 0$,
 \[
  (\omega_{\hat{t}^*}^{})^2~
  \equiv~\omega_{\hat{t}^*}^{} \>\smwedge\, \omega_{\hat{t}^*}^{}~
  =~\bigl( (\hat{t}_1^*)^2 + (\hat{t}_2^*)^2 \bigr) \;
    dx \>\smwedge\, dy \>\smwedge\, du \>\smwedge\, dv~\neq~0~,
 \]
 whereas, due to the fact that we are in a four-dimensional space,
 \[
  (\omega_{\hat{t}^*}^{})^3~
  \equiv~\omega_{\hat{t}^*}^{} \>\smwedge\, \omega_{\hat{t}^*}^{}
                               \>\smwedge\, \omega_{\hat{t}^*}~=~0~,
 \]
 which guarantees that $\fp$ has constant rank~$2$. However, $\fp$ does not
 have uniform rank~$2$, since
 \[
  \omega^1 \smwedge\, \omega^2~=~0~.
 \]
\end{exem}
On the other hand, polysymplectic forms do have uniform rank:
\begin{prp}~ \label{prp:PSUNR} 
 Let $V$ and $\hat{T}$ be finite-dimensional vector spaces and let $\,\fp$ be
 a $\hat{T}$-valued polysymplectic form of rank~$N$ on~$V$. Then $\,\fp$ has
 uniform rank~$N$.
\end{prp}
\proof~
 Introducing a (polysymplectic) canonical basis in which
 \[
  \fp~=~\bigl( e_i^a \,\smwedge\, e^i \bigr) \,\smotimes\, \hat{e}_a~,
 \]
 or equivalently
 \[
  \omega^a~=~e_i^a \,\smwedge\, e^i \quad (1 \leqslant a \leqslant \hat{n})~,
 \]
 suppose that  $\, \alpha = (\alpha_1^{},\ldots,\alpha_{\hat{n}}^{}) \smin\,
 \mathbb{N}^{\hat{n}} \,$ is a multi-index of degree~$k$ (i.e., such
 that $\, |\alpha| = \alpha_1 + \ldots + \alpha_{\hat{n}} = k$) and
 consider the form
 \[
  \begin{array}{l}
   \omega^\alpha~
   =~\pm \; \bigl( (e^{i_1^1} \,\smwedge \ldots\, \smwedge\,
                    e^{i_{\alpha_1}^1} ) \,\smwedge \ldots\, \smwedge\,
                   (e^{i_1^{\hat{n}}} \,\smwedge \ldots\, \smwedge\,
                    e^{i_{\alpha_{\hat{n}}}^{\hat{n}}} ) \bigr) \\[3mm]
   \hspace*{5em} \,\smwedge\,
            \bigl( e_{i_1^1}^1 \,\smwedge \ldots\, \smwedge\,
                   e_{i_{\alpha_1}^1}^1 \bigr) \,\smwedge \ldots\, \smwedge
            \bigl( e_{i_1^{\hat{n}}}^{\hat{n}} \,\smwedge \ldots\, \smwedge\,
                   e_{i_{\alpha_{\hat{n}}}^{\hat{n}}}^{\hat{n}} \bigr)~.
  \end{array}
 \]
 Obviously any such form vanishes when $\, k = N+1 \,$ since it then contains
 an exterior product of ($N+1$) $1$-forms $e^i$ belonging to an $N$-dimensional
 subspace. On the other hand, all these forms are linearly independent when
 $\, k = N \,$ since $\omega^\alpha$ then contains the exterior product
 $\, e^1 \,\smwedge \ldots\, \smwedge\, e^N \,$ multiplied by the exterior
 product of $\alpha_1^{}$ $1$-forms of type $e_i^1$ with $\ldots$ \linebreak
 with $\alpha_{\hat{n}}^{}$ $1$-forms of type $e_i^{\hat{n}}$; thus
 $\omega^\alpha$ and $\omega^\beta$ belong to different subspaces
 of~$\bwedge^{2N} V^*$ whenever $\, \alpha \neq \beta$.
\qed

\noindent
The converse statement, as we shall see shortly, is remote from being true.
In fact, if it were true, then if $\, \hat{n} \geqslant 2$, it should be
possible to construct the polylagrangian subspace as the sum of the kernels
of the projected forms, as required by Theorem~\ref{thm:POLILS1}. Therefore,
it should be possible to show that the subspace defined as the sum of
these kernels is isotropic. And indeed, as a partial result in this
direction, we have the following
\begin{prp}~
 Let $V$ and $\hat{T}$ be finite-dimensional vector spaces and let $\,\fp$
 be a $\hat{T}$-valued $2$-form of uniform rank~$N$ on~$V$. Then for any
 $\, \hat{t}_1^*,\hat{t}_2^* \smin\, \hat{T}^* \setminus \{0\} \,$, the
 kernel of~$\omega_{\hat{t}_1^*}$ is isotropic with respect to~%
 $\omega_{\hat{t}_2^*}$.
\end{prp}
\proof~
 Given $\, u,v \smin \ker \, \omega_{\hat{t}_1^*}$, we have
 \[
  \mathrm{i}_u^{} \omega_{\hat{t}_1^*}^N~
  =~N \; \mathrm{i}_u \omega_{\hat{t}_1^*}^{} \,\smwedge\,
    \omega_{\hat{t}_1^*}^{N-1}~=~0~~~,~~~
  \mathrm{i}_v^{} \omega_{\hat{t}_1^*}^N~
  =~N \; \mathrm{i}_v \omega_{\hat{t}_1^*}^{} \,\smwedge\,
    \omega_{\hat{t}_1^*}^{N-1}~=~0~,
 \]
 and therefore
 \[
  \omega_{\hat{t}_2^*}^{}(u,v) \, \omega_{\hat{t}_1^*}^N~
  =~\mathrm{i}_u^{} \mathrm{i}_v^{}
    \bigl( \omega_{\hat{t}_2^*}^{} \,\smwedge\, \omega_{\hat{t}_1^*}^N \bigr)~
  =~0~.
 \]
 Using that $\; \omega_{\hat{t}_1^*}^N \neq 0$, it follows that
 $\; \omega_{\hat{t}_2^*}^{}(u,v) = 0$.
\qed

\noindent
However, isotropy of the subspace defined as the sum of the kernels of all the
projected forms, which is equivalent to the (stronger) condition that for any
$\, \hat{t}_1^*,\hat{t}_2^*,\hat{t}_3^* \smin\, \hat{T}^* \setminus \{0\} \,$,
$\ker \, \omega_{\hat{t}_1^*} \,$ and $\ker \, \omega_{\hat{t}_2^*} \,$ are
orthogonal under $\omega_{\hat{t}_3^*}$, i.e., that
\[
 \omega_{\hat{t}_3^*}^{}(u_1^{},u_2^{})~=~0
 \qquad \mbox{for $\, u_1^{} \smin\, \ker \, \omega_{\hat{t}_1^*} \,$ and
                  $\, u_2^{} \smin\, \ker \, \omega_{\hat{t}_2^*}$}~,
\]
cannot be derived from the condition of uniform rank. A nice counterexample
is obtained by choosing $V$ and $\hat{T}$ to be the same space, supposing it
to be a Lie algebra $\mathfrak{g}$ and defining $\fp$ to be the commutator
in~$\mathfrak{g}$. Then for $\, \hat{t}^* \smin \mathfrak{g}^*$, the kernel
$\, \ker \omega_{\hat{t}^*} \,$ and the support $\; \mathrm{supp} \;
\omega_{\hat{t}^*} \,$ of the projected form $\omega_{\hat{t}^*}$ are the
isotropy algebra of $\hat{t}^*$ and the tangent space to the coadjoint
orbit passing through $\hat{t}^*$, respectively. There is one and only
one semisimple Lie algebra for which $\fp$ has constant rank, since this
condition states that all coadjoint orbits except the trivial one, $\{0\}$,
should have the same dimension: this is the algebra of type $A_1$, that is,
$\, \mathbb{R}^3 \,$ equipped with the vector product $\times$.
\begin{exem}~ 
 ($\, \hat{n} = 3 \,$, $\, N = 1 \,$, $\, \dim V = 3 \,$,
  $\, \ker \, \fp = \{0\} \,$)\,: \\[2mm]
 Let $\, V = \hat{T} = \mathbb{R}^3 \,$ and consider the $\mathbb{R}^3$-%
 valued $2$-form $\fp$ built from the following three ordinary $2$-forms:
 \[
  \omega^1~=~dy \>\smwedge\, dz~~~,~~~
  \omega^2~=~dz \>\smwedge\, dx~~~,~~~
  \omega^3~=~dx \>\smwedge\, dy~.
 \]
 Then for $\; \hat{t}^* = \hat{t}_a^* \, \hat{e}^a \smin\, (\mathbb{R}^3)^*$,
 we have
 \[
  \omega_{\hat{t}^*}^{}~=~\hat{t}_a^* \, \omega^a~
                       =~\hat{t}_1^* \; dy \>\smwedge\, dz \, + \,
                         \hat{t}_2^* \; dz \>\smwedge\, dx \, + \,
                         \hat{t}_3^* \; dx \>\smwedge\, dy~.
 \]
 Obviously, $\omega^1$, $\omega^2$ and $\omega^3$ are linearly independent
 and hence $\fp$ has uniform rank~$1$, since there exists no non-zero
 $4$-form on a three-dimensional space. On the other hand, we have
 \[
  \ker \, \omega_{\hat{t}^*}^{}~
  =~\langle \, \hat{t}_1^* \; \frac{\partial}{\partial x} \, + \,
               \hat{t}_2^* \; \frac{\partial}{\partial y} \, + \,
               \hat{t}_3^* \; \frac{\partial}{\partial z} \, \rangle~,
 \]
 Therefore, the intersection of the three kernels $\, \ker \, \omega^1$,
 $\ker \, \omega^2 \,$ and $\, \ker \, \omega^3 \,$ is $\{0\}$ (i.e.,
 $\fp$ is non-degenerate). However, $\, \ker \, \omega^1 \,$ and
 $\, \ker \, \omega^2 \,$ are orthogonal under~$\omega^1$ and under~%
 $\omega^2$ but not under~$\omega^3$. Now if there existed a poly%
 lagrangian subspace it would have to coincide with the sum of the
 kernels of all the projected forms, but that is the whole space
 $\mathbb{R}^3$, which is not isotropic. Thus $\fp$ does not admit
 a polylagrangian subspace.
\end{exem}
Finally, we observe that even if the sum of the kernels of all the projected
forms is an isotropic subspace with respect to $\fp$, it may still fail to be
a polylagrangian subspace, as shown by the following
\begin{exem}~
 ($\, \hat{n} = 2 \,$, $\, N = 2 \,$, $\, \dim V = 5 \,$,
  $\, \ker \, \fp = \{0\} \,$)\,: \\[2mm]
 Let $\, V = \mathbb{R}^5$, $\hat{T} = \mathbb{R}^2 \,$  and consider the
 $\mathbb{R}^2$-valued $2$-form $\fp$ built from the following two ordinary
 $2$-forms:
 \[
  \omega^1~=~dx^1 \,\smwedge\, dx^4 \, + \, dx^2 \,\smwedge\, dx^3~~~,~~~
  \omega^2~=~dx^1 \,\smwedge\, dx^3 \, - \, dx^2 \,\smwedge\, dx^5~.
 \]
 Then for $\; \hat{t}^* = \hat{t}_a^* \, \hat{e}^a \smin\, (\mathbb{R}^2)^*$,
 we have
 \[
  \omega_{\hat{t}^*}^{}~=~\hat{t}_a^* \, \omega^a~
  =~dx^1 \,\smwedge\,
    \bigl( \hat{t}_1^* \, dx^4 \, + \, \hat{t}_2^* \, dx^3 \bigr) \, + \,
    dx^2 \,\smwedge\,
    \bigl( \hat{t}_1^* \, dx^3 \, - \, \hat{t}_2^* \, dx^5 \bigr)~.
 \]
 Obviously, $\omega^1$, $\omega^2$ and the forms
 \[
  \begin{array}{c}
   (\omega^1)^2~\equiv~\omega^1 \,\smwedge\, \omega^1~
                =~2 \; dx^1 \>\smwedge\, dx^2 \>\smwedge\,
                       dx^3 \>\smwedge\, dx^4~, \\[2mm]
   \omega^1 \,\smwedge\, \omega^2~
   =~dx^1 \>\smwedge\, dx^2 \>\smwedge\, dx^4 \>\smwedge\, dx^5~, \\[2mm]
   (\omega^2)^2~\equiv~\omega^2 \,\smwedge\, \omega^2~
                    =~2 \; dx^1 \>\smwedge\, dx^2 \>\smwedge\,
                           dx^3 \>\smwedge\, dx^5~,
  \end{array}
 \]
 are linearly independent and hence $\fp$ has uniform rank~$2$, since there
 exists no non-zero $6$-form on a five-dimensional space. On the other hand,
 we have
 \[
  \ker \, \omega_{\hat{t}^*}^{}~
  =~\langle \, \hat{t}_1^* \hat{t}_2^* \, \frac{\partial}{\partial x^3} \, - \,
               (\hat{t}_2^*)^2 \, \frac{\partial}{\partial x^4} \, + \,
               (\hat{t}_1^*)^2 \, \frac{\partial}{\partial x^5} \, \rangle~.
 \]
 The intersection of the two kernels $\, \ker \, \omega^1 \,$ and~%
 $\, \ker \, \omega^2 \,$ is $\{0\}$ (i.e., $\fp$ is non-degenerate).
 Note that their (direct) sum is the two-dimensional subspace of~$V$,
 say~$L'$, spanned by $\partial/\partial x^4$ and $\partial/\partial
 x^5$, whereas the subspace spanned by all the kernels $\, \ker \,
 \omega_{\hat{t}^*}$ ($\hat{t}^* \smin\, \hat{T}^* \setminus \{0\}$)
 is the three-dimensional subspace of~$V$, say $L''$, spanned by
 $\partial/\partial x^i$ with $\, i = 3,4,5$, and this is isotropic
 with respect to all the forms $\omega_{\hat{t}^*}$ ($\hat{t}^* \smin\,
 \hat{T}^* \setminus \{0\}$). More than that: since its codimension
 is $2$, it is maximal isotropic with respect to all the forms
 $\omega_{\hat{t}^*}$ ($\hat{t}^* \smin\, \hat{T}^* \setminus \{0\}$).
 Now if there existed a polylagrangian subspace it would have to
 coincide with $L'$ and also with $L''$, but these two are not equal
 and do not have the right dimension, which according to equation~%
 (\ref{eq:POLID2}) would have to be $4$: both of them are too small.
 Thus $\fp$ does not admit a polylagrangian subspace.
\end{exem}

\noindent
To summarize, the examples given above show that the hypothesis of existence
of a poly\-lagrangian subspace is highly non-trivial and very restrictive: as
it seems, it cannot be replaced by any other hypothesis that is not obviously
equivalent. The examples also show the great variety of possibilities for
the ``relative positions'' of the kernels of the various projected forms
that prevails when such a subspace does not exist. In this sense, the
definition adopted in Ref.~\cite{Gr} is quite inconvenient, since it
makes no reference to this subspace, thus hiding the central aspect
of the theory.

To conclude, we want to add some remarks about the relation between the
polylagrangian subspace, when it exists, and the more general class of
maximal isotropic subspaces. First, we emphasize that in contrast with
a polylagrangian subspace, maximal isotropic subspaces always exist.
To prove this, it suffices to start out from an arbitrary one-dimensional
subspace $L_1^{}$, which is automatically isotropic, and construct a chain
$\, L_1^{} \smsubset L_2^{} \smsubset \ldots \,$ of subspaces where $L_{p+1}$
is defined as the direct sum of $L_p^{}$ and the one-dimensional subspace
spanned by some non-zero vector in its $1$-orthogonal complement $L_p^{\fp,1}$.
For dimensional reasons, this process must stop at some point, which means
that at this point we have succeeded in constructing a maximal isotropic
subspace. However, nothing guarantees that maximal isotropic subspaces
resulting from different chains must have the same dimension, nor that
there must exist some chain leading to a maximal isotropic subspace
of sufficiently high dimension to be polylagrangian: this happens
only in the special case of ordinary forms ($\hat{n} = 1$), where
all maximal isotropic subspaces have the same dimension and where
the notions of a poly\-lagrangian subspace (or simply lagrangian
subspace, in this case) and of a maximal isotropic subspace coincide.

Another important point concerns the relation between the notions of
isotropic subspace and maximal isotropic subspace with respect to the
form $\,\fp$ and with respect to its projections. First, it is obvious
that a subspace of~$V$ is isotropic with respect to $\,\fp$ if and
only if it is isotropic with respect to each of the projected forms
$\,\omega_{\hat{t}^*}^{}$ ($\hat{t}^* \smin\, T^* \setminus \{0\}$) or
$\,\omega^a$ ($1 \leqslant a \leqslant \hat{n}$). However, this no
longer holds when we substitute the term ``isotropic'' by the term
``maximal isotropic'': a subspace of~$V$ that is maximal isotropic
with respect to each of the projections of~$\,\fp$ is also maximal
isotropic with respect to~$\,\fp$, but conversely, it can very well
be maximal isotropic with respect to $\,\fp$ (and hence isotropic
with respect to each of the projections of~$\,\fp$) but fail to
be maximal isotropic with respect to some of them. And finally,
a polylagrangian subspace of~$V$ is maximal isotropic with
respect to each of the projections of~$\,\fp$ (this follows from
Theorem~\ref{thm:POLILS1}), but as we have seen in the last example
above, the converse is false: a subspace can be maximal isotropic
with respect to each of the projections of~$\,\fp$ without being
polylagrangian. All these statements illustrate the special nature
of the polylagrangian subspace, already in the case of vector-valued
$2$-forms.

\section{Poincar\'e Lemma}

In this appendix we want to state the Poincar\'e lemma in the form in which
it is used in the proof of the Darboux theorem in Section~\ref{sec:DT}.
\begin{thm}~ \label{thm:POINCL}
 Let $\, \omega \smin \Omega^k(P,\hat{T}) \,$ be a closed form on
 a manifold $P$ taking values in a fixed vector space $\hat{T}$ and
 let $L$ be an involutive distribution on~$P$. Suppose that $\,\omega$
 is $(k-r)$-horizontal (with respect to $L$), i.e., such that for any
 $\, p \,\smin P \,$ and all $\, v_1,\dots,v_{r+1} \smin L_p^{}$,
 we have
 \[
  \mathrm{i}_{v_1} \ldots\, \mathrm{i}_{v_{r+1}} \omega_p^{}~=~0~.
 \]
 Then for any point of~$P$ there exist an open neighborhood $U$ of that
 point and a ($k-1$)-form $\, \theta \smin \Omega^{k-1}(U,\hat{T}) \,$
 on~$U$ which is also $(k-r)$-horizontal (with respect to $L$), i.e.,
 such that for any $\, p \,\smin U \,$ and all $\, v_1,\dots,v_r \smin
 L_p^{}$, we have
 \[
  \mathrm{i}_{v_1} \ldots\, \mathrm{i}_{v_r} \theta_p^{}~=~0~,
 \]
 and such that $\, \omega = d\theta \,$ on~$U$.
\end{thm}
\proof~
 Due to the local character of this theorem and since the subbundle~$L$
 of~$TP$ is involutive, we can without loss of generality work in a
 local chart of the manifold~$P$ around the chosen reference point in
 which the foliation defined by~$L$ is ``straightened out'', so we may
 assume that $\, P \cong K_0^{} \oplus L_0^{} \,$ with $\, L \cong L_0^{}$
 with fixed subspaces $K_0^{}$ and $L_0^{}$ and such that the aforementioned
 reference point corresponds to the origin. (In what follows, we shall omit
 the index $0$.) We also suppose that $\, \hat{T} = \mathbb{R}$, since
 we may prove the theorem separately for each component of~$\omega$
 and~$\theta$, with respect to some fixed basis of~$\hat{T}$.
 \\[2mm]
 For $\, t \smin \mathbb{R}$, define the ``$K$-contraction'' $\, F_t^K : P
 \rightarrow P \,$ and the ``$L$-contraction'' $\, F_t^L : P \rightarrow P \,$
 by $\, F_t^K(x,y) = (tx,y) \;$ and $\, F_t^L(x,y)=(x,ty) \,$; obviously,
 $F_t^K$ and $F_t^L$ are diffeomorphisms if $\, t \neq 0 \,$ and are
 projections if $\, t = 0$. Associated with each of these families of
 mappings there is a time dependent vector field which generates it
 in the sense that, for $\, t \neq 0$,
 \[
  X_t^K \bigl( F_t^K(x,y) \bigr)~=~\frac{d}{ds} \, F_s^K(x,y) \, \bigg|_{s=t}
  \qquad \mbox{and} \qquad
  X_t^L \bigl( F_t^L(x,y) \bigr)~=~\frac{d}{ds} \, F_s^L(x,y) \, \bigg|_{s=t}~.
 \]
 Explicitly, for $\, t \neq 0$,
 \[
  X_t^K(x,y)~=~t^{-1}(x,0) \qquad \mbox{and} \qquad X_t^L(x,y)~=~t^{-1}(0,y)~.
 \]
 Define $\; \omega_0^{} = (F_0^L)^* \omega \,$ and, for $\, \epsilon > 0$,
 \[
  \theta_\epsilon~
  =~\int_\epsilon^1 dt~
    \Bigl( (F_t^L)^*(\mathrm{i}_{X_t^L}^{} \omega) \, + \,
                    (F_t^K)^*(\mathrm{i}_{X_t^K}^{} \omega_0^{}) \Bigr)~,
 \]
 as well as
 \[
  \theta~=~\lim_{\epsilon \rightarrow 0} \; \theta_\epsilon~
  =~\int_0^1 dt~\Bigl( (F_t^L)^*(\mathrm{i}_{X_t^L}^{} \omega) \, + \,
                       (F_t^K)^*(\mathrm{i}_{X_t^K}^{} \omega_0^{}) \Bigr)~.
 \]
 To see that the $(k-1)$-forms $\theta$ and $\theta_\epsilon$ are well
 defined, consider $k-1$ vectors $\, (u_i^{},v_i^{}) \smin K \oplus L$
 ($1 \leqslant i \leqslant k-1$) and observe that, for $\, t \neq 0$,
 \[
  (F_t^L)^*(i_{X_t^L}^{} \omega)_{(x,y)}^{}
  \bigl( (u_1^{},v_1^{}) , \ldots , (u_{k-1}^{},v_{k-1}^{}) \bigr)~
  =~\omega_{(x,ty)}^{}
    \bigl( (0,y) , (u_1,tv_1) , \ldots , (u_{k-1},tv_{k-1}) \bigr)~,
 \]
 and
 \[    
  (F_t^K)^*(i_{X_t^K}^{} \omega_0^{})_{(x,y)}^{}
  \bigl( (u_1^{},v_1^{}) , \ldots , (u_{k-1}^{},v_{k-1}^{}) \bigr)~
  =~t^{k-1} \, \omega_{(tx,0)}^{}
    \bigl( (x,0) , (u_1^{},0) , \ldots, (u_{k-1}^{},0) \bigr)~.
 \]
 Here we see easily that both expressions are differentiable in~$t$ and
 provide $(k-1)$-forms which are $(k-r)$-horizontal and $(k-1)$-horizontal
 with respect to~$L$, respectively. Thus, $\theta_\epsilon$ and $\theta$
 are $(k-1)$-forms which are $(k-r)$-horizontal with respect to~$L$.
 Moreover, since $\, d\omega = 0 \,$ and $\, d\omega_0^{} = 0$,
 \begin{eqnarray*}
  d\theta_{\epsilon} \!\!
  &=&\!\! \int_\epsilon^1 dt~
          \Bigl( (F_t^L)^*(\mathrm{i}_{X_t^L}^{} \omega) \, + \,
                 (F_t^K)^*(\mathrm{i}_{X_t^K}^{} \omega_0^{}) \Bigr) \\
  &=&\!\! \int_\epsilon^1 dt~\Bigl( (F_t^L)^*(L_{X_t^L}^{} \omega) \, + \,
                                   (F_t^K)^*(L_{X_t^K}^{} \omega_0^{}) \Bigr) \\
  &=&\!\! \int_\epsilon^1 dt~
          \Bigl( \frac{d}{dt} \bigl( (F_t^L)^* \omega \bigr) \, + \,
                 \frac{d}{dt} \bigl( (F_t^K)^* \omega_0^{} \bigr) \Bigr) \\[3mm]
  &=&\!\! \omega \, - \, (F_\epsilon^L)^* \omega \, + \,
          \omega_0^{} \, - \, (F_\epsilon^K)^* \omega_0^{}~.
 \end{eqnarray*}
 Taking the limit $\, \epsilon \rightarrow 0$, we get $\, (F_\epsilon^L)^*
 \omega \rightarrow \omega_0^{} \,$ and $\, (F_\epsilon^K)^* \omega_0^{}
 \rightarrow 0 \,$ and hence $d\theta = \omega$.
\qed

\end{appendix}


\end{document}